\documentclass[10pt,a4paper]{article}

\usepackage{amstext}
\usepackage{array}
\usepackage{amsfonts}
\usepackage{amssymb, bm}
\usepackage{graphicx}
\usepackage{epstopdf}
\usepackage{algorithm}     
\usepackage{algpseudocode} 
\usepackage{url}
\usepackage{caption}
\usepackage{subfig}
\usepackage{csquotes}
\usepackage{tikz}
\usetikzlibrary{arrows}
\usetikzlibrary{decorations.markings}
\usetikzlibrary{shapes.geometric}
\usetikzlibrary{shapes.arrows}
\usetikzlibrary{fit,calc}
\usetikzlibrary{quotes,angles}
\usetikzlibrary{positioning, arrows.meta}
\usepackage{animate}  
\usepackage{graphicx}
\usepackage{pgf}
\usepackage{bbm}
\usepackage{pgfplots}
\usepackage{mathtools}
\usepackage{physics}
\usepackage{amsthm}
\usepackage{hyperref}
\usepackage{stmaryrd}

\usepackage{diagbox}

\usepackage{amsmath,amssymb}
\usepackage{booktabs}
\usepackage{hyperref}
\usepackage{enumitem}
\usepackage{graphicx}

\usepackage{algorithm}     
\usepackage{algpseudocode} 

\usepackage[utf8]{inputenc} 
\usepackage[T1]{fontenc}

\usepackage{tikz}
\usetikzlibrary{decorations.markings}
\usetikzlibrary{shapes.geometric}
\usetikzlibrary{shapes.arrows}
\usetikzlibrary{fit,calc}
\usetikzlibrary{quotes,angles}
\usetikzlibrary{arrows.meta}
\usepackage{animate}  
\usepackage{graphicx}
\usepackage{pgf}

\usepackage{pgfplots}
\usepackage[normalem]{ulem}

\newcommand{\vertiii}[1]{{\left\vert\kern-0.25ex\left\vert\kern-0.25ex\left\vert #1 
   \right\vert\kern-0.25ex\right\vert\kern-0.25ex\right\vert}}

\usepackage{geometry}
\makeatletter
\def\blfootnote{\xdef\@thefnmark{$\star$}\@footnotetext}
\makeatother
\newenvironment{Authors}%
  {\begin{center}\begin{bfseries}}%
  {\end{bfseries}\end{center}}
\newenvironment{Addresses}%
  {\begin{flushleft}\begin{itshape}}%
  {\end{itshape}\end{flushleft}}
  \newcommand{\email}[1]{\hspace*{\stretch{1}}\emph{\texttt{#1}}}

 \usepackage{fancyhdr}  
  \fancypagestyle{plain}{
\fancyhead{}
\fancyhead[C]{\hfill Submitted to Elsevier, December 2025}
 }

\newcommand{\bs}[1]{\boldsymbol{#1}}

\newtheorem{theorem}{Theorem}[section]

\newtheorem{lemma}[theorem]{Lemma}

\newtheorem{remark}[theorem]{Remark}
\newtheorem*{remark*}{Remark}  

\begin{document}

\title{
High-order implicit Runge-Kutta time integrators for component-based model reduction of FSI problems}

 \date{}

 \maketitle
\vspace{-50pt} 
 
\begin{Authors}
Tommaso Taddei$^{1}$, Xuejun Xu$^{2,3}$
Lei Zhang$^2$.
\end{Authors}

\begin{Addresses}
$^1$
Dipartimento di Matematica, Sapienza University of Rome, 00185 Rome, Italy \\
\email{tommaso.taddei@uniroma1.it} \\
$^2$
School of Mathematical Sciences, Key Laboratory of Intelligent Computing and Applications (Ministry of Education), Tongji University, Shanghai 200092, China \\ \email{xxj@lsec.cc.ac.cn,22210@tongji.edu.cn} \\ 
$^3$
Institute of Computational Mathematics, AMSS, Chinese Academy of Sciences, Beijing 100190, China \\
\email{xxj@lsec.cc.ac.cn} \\[3mm]

\end{Addresses}

\begin{abstract}
 We propose a model order reduction framework for incompressible fluid-structure interaction (FSI) problems based on high-order implicit Runge-Kutta (IRK)  methods. We consider separate reduced spaces for fluid velocity, fluid pressure and solid displacement; we enrich the velocity space with supremizer modes to ensure the inf-sup stability of the fluid subproblem; we consider bubble-port decomposition of fluid velocity and solid displacement to satisfy the kinematic conditions at the fluid structure interface. 
 We resort to Galerkin projection to define the semi-discrete reduced-order model and we consider a Radau-IIA IRK method for time integration: the resulting algebraic system is solved using static condensation of the interface degrees of freedom.
 The reduced-order model  preserves a semi-discrete energy balance inherited from the full-order model, and avoids the need for additional interface enrichment. Numerical experiments demonstrate that the proposed combination of high-order IRK schemes with  bubble-port decoupling of velocity and displacement degrees of freedom yields stable and accurate  reduced-order model for  long-time integration of strongly-coupled  parametric FSI problems.
\end{abstract}


\section{Introduction}

\subsection{Model reduction of fluid-structure interaction problems}

Fluid-structure interaction (FSI) problems arise in a wide range of engineering and scientific applications, such as  the aeroelastic behavior  of wings and turbine blades \cite{piperno1995partitioned,
piperno2001partitioned,wang2016fluid}, wave-induced dynamics of offshore structures \cite{battisti2025multi,agamloh2008application}, hemodynamics in compliant vessels \cite{crosetto2011parallel,sun2021advanced}, or micro-electro-mechanical systems (MEMS) \cite{ostasevicius2007numerical,dutta2024monolithic}. 
Accurate numerical simulations of   coupled problems are essential for understanding, predicting, and controlling the behavior of the underlying physical system. 
However, despite the many recent advances, high-fidelity (HF) computational models of FSI problems are extremely computationally demanding: they are hence ill-suited to deal with 
real-time (e.g., control) and  many-query (e.g., parametric studies, uncertainty quantification) tasks and also with long time integration.

Over the past two decades,
model order reduction (MOR) \cite{Quarteroni2015_ROMbook,Hesthaven_ROM_book2016}  
 has emerged as an effective strategy for  alleviating the computational burden associated with 
 parametric studies of engineering  systems governed by partial differential equations (PDEs). By approximating the HF dynamics within a suitably chosen low-dimensional subspace, reduced-order models (ROMs) significantly lower computational demands while maintaining predictive accuracy. 
In particular, projection-based methods exploit (Petrov-)Galerkin projection to devise model-based ROMs of parametric systems.
Traditional (monolithic) MOR methods have been successfully applied to a broad range of parametric problems, but they are ill-suited to deal with coupled and/or large-scale systems 
with parameter-induced topology changes for which it is not possible to define a single parameter-independent low-dimensional subspace.
To address this issue, several authors have proposed component-based (CB) MOR techniques that combine traditional MOR strategies with   domain decomposition (DD) methods to handle  large-scale systems 
\cite{buhr2020localized,huynh2013static,Pegolotti2021}.

A variety of projection-based ROMs have been proposed for FSI, including monolithic POD-Galerkin approaches \cite{Ballarin2016_fsi}, partitioned semi-implicit schemes \cite{Ballarin2017_fsi,Nonino2023_fsi}, and time-segmented techniques \cite{zhai2025projection}. Non-intrusive and segregated strategies have also been investigated \cite{xiao2016non,ngan2025reduced}. 
In the context of hemodynamics, Deparis and coauthors  developed reduced FSI problems based on the thin-walled assumption and transpiration condition, and they  applied CB-MOR techniques to further reduce the computational costs \cite{colciago2018reduced,deparis2019coupling,Pegolotti2021}.
We refer to \cite{colciago2017fluid} for a thorough overview of full-order and reduced-order methods for vascular FSI simulations;
we also refer to  \cite{riffaud2024low} for a recent application of projection-based MOR techniques to uncertainty quantification of FSI systems in hemodynamics.
Despite these advances, the development of stable and accurate ROMs for FSI remains challenging due to stringent requirements for preserving inter-domain coupling conditions and numerical stability.

A key factor that influences the stability and accuracy of FSI simulations, and thus the effectiveness of associated ROM techniques, is the choice of temporal discretization scheme. In practical FSI computations of incompressible flows, it is common to integrate the incompressible Navier–Stokes equations using the second-order backward differentiation formula (BDF2), while employing the Newmark method for structural dynamics, which results in an overall second-order accuracy in time \cite{chabannes2013high,deparis2016facsi,wood2010partitioned}. However, the standard Newmark method may exhibit numerical instability for geometrically nonlinear structures, where energy transfer occurs across an infinite number of deformation modes \cite{vazquez2007nonlinear,cori2015high}. 
Although alternative dissipative methods, such as Bossak, generalized–$\alpha$, or explicit singly diagonally implicit Runge–Kutta (ESDIRK) methods, can enhance numerical stability by effectively damping spurious high-frequency oscillations \cite{etienne2009geometric,etienne2009perspective}, their inherent numerical dissipation risks to excessively attenuate physically-relevant structural dynamics.

These observations highlight the need for temporal schemes that are both strongly stable and sufficiently accurate, without relying on excessively dissipative modifications. 
Rather than increasing numerical damping, 
a more attractive strategy is to employ higher-order temporal
schemes that can achieve a prescribed accuracy with larger time steps and result in superior computational efficiency, which is particularly beneficial for long-time FSI simulations 
\cite{bijl2002implicit,van2005higher,yang2007higher,cori2015high,averweg2024monolithic}. 
However, linear multistep methods such as BDF are subject to the second Dahlquist barrier \cite{hairer1991solving2,abu2023monolithic}, which prevents higher-order variants from achieving desirable stability properties such as A- and L-stability. On the other hand, 
fully implicit Runge-Kutta (IRK) methods \cite{hairer1991solving2,hairer2006geo} 
naturally overcome these stability limitations.

In this work, we rely on  Radau-IIA IRK methods  for time discretization of both full-order and reduced-order models.
Radau-IIA IRK methods
are high-order, 
A-stable, L-stable, and  stiffly accurate:
they  enable effective suppression of spurious numerical oscillations at high frequencies without compromising the accurate representation of physically relevant low-frequency dynamics \cite{hairer1991solving2,abu2023monolithic,cori2015high}.
Radau-IIA IRK methods 
can be applied to both fluid and structural equations and they
effectively cure 
numerical instabilities associated with nonlinear structural dynamics encountered in Newmark-based approaches, without introducing excessive dissipation
\cite{etienne2009geometric,etienne2009perspective,
hay2014high,cori2015high}.

\subsection{Application of high-order implicit Runge  Kutta time integrators to component-based model reduction}

In this work, we aim to devise energy-stable ROMs for FSI systems: towards this end, we combine the static condensation reduced basis element (scRBE) method with the IRK time integrators discussed above.
The  (port-reduced) scRBE method
exploits a port-bubble solution decomposition to ensure solution continuity at the interfaces (ports), and  relies on static condensation to eliminate internal degrees of freedom and effectively solve the reduced system.
The scRBE method was first introduced for linear coercive elliptic problems in
\cite{huynh2013static} and
\cite{eftang2013port} and then extended to mechanics problems with local nonlinearities in 
\cite{ballani2018component} and to fully-nonlinear problems in 
\cite{ebrahimi2024hyperreduced}.
We recall that the scRBE method shares many features with 
component mode synthesis (CMS) in structural dynamics \cite{craig1968coupling,seshu1997substructuring,castanier2001characteristic}.

 To our knowledge, the present work constitutes the first application of the scRBE method to FSI systems; we note, however, that CMS has already been applied to linearized FSI problems in \cite{magalhaes2005development,yao2023symmetric}.
The combination of projection-based MOR methods with IRK time integration is also new: while the need for solving a much larger and possibly poorly-conditioned nonlinear problem has limited so far the application of IRK methods to HF FSI simulations, the additional cost at the MOR level is  much less significant thanks to port reduction that enables the use of direct linear solvers.

We observe that our formulation can be interpreted as an optimization-based method where the solid displacement at the interface is treated as the control. In this respect, the method is related to previous efforts on optimization-based  coupling of fluid 
\cite{taddei2024non}
and
FSI \cite{taddei2025optimization} problems.
The methods
\cite{taddei2024non,taddei2025optimization}
(see also \cite{prusak2023optimisation,prusak2024optimisation})
exploit the normal flux as the control variable and were inspired by the works by Gunzburger and coauthors on optimization-based coupling of  elliptic problems and   the incompressible Navier-Stokes equations
\cite{Gunzburger_Lee_2000,Gunzburger_Peterson_Kwon_1999} and the subsequent works on FSI  
\cite{Kuberry_Lee_2013,Kuberry_Lee_2015,Kuberry_Lee_2016}.
We show that our method has two key advantages over the optimization-based coupling of \cite{taddei2024non,taddei2025optimization}: first, we prove a semi-discrete energy balance for the reduced system, which clearly illustrates the superior stability properties of the new formulation; second, thanks to the port-bubble decomposition the new formulation does not require any enrichment of the velocity basis to restore the full rank of the reduced Jacobian of the coupled problem.

The remainder of this paper is organized as follows. Section \ref{sec:gov} reviews the governing equations of the FSI problem in ALE form and briefly recalls the weak formulation. Section \ref{sec:discrete} presents the finite element spatial discretization, introduces the Radau-IIA IRK time integrator, and derives the resulting fully discrete coupled system, together with the optimization-based partitioned algorithm. 
Section \ref{sec:rom} presents the construction of bubble-port reduced spaces and derives the associated semi-discrete and fully discrete ROMs. Numerical results are presented in Section \ref{sec:numerics}. Section \ref{sec:conclusions} provides concluding remarks and outlines future research directions.

 \section{Governing equations and semi-discrete formulation}
 \label{sec:gov}
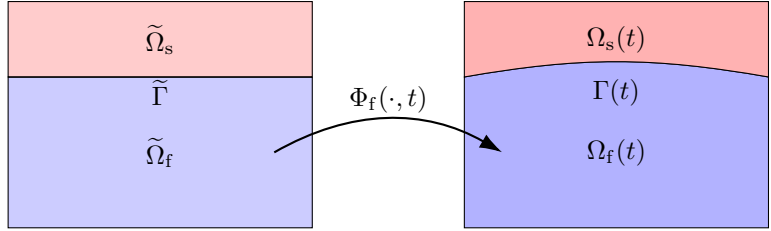
\begin{figure}[H]
    \centering
\begin{tikzpicture}
    \draw [fill=blue!20] (0,0) rectangle (4,2);
    \node at (2,1) {$\widetilde{\Omega}_{\rm f}$};

    \draw [fill=red!20] (0,2) rectangle (4,3);
    \node at (2,2.5) {$\widetilde{\Omega}_{\rm s}$};

    \draw (0,2) -- (4,2);
    \node at (2,1.8) {$\widetilde{\Gamma}$};

    \draw [fill=blue!30] (6,0) -- (10,0) -- (10,2) to [out=170,in=10] (6,2) -- cycle;
    \node at (8,1) {$\Omega_{\rm f}(t)$};

    \draw [fill=red!30] (6,2) to [out=10,in=170] (10,2) -- (10,3) -- (6,3) -- cycle;
    \node at (8,2.5) {$\Omega_{\rm s}(t)$};

    \node at (8,1.8) {$\Gamma(t)$};

    \draw[-{Latex[length=3mm, width=2mm]}, thick] (3.5,1) to [out=30,in=150] (6.5,1);
    	\node at (5.0,1.7) {${\Phi}_{\rm f}(\cdot, t)$};
\end{tikzpicture}

\caption{Reference  and current   domains of FSI problems. Here, $\Phi_{\rm f}(t)$ is the ALE bijection that maps the reference fluid domain $\widetilde{\Omega}_{\rm f}$ into the current domain $\Omega_{\rm f}(t)$.}
\label{fig:FSI_illustration}
\end{figure}

 We consider the interaction between an  incompressible viscous fluid and  an elastic solid in an Arbitrary Lagrangian Eulerian (ALE) setting.  
Let $\widetilde{\Omega}_{\mathrm f}$, $\widetilde{\Omega}_{\mathrm s}\subset\mathbb R^{\texttt{d}}$ ($\texttt{d}\in\{2,3\}$) denote the reference fluid and solid domains, respectively, that share the interface $\widetilde{\Gamma}=\partial\widetilde{\Omega}_{\mathrm f}\cap\partial\widetilde{\Omega}_{\mathrm s}$, as illustrated in Figure \ref{fig:FSI_illustration}.  The current fluid and solid configurations at time $t\ge0$ are obtained through the mappings
\begin{equation*}
  \Omega_{\mathrm f}(t)=\Phi_{\mathrm f}\bigl(\widetilde{\Omega}_{\mathrm f},t\bigr),\qquad
  \Omega_{\mathrm s}(t)=\bigl\{\widetilde{\mathbf x}+\widetilde{\mathbf d}_{\mathrm s}(\widetilde{\mathbf x},t):\widetilde{\mathbf x}\in\widetilde{\Omega}_{\mathrm s}\bigr\},
\end{equation*}
where $\widetilde{\mathbf d}_{\mathrm s}$ is the solid displacement, and the current interface  
$  \Gamma(t)=\Phi_{\mathrm f}\bigl(\widetilde{\Gamma},t\bigr)$.
The ALE velocity is $\widetilde{\boldsymbol\omega}_{\mathrm f}=\partial_t\Phi_{\mathrm f}$, and the ALE time derivative of any Eulerian field $q$ is
\begin{equation}
  \partial_t q\bigl|_{\Phi_{\mathrm f}}(\mathbf x,t):=\frac{d}{dt}\Bigl(q\bigl(\Phi_{\mathrm f}(\widetilde{\mathbf x},t),t\bigr)\Bigr),\quad \mathbf x=\Phi_{\mathrm f}(\widetilde{\mathbf x},t).
\end{equation}

We briefly recall the strong and weak formulations of the coupled FSI problem \cite{formaggia2010cardiovascular,
Barker2010,Wu_Cai_2014,taddei2025optimization}. 
\subsection{Strong form}
\paragraph{Fluid subproblem} 
The fluid is governed by 
the incompressible Navier–Stokes in the ALE framework
\begin{equation}
  \label{eq:fluid_model}
\left\{
\begin{array}{ll}
\displaystyle{
\rho_{\mathrm f}\partial_t\mathbf u_{\mathrm f}\bigl|_{\Phi_{\mathrm f}} 
+ \rho_{\rm f}\bigl[(\mathbf{u}_{\rm f}-\bs{\omega}_{\rm f})\cdot \nabla\bigr]\mathbf{u}_{\rm f}  -\nabla\cdot \bs{\sigma}_{\rm f} =0, }&  \mbox{in}\; \Omega_{\rm f}(t), \\[2mm]
\nabla \cdot \mathbf{u}_{\rm f} = 0, & \mbox{in}\; \Omega_{\rm f}(t), \\[2mm]
\mathbf{u}_{\rm f}|_{\Gamma^{\rm dir}_{\rm f}}=\mathbf{u}_{\rm f}^{\rm dir},
\;\;
\bs{\sigma}_{\rm f} \mathbf{n}_{\rm f}
|_{\Gamma^{\rm neu}_{\rm f}}=\mathbf{g}^{\rm neu}_{\rm f},
&
\end{array}
\right.
\end{equation} 
where $\mathbf u_{\mathrm f}$ denotes the fluid velocity, $\rho_{\rm f}$ the density, 
$$\boldsymbol\sigma_{\mathrm f}=2\mu_{\mathrm f}\boldsymbol\varepsilon_{\mathrm f}-p_{\mathrm f}\mathbf I,\quad (\boldsymbol\varepsilon_{\mathrm f}=\tfrac12(\nabla\mathbf u_{\mathrm f}+\nabla\mathbf u_{\mathrm f}^{\mathsf T}))$$ 
is the  Newtonian Cauchy stress tensor, 
$p_{\rm f}$ is the pressure,  and $\mu_{\mathrm f}$ is the dynamic viscosity.
The ALE velocity
$\boldsymbol\omega_{\mathrm f}$ enters in the advection term of the momentum equation to
 account for the fluid domain motion induced by the structural displacement.

\paragraph{Solid subproblem}  The solid is modeled by the elastodynamic equations for the displacement in the Lagrangian framework
\begin{equation}
\left\{
\begin{array}{ll}
\displaystyle{
\widetilde{\rho}_{\mathrm s}\partial_{tt}\widetilde{\mathbf d}_{\mathrm s}-\widetilde{\nabla}\cdot\widetilde{\boldsymbol\Pi}_{\mathrm s}=\mathbf0} &\quad\text{in }\widetilde{\Omega}_{\mathrm s}, \\[3mm]
\widetilde{\mathbf{d}}_{\rm s}|_{\widetilde{\Gamma}^{\rm dir}_{\rm s}}=\widetilde{\mathbf{d}}_{\rm s}^{\rm dir},
\;\;
\widetilde{\bs{\Pi}}_{\rm s} \widetilde{\mathbf{n}}_{\rm s}
|_{\widetilde{\Gamma}^{\rm neu}_{\rm s}}=\widetilde{\mathbf{g}}_{\rm s}^{\rm neu}. 
\end{array}
\right.
  \label{eq:solid_model}
\end{equation}
Here $\widetilde{\rho}_{\mathrm s}$ denotes the solid density, $\widetilde{\boldsymbol\Pi}_{\mathrm s}$ is the first Piola-Kirchhoff stress tensor. 
We consider both the small strain linear elasticity and the Saint-Venant Kirchhoff (SVK) hyperelasticity \cite[pp.98-100]{formaggia2010cardiovascular}\cite{chabannes2013high,Shamanskiy2021_mesh}. 
In the linear elasticity case, 
the first Piola-Kirchhoff stress tensor is given by
 $$\widetilde{\boldsymbol\Pi}_{\mathrm s}=2\mu_{\mathrm s}\widetilde{\boldsymbol\varepsilon}_{\mathrm s}+\lambda_{\mathrm s}(\widetilde\nabla\cdot\widetilde{\mathbf d}_{\mathrm s})\mathbf I,$$ 
where $\lambda_{\mathrm s}$, $\mu_{\mathrm s}$ are the  first and the second Lam\'{e} coefficients, 
 and $\widetilde{\bs\varepsilon}_{\rm s}= \frac{1}{2} (\widetilde\nabla \widetilde{\mathbf{d}}_{\rm s} + \widetilde\nabla \widetilde{\mathbf{d}}_{\rm s}^\top)$ is the strain tensor in the Lagrangian frame. In the SVK case, the first Piola-Kirchhoff stress tensor 
 is given by
 $$\widetilde{\boldsymbol\Pi}_{\mathrm s}=\widetilde{\mathbf F}_{\mathrm s}\,\widetilde{\boldsymbol\Sigma}_{\mathrm s},$$ where $\widetilde{\mathbf F}_{\mathrm s}=\mathbf I+\widetilde\nabla\widetilde{\mathbf d}_{\mathrm s}$ is
 the deformation gradient, $\widetilde{\bs\Sigma}_{\rm s} =  2\mu_{\rm s} \widetilde{\mathbf{E}}_{\rm s} + \lambda_{\rm s} {\rm tr}(\widetilde{\mathbf{E}}_{\rm s}) \mathbf{I}$ 
 is the second Piola-Kirchhoff stress tensor,  
 and  $\widetilde{\mathbf{E}}_{\rm s} = \frac{1}{2}(\widetilde{\mathbf{F}}_{\rm s}^\top  \widetilde{\mathbf{F}}_{\rm s}  - \mathbf{I})$ is the Green-Lagrange tensor. 
 
 \paragraph{Coupling conditions} On the interface $\widetilde\Gamma$ we impose: \\
 (i) geometric compatibility $\Phi_{\mathrm f}=\mathrm{id}+\widetilde{\mathbf d}_{\mathrm s}$, where
$\texttt{id}$ is the identity map, \\
(ii) kinematic continuity $\mathbf u_{\mathrm f}\circ\Phi_{\mathrm f}=\partial_t\widetilde{\mathbf d}_{\mathrm s}$,  and 
\\ (iii) dynamic equilibrium $\widetilde{\boldsymbol\Pi}_{\mathrm s}\widetilde{\mathbf n}_{\mathrm s}+J_{\Phi}\bigl(\boldsymbol\sigma_{\mathrm f}\circ\Phi_{\mathrm f}\bigr)\widetilde\nabla\Phi_{\mathrm f}^{-\mathsf T}\widetilde{\mathbf n}_{\mathrm f}=\mathbf0$, where  
$J_{\Phi}  = {\rm det} (\widetilde{\nabla} \Phi_{\rm f})$.

\paragraph{ALE mapping}
To close the system, we need to determine the ALE mapping. 
Following \cite{Wick2011,Shamanskiy2021_mesh},
we define the ALE mapping
$\Phi_{\rm f}$ using a pseudo-elasticity model.
Specifically, we define $\Phi_{\rm f} = \texttt{id}+\widetilde{\mathbf{d}}_{\rm f}$, where the pseudo-displacement $\widetilde{\mathbf{d}}_{\rm f}$ solves a steady linear elasticity problem, with suitably-chosen pseudo material parameters. As for boundary conditions, we impose $\widetilde{\mathbf{d}}_{\rm f}=0$  on $\partial\widetilde{\Omega}_{\mathrm f}\setminus\widetilde\Gamma$, and 
$\widetilde{\mathbf{d}}_{\rm f} =  \widetilde{\mathbf{d}}_{\rm s}$  on  $\widetilde{\Gamma}$ to 
 match the solid displacement on the fluid-solid interface. 

\subsection{Weak formulation}
 We introduce the function spaces for the fluid:
\begin{equation}
\label{eq:fluid_spaces_continuous}
V_{\rm f}^{\rm dir}(t):=\{\mathbf{v}\in [H^1(\Omega_{\rm f}(t)
)]^{\texttt{d}}  \, : \, \mathbf{v}\big|_{\Gamma_{\rm f}^{\rm dir}(t)} = \mathbf{u}_{\rm f}^{\rm dir}(t)  \},
\quad
V_{{\rm f},0}(t):= [H_{0,\Gamma_{\rm f}^{\rm dir}(t)}^1(\Omega_{\rm f}(t)
)]^{\texttt{d}} ,
\quad
Q_{\rm f}(t):= L^2  (\Omega_{\rm f}(t)
),
\end{equation}
and for the solid
\begin{equation}
\label{eq:solid_spaces_continuous}
V_{\rm s}^{\rm dir}(t):=\{\widetilde{\mathbf{v}}\in [H^1(\widetilde{\Omega}_{\rm s})]^{\texttt{d}}  \, : \, \widetilde{\mathbf{v}}\big|_{\widetilde{\Gamma}_{\rm s}^{\rm dir}} = \widetilde{\mathbf{d}}_{\rm s}^{\rm dir}(t) \},
\quad
V_{{\rm s},0} := [H_{0,\widetilde{\Gamma}_{\rm s}^{\rm dir}(t)}^1(\widetilde{\Omega}_{\rm s})]^{\texttt{d}}.
\end{equation}
The weak formulation reads: 
For any $t>0$, we seek $\mathbf{u}_{\rm f} \in V_{\rm f}^{\rm dir}$, $p_{\rm f} \in Q_{\rm f}$,  
$\widetilde{\mathbf{d}}_{\rm s} \in V_{\rm s}^{\rm dir}$ such that
\begin{subequations}
\label{eq:weak_formulation}
\begin{equation}
\label{eq:weak_formulation_a}
\left\{
\begin{array}{l}
\displaystyle{
\frac{d}{dt}\int_{ \Omega_{\rm f}(t)   }
\rho_{\rm f} 
 \mathbf{u}_{\rm f} 
 \cdot \mathbf{v} \, dx
+
R_{\rm f}'(\mathbf{u}_{\rm f}, p_{\rm f}, \bs\omega_{\rm f}, \mathbf{v}) 
+
\int_{ \widetilde{\Omega}_{\rm s}   }
\widetilde\rho_{\rm s} 
\frac{\partial^2 \widetilde{\mathbf{d}}_{\rm s} }{\partial t^2}
\cdot \widetilde{\mathbf{w}} \, dx
+
R_{\rm s}'(\widetilde{\mathbf{d}}_{\rm s},  \widetilde{\mathbf{w}})
 = 0
};
\\[3mm]
\displaystyle{
b_{\rm f}(\mathbf{u}_{\rm f}, q)
 = 0
};
\\[3mm]
\displaystyle{
\Phi_{\rm f}  \big|_{\widetilde{\Gamma}} = \texttt{id} + \widetilde{\mathbf{d}}_{\rm s}\big|_{\widetilde{\Gamma}},
\quad
\mathbf{u}_{\rm f} \circ \Phi_{\rm f} \big|_{\widetilde{\Gamma}}
=
\frac{\partial \widetilde{\mathbf{d}}_{\rm s}}{\partial t}
\big|_{\widetilde{\Gamma}}
},
  \\
\end{array}
\right.
\end{equation}
for any $\mathbf{v}\in V_{{\rm f},0}$, $ \widetilde{\mathbf{w}}\in V_{{\rm s},0}$ such that $\mathbf{v}\circ \Phi_{\rm f} = \widetilde{\mathbf{w}}$ on $\widetilde{\Gamma}$ and any $ q\in Q_{{\rm f}}$.  
 We emphasize that the test functions for the fluid and the solid velocities are chosen to match at the interface, 
which is consistent with the kinematic continuity condition in the strong formulation. 
Here, $R_{\mathrm f}'$, $R_{\mathrm s}'$ and $b_{\mathrm f}$ denote the spatial parts of the fluid and solid residuals and the fluid  incompressibility constraint, respectively, and are defined by
\begin{equation}
\label{eq:weak_formulation_b}
\left\{
\begin{array}{l}
\displaystyle{
R_{\rm f}'(\mathbf{u}, p, \bs\omega, \mathbf{v})
=
\int_{\Omega_{\rm f}} 
\left(
\bs\sigma_{\rm f}(\mathbf{u},p): \nabla \mathbf{v}  \, + \,
\rho_{\rm f} (\mathbf{u} - \bs\omega) \cdot \nabla \mathbf{u} \cdot \mathbf{v} 
\, - \,
\rho_{\rm f} (\nabla \cdot  \bs\omega)  \mathbf{u} \cdot \mathbf{v}
+
\frac{\rho_{\rm f}}{2}
(\nabla \cdot \mathbf{u}) \mathbf{u}\cdot \mathbf{v} 
\right) \, dx
-
\int_{\Gamma_{\rm f}^{\rm neu}} \mathbf{g}_{\rm f}^{\rm neu} \cdot \mathbf{v} \, dx;
} \\[3mm]
\displaystyle{
b_{\rm f}(\mathbf{u}, q)
= - 
\int_{\Omega_{\rm f}}  
(\nabla \cdot \mathbf{u})  \, q  \, dx
;} \\[3mm]
\displaystyle{
R_{\rm s}'(\widetilde{\mathbf{d}}, \widetilde{\mathbf{w}})
=
\int_{\widetilde{\Omega}_{\rm s}} 
\left(
\widetilde{\boldsymbol{\Pi}}_{\rm s}(\widetilde{\mathbf{d}}): \widetilde\nabla \widetilde{\mathbf{w}} 
\right) \, dx
-
\int_{\widetilde{\Gamma}_{\rm s}^{\rm neu}} \widetilde{\mathbf{g}}_{\rm s}^{\rm neu} \cdot \widetilde{\mathbf{w}} \, dx. 
} \\
\end{array}
\right.    
\end{equation}
\end{subequations}
We use a prime to indicate that the residuals $R_{\rm f}'$ and $R_{\rm s}'$  contain only spatial contributions, that is, no time derivatives. 
  We shall use $R_{\rm f}$ and $R_{\rm s}$ to refer to  the corresponding full residuals that include both temporal and spatial terms. 
  We add  the strongly consistent term 
$\frac{\rho_{\rm f}}{2}
(\nabla \cdot \mathbf{u}) \mathbf{u}\cdot \mathbf{v} $  to the fluid residual  in \eqref{eq:weak_formulation_b}$_1$ to ensure the energy balance at the semi-discrete level  
\cite{nobile2001numerical,taddei2025optimization}.

\subsection{Semi-discrete formulation}
\label{sec:semi_discrete_FSI}
We now introduce a finite element semi-discretization of the weak FSI problem. 
We adopt a nodal $\mathbb{P}_{\kappa}-\mathbb{P}_{\kappa-1}$ Taylor-Hood finite element (FE) discretization for the fluid and a $\mathbb{P}_{\kappa}$ FE discretization for the solid.
Let $ \{ \boldsymbol\phi^{\mathrm f}_i \}_{i=1}^{N_{\mathrm u}}$,
$ \{ \theta^{\mathrm f}_i \}_{i=1}^{N_{\mathrm p}}$ and
$ \{ \boldsymbol\varphi^{\mathrm s}_i \}_{i=1}^{N_{\mathrm s}}$
denote the nodal basis functions for the fluid velocity, fluid pressure and solid
displacement, respectively. We introduce the vectors of finite element coefficients
${\mathbf{u}}_{\rm f}: (0,T] \to \mathbb{R}^{N_{\rm u}},\quad
{\mathbf{p}}_{\rm f}: (0,T] \to \mathbb{R}^{N_{\rm p}}, \quad
{\mathbf{d}}_{\rm s}: (0,T] \to \mathbb{R}^{N_{\rm s}}$, 
and define the associated finite element approximations by
\begin{equation}
\begin{array}{ll}
\displaystyle{
{\mathbf{u}}_{\rm f,h}({\mathbf{x}},t)
=\sum_{j=1}^{N_{\rm u}}
\left( {\mathbf{u}}_{\rm f}(t) \right)_j \bs\phi_j^{\rm f} ({\mathbf{x}},t),
\qquad
{p}_{\rm f,h}({\mathbf{x}},t)
=\sum_{j=1}^{N_{\rm p}}
\left( {\mathbf{p}}_{\rm f}(t) \right)_j \theta_j^{\rm f} ({\mathbf{x}},t),
\qquad}
&
{\mathbf{x}}\in {\Omega}_{\rm f}(t); \\
\displaystyle{
\widetilde{\mathbf{d}}_{\rm s,h}(\widetilde{\mathbf{x}},t)
=\sum_{j=1}^{N_{\rm s}}
\left(  {\mathbf{d}}_{\rm s}(t) \right)_j \bs\varphi_j^{\rm s} (\widetilde{\mathbf{x}}),}
\qquad
&
\widetilde{\mathbf{x}} \in \widetilde{\Omega}_{\rm s}.
\end{array}
\end{equation}
Next, we define the fluid and solid mass matrices $\mathbf M_{\mathrm f}\in \mathbb{R}^{N_{\mathrm u}\times N_{\mathrm u}}$ and $\mathbf M_{\mathrm s} \in \mathbb{R}^{N_{\mathrm s}\times N_{\mathrm s}}$ as
\begin{align*}
(\mathbf M_{\mathrm f})_{ij}(t)
&= \int_{\Omega_{\mathrm f}(t)}
\rho_{\mathrm f}
\boldsymbol\phi^{\mathrm f}_j(\mathbf x,t)\cdot \boldsymbol\phi^{\mathrm f}_i(\mathbf x,t)\, dx,
&& i,j=1,\dots,N_{\mathrm u}, \\
(\mathbf M_{\mathrm s})_{ij}
&= \int_{\widetilde{\Omega}_{\mathrm s}}
\widetilde\rho_{\mathrm s}
\boldsymbol\varphi^{\mathrm s}_j(\widetilde{\mathbf x})\cdot \boldsymbol\varphi^{\mathrm s}_i(\widetilde{\mathbf x})\, dx,
&& i,j=1,\dots,N_{\mathrm s},
\end{align*}
and the discrete incompressibility constraint matrix $\mathbf B_{\mathrm f}\in \mathbb{R}^{N_{\mathrm p}\times N_{\mathrm u}}$ 
\begin{equation*}
(\mathbf B_{\mathrm f})_{ij}
= -\int_{\Omega_{\mathrm f}(t)}
\bigl(\nabla\cdot \boldsymbol\phi^{\mathrm f}_j(\mathbf x,t)\bigr)\,
\theta^{\mathrm f}_i(\mathbf x,t)\, dx,
\qquad i=1,\dots,N_{\mathrm p},\; j=1,\dots,N_{\mathrm u}.
\label{eq:Bf_def}
\end{equation*}
We notice that the matrices $\mathbf{B}_{\rm f}$ and 
$\mathbf{M}_{\rm f}$ depend on the interface displacement field
$\widetilde{\mathbf{d}}_{\rm s,h} \big|_{\tilde{\Gamma}}$ through the morphing.
We also define Boolean mask matrices
$\mathbf P_{\bullet,\Gamma},\mathbf P_{\bullet,\rm{in}},\mathbf P_{\bullet,\rm{dir}}$
associated with the interface, the internal and the Dirichlet degrees of freedom, respectively, for the fluid ($\bullet={\rm f}$) and the solid ($\bullet={\rm s}$) subdomains.  For example, $
  \mathbf P_{\mathrm f,\Gamma} \mathbf u_{\mathrm f} \in \mathbb R^{N_{\mathrm c}}
$
collects the fluid velocity DOFs at the fluid-structure interface $\Gamma(t)$, where $N_{\mathrm c}$ is the number of DOFs at the interface. 
Moreover, we introduce the fluid residual vector $\mathbf{R}_\mathrm{f}^{\prime}(\mathbf{u}_\mathrm{f},\mathbf{p}_\mathrm{f},
\mathbf{P}_{\rm s, \Gamma}   \mathbf{d}_\mathrm{s})\in \mathbb{R}^{N_{\rm u}}$, which collects all spatial contributions of the semi-discrete fluid momentum equation (diffusion, pressure, convection, etc.)
\begin{equation}
\begin{array}{ll}
\displaystyle{
\left(\mathbf{R}_\mathrm{f}^{\prime}(\mathbf{u}_\mathrm{f},\mathbf{p}_\mathrm{f},\
\mathbf{P}_{\rm s, \Gamma}   \mathbf{d}_\mathrm{s}
)\right)_i
}
&
\displaystyle{
=
\int_{\Omega_\mathrm{f}(t)}\left(2\mu_\mathrm{f}\boldsymbol{\varepsilon}_\mathrm{f}(\mathbf{u}_\mathrm{f,h}):\boldsymbol{\varepsilon}_\mathrm{f}(\boldsymbol{\phi}_i^\mathrm{f})-p_\mathrm{f,h}\nabla\cdot\boldsymbol{\phi}_i^\mathrm{f}+\rho_\mathrm{f}\left[\left(\mathbf{u}_\mathrm{f,h}-\boldsymbol{\omega}_\mathrm{f,h}\right)\cdot\nabla\right]\mathbf{u}_\mathrm{f,h}\cdot\boldsymbol{\phi}_i^\mathrm{f}
\right.
}\\
& 
\displaystyle{
- \rho_\mathrm{f}\left(\nabla \cdot \boldsymbol{\omega}
_\mathrm{f,h}\right)
\mathbf{u}_\mathrm{f,h}\cdot\boldsymbol{\phi}_i^\mathrm{f}
+\frac{\rho_\mathrm{f}}{2}(\nabla\cdot\mathbf{u}_\mathrm{f,h})\left.\mathbf{u}_\mathrm{f,h}\cdot\boldsymbol{\phi}_i^\mathrm{f}\right)\,d{x}-\int_{\Gamma_\mathrm{f}^\mathrm{neu}(t)}\mathbf{g}_\mathrm{f}^\mathrm{neu}\cdot\boldsymbol{\phi}_i^\mathrm{f}\,dx,
}
\\
&
\qquad  \qquad \qquad
\qquad
\hfill i=1,\dots,N_{\mathrm u}. 
\end{array}
\label{eq:Rf_prime_discrete}
\end{equation}
Similarly, the solid residual vector $\mathbf{R}_\mathrm{s}^{\prime}(\mathbf{d}_\mathrm{s})\in \mathbb{R}^{N_{\rm s}}$ is defined as
\begin{equation}
\left(\mathbf{R}_{\mathrm{s}}^{\prime}(\mathbf{d}_{\mathrm{s}})\right)_{i}=\int_{\widetilde{\Omega}_{\mathrm{s}}}\widetilde{\boldsymbol{\Pi}}_{\mathrm{s}}(\widetilde{\mathbf{d}}_{\mathrm{s,h}}):\widetilde{\nabla}\boldsymbol{\varphi}_{i}^{\mathrm{s}}\,d{x}
-\int_{\widetilde{\Gamma}_{\mathrm{s}}^{\mathrm{neu}}}\widetilde{\mathbf{g}}_{\mathrm{s}}^{\mathrm{neu}}\cdot\boldsymbol{\varphi}_{i}^{\mathrm{s}}\,d{x},\quad
i=1,\ldots,N_{\mathrm{s}}. 
\end{equation}

The semi-discrete FSI equations in algebraic form
can now be expressed as
\begin{equation}
\label{eq:FSI_semidiscrete}
\left\{
\begin{array}{l}
\displaystyle{
\mathbf{P}_{\rm f,\Gamma} \left[
 \frac{d (\mathbf{M}_{\rm f}  \mathbf{u}_{\rm f})}{dt}
+ 
 \mathbf{R}_{\rm f}' (  \mathbf{u}_{\rm f}, \mathbf{p}_{\rm f}, 
  \mathbf{P}_{\rm s,\Gamma}
 \mathbf{d}_{\rm s} ) \right] 
 + 
 \mathbf{P}_{\rm s,\Gamma}
\left[
 \mathbf{M}_{\rm s}  \ddot{\mathbf{d}}_{\rm s}
+ 
 \mathbf{R}_{\rm s}' ({\mathbf{d}}_{\rm s})
 \right]
= 0;
}
\\[3mm]
\displaystyle{
\mathbf{P}_{\rm f,\Gamma}    {\mathbf{u}}_{\rm f}
=
 \mathbf{P}_{\rm s,\Gamma}    \dot{\mathbf{d}}_{\rm s};

}
\\[3mm]
\displaystyle{
\mathbf{P}_{\rm f,in} \left[
 \frac{d (\mathbf{M}_{\rm f}  \mathbf{u}_{\rm f})}{dt}
+ 
 \mathbf{R}_{\rm f}' (  
 \mathbf{u}_{\rm f}, 
 \mathbf{p}_{\rm f}, 
 \mathbf{P}_{\rm s,\Gamma}
 \mathbf{d}_{\rm s} ) \right]
= 0;
}
\\[3mm]
\displaystyle{
 \mathbf{B}_{\rm f}  \mathbf{u}_{\rm f}  
=
0;
}
\\[3mm]
\displaystyle{
\mathbf{P}_{\rm s,in}
\left[
 \mathbf{M}_{\rm s}  \ddot{\mathbf{d}}_{\rm s}
+ 
 \mathbf{R}_{\rm s}' ({\mathbf{d}}_{\rm s})
 \right]
 = 0;
}
\\[3mm]
\displaystyle{
  \mathbf{P}_{\rm f,dir}    {\mathbf{u}}_{\rm f} = \mathbf{u}_{\rm f,dir},
  \;\;
  \mathbf{P}_{\rm s,dir}    {\mathbf{d}}_{\rm s} = \mathbf{d}_{\rm s,dir}.
}
\\
\end{array}
\right.
\end{equation}
Equations \eqref{eq:FSI_semidiscrete}$_1$-\eqref{eq:FSI_semidiscrete}$_2$ correspond to the weak form
of the dynamic equilibrium and the strong form of the kinematic continuity of velocities at the interface.
Equations \eqref{eq:FSI_semidiscrete}$_3$-\eqref{eq:FSI_semidiscrete}$_5$
represent the semi-discrete fluid and solid balance equations within their respective interior domains,
and \eqref{eq:FSI_semidiscrete}$_6$ enforces the Dirichlet boundary conditions.

We now rewrite the coupled semi-discrete FSI system \eqref{eq:FSI_semidiscrete} into an abstract differential-algebraic form, which will allow us to directly apply the implicit Runge-Kutta scheme introduced later in Section \ref{sec:full_discrete}. To this end, we first introduce the structural velocity
\[
  \mathbf u_{\mathrm s}(t) := \dot{\mathbf d}_{\mathrm s}(t),
\]
so that the second-order structural equation is transformed into an equivalent  first-order form. We define the masked momentum variables by
\begin{equation*}
  \mathbf y_\Gamma := \mathbf P_{\mathrm f,\Gamma}\,\mathbf M_{\mathrm f}\,\mathbf u_{\mathrm f} + \mathbf P_{\mathrm s,\Gamma}\,\mathbf M_{\mathrm s}\,{\mathbf u}_{\mathrm s},\qquad
  \mathbf y_{\rm{f,in}} := \mathbf P_{\rm{f,in}}\,\mathbf M_{\mathrm f}\,\mathbf u_{\mathrm f},\qquad
  \mathbf y_{\rm{s,in}} := \mathbf P_{\rm{s,in}}\,\mathbf M_{\mathrm s}\,{\mathbf u}_{\mathrm s},
\end{equation*}
and collect them, along with the structural displacement, into the differential state vector
\begin{equation*}
  \mathbf y(t)
  := \operatorname{col}\bigl(\mathbf y_\Gamma(t),\,\mathbf y_{\rm{f,in}}(t),\,
\mathbf y_{\rm{s,in}}(t),\,\mathbf d_{\mathrm s}(t)\bigr).
\end{equation*}
The fluid pressure degrees of freedom are treated as algebraic variables and collected into
\[
  \mathbf z(t) := \mathbf p_{\mathrm f}(t).
\]

With these definitions, the coupled semi-discrete FSI system \eqref{eq:FSI_semidiscrete} can be compactly expressed as the differential-algebraic equations
\begin{equation}
\left\{
\begin{aligned}
  \dot{\mathbf y}(t) &= \mathbf f\bigl(t,\mathbf y(t),\mathbf z(t)\bigr), \\[2pt]
  \mathbf 0 &= \mathbf g\bigl(t,\mathbf y(t)\bigr),
\end{aligned}
\right.
\label{eq:FSI_DAE_abstract}
\end{equation}
where the right-hand side of the differential equation is given by
\begin{align*}
\mathbf{f}_1 &= - \mathbf P_{\mathrm f,\Gamma}\,\mathbf R'_{\mathrm f}\bigl(\mathbf u_{\mathrm f}(t,\mathbf y),\mathbf z,
 \mathbf{P}_{\rm s,\Gamma}
\mathbf d_{\mathrm s}\bigr)
    - \mathbf P_{\mathrm s,\Gamma}\,\mathbf R'_{\mathrm s}(\mathbf d_{\mathrm s}), \\[2pt]
\mathbf{f}_2 &= - \mathbf P_{\rm{f,in}}\,\mathbf R'_{\mathrm f}\bigl(\mathbf u_{\mathrm f}(t,\mathbf y),\mathbf z,
 \mathbf{P}_{\rm s,\Gamma}
\mathbf d_{\mathrm s}\bigr), \\[2pt]
\mathbf{f}_3 &= - \mathbf P_{\rm{s,in}}\,\mathbf R'_{\mathrm s}(\mathbf d_{\mathrm s}), \\[2pt]
\mathbf{f}_4 &= \mathbf u_{\mathrm s}(t,\mathbf y),
\end{align*}
and the algebraic constraints $\mathbf g(t,\mathbf y) = \mathbf 0$ contains the remaining equations \eqref{eq:FSI_semidiscrete}$_2$, \eqref{eq:FSI_semidiscrete}$_4$, and \eqref{eq:FSI_semidiscrete}$_6$, which do not involve time derivatives.

\begin{remark}
The Boolean mask matrices $(\mathbf P_{\rm f,\Gamma},\mathbf P_{\mathrm{f,in}},\mathbf P_{\mathrm{f,dir}})$ and
$(\mathbf P_{\rm s,\Gamma},\mathbf P_{\rm{s,in}},\mathbf P_{\mathrm{s,dir}})$ define a partition of the fluid and solid
degrees of freedom. For fixed Dirichlet boundary data
$\mathbf u_{\mathrm{f,dir}}(t)$ and $\mathbf d_{\mathrm{s,dir}}(t)$, the masked momentum variables
$(\mathbf y_\Gamma,\mathbf y_{\mathrm{f,in}},\mathbf y_{\mathrm{s,in}})$ uniquely determine the free fluid and structural
velocities. Therefore, we  regard $\mathbf u_{\rm f}$ and ${\mathbf u}_{\rm s}$ as functions of the differential
state $\mathbf y(t)$ and of time, and we write
\begin{equation*}
  \mathbf u_{\rm f} = \mathbf u_{\rm f}(t,\mathbf y),\qquad {\mathbf u}_{\rm s} = \mathbf u_{\rm s}(t,\mathbf y)
\end{equation*}
for brevity.
\end{remark}

\section{Implicit Runge-Kutta time integration}
\label{sec:discrete}

\subsection{Fully-discrete formulation} 
\label{sec:full_discrete} 

This section introduces the fully discrete formulation of the FSI problem based on an implicit Runge-Kutta (IRK) temporal discretization. We first recall the general IRK framework and introduce the relevant notation. 
We then apply the IRK scheme to the coupled differential-algebraic system introduced in Section \ref{sec:semi_discrete_FSI}. Finally, we assemble the fully discrete coupled system into a compact monolithic form suitable for numerical solution. 

Consider a general ordinary differential equation (ODE) of the form $\frac{d\mathbf y}{dt}=\mathbf f(t,\mathbf y)$. An $s$-stage IRK discretization is characterized by the Butcher tableau $(\mathbf A, \mathbf b, \mathbf c)$ 
\begin{equation}
\begin{array} {c|ccccc}
    c_1 & a_{11} & a_{12} & \cdots & a_{1s} \\
    c_2 & a_{21} & a_{22} & \cdots & a_{2s} \\
    \vdots & \vdots & & \ddots & \vdots \\
    c_s & a_{s1} & a_{s2} & \cdots & a_{ss} \\
    \hline
    & b_1 & b_2 & \cdots & b_s
\end{array}.
\end{equation}
Given the state $\mathbf{y}^{(n)}$, we define $\mathbf y^{(n+1)}$ at   time   $t_{n+1} = t_n + \Delta t$ as follows
\begin{equation}
\mathbf y^{(n+1)}=\mathbf y^{(n)}+\Delta t\sum_{i=1}^sb_i\mathbf k_i,\quad 
\mathbf k_i=\mathbf f\left(t_n+c_i\Delta t,\mathrm{~}\mathbf y^{(n)}+\Delta t\sum_{j=1}^sa_{ij}\mathbf k_j\right),\quad i=1,\ldots,s.
\label{eq:irk_ynp1}
\end{equation}
In principle, \eqref{eq:irk_ynp1} can be directly solved as a nonlinear system involving the stage increments $\{ \mathbf{k}_i\}_{i=1}^s$. However, it is beneficial for computational implementation to reformulate the method in terms of intermediate stage states explicitly. This equivalent formulation  facilitates the definition of consistent residuals at the discrete level in the FSI problem. Towards this end, 
we introduce the stage times $t^{(n+c_i)} := t_n+c_i\Delta t$ and the stage states 
 $\mathbf y^{(n+c_i)} := \mathbf y^{(n)}+\Delta t\sum\limits_{j=1}^sa_{ij}\mathbf k_j$. The elimination of the stage increments $\mathbf k_j$ in  \eqref{eq:irk_ynp1} yields the equivalent stage form
\begin{equation}
\mathbf y^{(n+c_i)}=\mathbf y^{(n)}+\Delta t\sum_{j=1}^sa_{ij} 
\mathbf f\left(t^{(n+c_j)},\mathbf y^{(n+c_j)}\right)\quad\mathrm{~for~}i=1,\ldots,s,
\label{eq:irk_yni}
\end{equation} 
for the intermediate states
$\{ \mathbf y^{(n+c_i)} \}_i$, which is 
followed by the projection step to define the next state $\mathbf y^{(n+1)}$,
\begin{equation}
\mathbf y^{(n+1)}=\mathbf y^{(n)}+\Delta t\sum_{i=1}^sb_i \mathbf f\left(t^{(n+c_i)},\mathbf y^{(n+c_i)}\right).
\label{eq:irk_ynp1_new}
\end{equation}
In this work,  we adopt Radau-IIA IRK methods \cite{hairer2006geo,etienne2009perspective,cori2015high}, which are stiffly accurate: $a_{si}=b_i$ ($i=1,\ldots, s$) and $c_s=1$.
Consequently, the final stage coincides with the step solution, $\mathbf y^{(n+1)}=\mathbf y^{(n+c_s)}$, and the separate projection step \eqref{eq:irk_ynp1_new} can be omitted.

We now apply the IRK time integration to the semi-discrete FSI system \eqref{eq:FSI_DAE_abstract}, and obtain
the component-wise stage equations
\begin{equation}
\left\{
\begin{array}{rl}
\mathbf y_\Gamma^{(n+c_i)}
&= \mathbf y_\Gamma^{(n)}
   - \Delta t \sum_{j=1}^s a_{ij}
     \Bigl[
       \mathbf P_{\mathrm f,\Gamma}\,
       \mathbf R'_{\mathrm f}\bigl(\mathbf u_{\mathrm f}^{(n+c_j)},\mathbf p_{\mathrm f}^{(n+c_j)},
        \mathbf{P}_{\rm s,\Gamma}                                   \mathbf d_{\mathrm s}^{(n+c_j)}\bigr)
       + \mathbf P_{\mathrm s,\Gamma}\,
         \mathbf R'_{\mathrm s}\bigl(
                   \mathbf d_{\mathrm s}^{(n+c_j)}\bigr)
     \Bigr], \\[1pt]
\mathbf y_{\rm{ f,in}}^{(n+c_i)}
&= \mathbf y_{\rm{ f,in}}^{(n)}
   - \Delta t \sum_{j=1}^s a_{ij}
     \mathbf P_{\rm{ f,in}}\,
     \mathbf R'_{\mathrm f}\bigl(\mathbf u_{\mathrm f}^{(n+c_j)},\mathbf p_{\mathrm f}^{(n+c_j)},
                                  \mathbf{P}_{\rm s,\Gamma}
                                 \mathbf d_{\mathrm s}^{(n+c_j)}\bigr), \\[1pt]
\mathbf y_{\rm{ s,in}}^{(n+c_i)}
&= \mathbf y_{\rm{ s,in}}^{(n)}
   - \Delta t \sum_{j=1}^s a_{ij}
     \mathbf P_{\rm{ s,in}}\,
     \mathbf R'_{\mathrm s}\bigl(\mathbf d_{\mathrm s}^{(n+c_j)}\bigr), \\[1pt]
\mathbf d_{\mathrm s}^{(n+c_i)}
&= \mathbf d_{\mathrm s}^{(n)}
   + \Delta t \sum_{j=1}^s a_{ij}
     \mathbf u_{\mathrm s}^{(n+c_j)}, \\[1pt]
 \mathbf P_{\mathrm f,\Gamma}\,\mathbf u_{\mathrm f}^{(n+c_i)} &= 
\mathbf P_{\mathrm s,\Gamma}\,\mathbf u_{\mathrm s}^{(n+c_i)}, \\[1pt]
\mathbf B_{\rm f}\,\mathbf u_{\rm f}^{(n+c_i)} &= \mathbf 0, \\[1pt]
\mathbf P_{\rm{ f,dir}}\,\mathbf u_{\mathrm f}^{(n+c_i)} & =\mathbf u_{\rm{ f,dir}}^{(n+c_i)},\quad
\mathbf P_{\rm{ f,dir}}\,\mathbf d_{\mathrm s}^{(n+c_i)} =\mathbf d_{\rm{ f,dir}}^{(n+c_i)}, 
\end{array}
\right.
\label{eq:FSI_fully_discrete_stage_full}
\end{equation}
If we substitute the definitions of the masked momentum variables (cf. section \ref{sec:semi_discrete_FSI}),
\[
  \mathbf y_\Gamma^{(\cdot)} =
    \mathbf P_{\mathrm f,\Gamma}\mathbf M_{\mathrm f}^{(\cdot)}\mathbf u_{\mathrm f}^{(\cdot)}
    + \mathbf P_{\mathrm s,\Gamma}\mathbf M_{\mathrm s}\mathbf u_{\mathrm s}^{(\cdot)},\quad
  \mathbf y_{\rm{ f,in}}^{(\cdot)} =
    \mathbf P_{\rm{ f,in}}\mathbf M_{\mathrm f}^{(\cdot)}\mathbf u_{\mathrm f}^{(\cdot)},\quad
  \mathbf y_{\rm{ s,in}}^{(\cdot)} =
    \mathbf P_{\rm{ s,in}}\mathbf M_{\mathrm s}\mathbf u_{\mathrm s}^{(\cdot)},
\]
into the first three  relations of \eqref{eq:FSI_fully_discrete_stage_full}, and we rearrange the terms, we obtain the discrete fully discrete formulation:
at each stage $i=1,\ldots,s$ 
\begin{equation}
\left\{
\begin{aligned}
&\mathbf P_{\mathrm f,\Gamma}\,
\mathbf R_{\rm f}^{(n+c_i)}+
\mathbf P_{\mathrm s,\Gamma}\,
\mathbf R_{\rm s}^{(n+c_i)}= 0, \\[1mm]
& \mathbf P_{\mathrm f,\Gamma}\,\mathbf u_{\mathrm f}^{(n+c_i)} = 
\mathbf P_{\mathrm s,\Gamma}\,\mathbf u_{\mathrm s}^{(n+c_i)},\\[1mm]
&\mathbf P_{\mathrm f,\mathrm{in}}\,
\mathbf R_{\rm f}^{(n+c_i)}= 0,\\[1mm]
&\mathbf B_{\rm f}^{(n+c_i)}\mathbf u_{\rm f}^{(n+c_i)} = 0,\\[1mm]
&\mathbf P_{\mathrm s,\mathrm{in}}\,
\mathbf R_{\rm s}^{(n+c_i)}= 0,\\[1mm]
&
\mathbf P_{\mathrm f,\mathrm{dir}}\mathbf u_{\mathrm f}^{(n+c_i)}=\mathbf u_{\mathrm f,\mathrm{dir}}^{(n+c_i)},\quad
\mathbf P_{\mathrm s,\mathrm{dir}}\mathbf d_{\mathrm s}^{(n+c_i)}=\mathbf d_{\mathrm s,\mathrm{dir}}^{(n+c_i)},
\end{aligned}
\right.
\label{eq:FSI_fully_discrete_stage}
\end{equation}
where $\mathbf R_{\rm f}^{(n+c_i)}$ and $\mathbf R_{\rm s}^{(n+c_i)}$ denote the fully
discrete fluid and solid momentum residuals at the stage time $t^{(n+c_i)}$. 
Their explicit expressions are provided in
 \ref{sec:alt_derive_fully_discrete}.  
For completeness, in
 \ref{sec:alt_derive_fully_discrete} we prove that \eqref{eq:FSI_fully_discrete_stage} can be 
obtained by first discretizing the fluid and solid equations separately with the IRK scheme and subsequently enforcing the interface coupling conditions. Both approaches are algebraically equivalent and yield the same fully discrete coupled system.

We further note that the fully discrete formulation \eqref{eq:FSI_fully_discrete_stage_full} involves four
differential components in the state vector $\mathbf y$. After
applying the IRK scheme, the kinematic stage equations for the displacement
$\mathbf d_{\mathrm s}^{(n+c_i)}$ are used to eliminate the solid stage velocities
$\mathbf u_{\mathrm s}^{(n+c_i)}$ from the solid momentum residual
$\mathbf R_{\mathrm s}^{(n+c_i)}$, for $i=1,\ldots, s$. As a consequence,
\eqref{eq:FSI_fully_discrete_stage_full}$_4$  is implicitly encoded in the definition of
$\mathbf R_{\mathrm s}^{(n+c_i)}$, and the resulting fully discrete stage system
\eqref{eq:FSI_fully_discrete_stage} is expressed in terms
of three momentum blocks plus the algebraic constraints, in direct analogy with the
semi-discrete HF formulation \eqref{eq:FSI_semidiscrete}.

Note that the equations at each IRK stage are not independent; rather, they are inherently coupled across stages through the off-diagonal Butcher tableau coefficients. This inter-stage coupling implies that we need to solve a large implicit nonlinear system that encompass all stages. To explicitly illustrate and address this coupling structure, we introduce stage-stacked vectors. Specifically, we define stage–stacked unknowns as:
\[
\mathbf U_{\mathrm f}:=\operatorname{col}\!\big(\mathbf u_{\mathrm f}^{(n+c_1)},\ldots,\mathbf u_{\mathrm f}^{(n+c_s)}\big)\in\mathbb{R}^{sN_{\mathrm u}},\quad
\mathbf p_{\mathrm f}:=\operatorname{col}\!\big( \mathbf p_{\mathrm f}^{(n+c_1)},\ldots,
\mathbf p_{\mathrm f}^{(n+c_s)}\big)\in\mathbb{R}^{sN_{\mathrm p}},
\]
\[
\mathbf D_{\mathrm s}:=\operatorname{col}\!\big(\mathbf d_{\mathrm s}^{(n+c_1)},\ldots,\mathbf d_{\mathrm s}^{(n+c_s)}\big)\in\mathbb{R}^{sN_{\mathrm s}},
\quad
\mathbf U_{\mathrm s}:=\operatorname{col}\!\big(\mathbf u_{\mathrm s}^{(n+c_1)},\ldots,\mathbf u_{\mathrm s}^{(n+c_s)}\big)\in\mathbb{R}^{sN_{\mathrm s}},
\]
where the operator $\operatorname{col}(\cdot)$ denotes vertical concatenation. Correspondingly, we also define stage-stacked residual vectors
\[
          \mathcal R_{\rm f}(\mathbf U_{\rm f},\mathbf p_{\rm f},\mathbf D_{\rm s})
          :=\mathrm{col}\bigl(
          \mathbf R_{\rm f}^{(n+c_1)},\ldots,
          \mathbf R_{\rm f}^{(n+c_s)}\bigr),
          \qquad
          \mathcal R_{\rm s}(\mathbf D_{\rm s})
          :=\mathrm{col}\bigl(
          \mathbf R_{\rm s}^{(n+c_1)},\ldots,
          \mathbf R_{\rm s}^{(n+c_s)}\bigr).
        \]
We further define stage-stacked Boolean mask matrices and the discrete divergence operator:
\[
\mathcal P_{\bullet}:=I_s\otimes \mathbf P_{\bullet},\qquad
\mathcal B_{\rm f}:=\operatorname{blkdiag}\!\big(\mathbf B_{\rm f}^{(n+c_1)},\ldots,\mathbf B_{\rm f}^{(n+c_s)}\big),
\]
where $I_s\otimes \mathbf P_{\bullet} \in \mathbb{R}^{s N_{\bullet}^{\mathrm{row}} \times s N_{\bullet}^{\mathrm{col}}}$ is the Kronecker product of the identity matrix $I_s$ and $\mathbf P_\bullet \in \mathbb{R}^{N_{\bullet}^{\mathrm{row}} \times
N_{\bullet}^{\mathrm{col}}}$,
\[
\bigl(I_s\otimes \mathbf P_{\bullet}\bigr)_{i,j} =
\bigl(I_s\bigr)_{r,q}
\bigl(\mathbf P_{\bullet}\bigr)_{v,w},
\qquad
i = (r-1) N_{\bullet}^{\mathrm{row}} + v,\quad
j = (q-1) N_{\bullet}^{\mathrm{col}} + w,
\]
with
\[
1 \le r,q \le s,\qquad
1 \le v \le N_{\bullet}^{\mathrm{row}},\qquad
1 \le w \le N_{\bullet}^{\mathrm{col}}.
\]

As mentioned previously, the solid velocity variables can be eliminated by substituting the kinematic relation \eqref{eq:FSI_fully_discrete_stage_full}$_4$ into the solid momentum equations \eqref{eq:FSI_fully_discrete_stage_full}$_3$ and the interface condition \eqref{eq:FSI_fully_discrete_stage_full}$_1$: this substitution leads to a formulation that involves only solid displacements as unknowns. In more detail, we find
\[
\mathbf U_{\rm s}
=\frac{1}{\Delta t}\,(\mathbf A^{-1}\otimes I_{N_{\rm s}})\,\mathbf D_{\rm s}
-\frac{1}{\Delta t}\,\big((\mathbf A^{-1}\mathbf 1_{s})\otimes I_{N_{\rm s}}\big)\,\mathbf d_{\rm s}^{(n)},
\]
where
$\mathbf A \in \mathbb{R}^{s\times s}$ is the Butcher matrix associated with the RK scheme and 
$\mathbf 1_{s}\in\mathbb{R}^s$ is the vector of ones. For simplicity, we rewrite this relation compactly as
\begin{equation}
\mathbf U_{\rm s}
=\mathcal D_{\rm s,\Delta t}\,\mathbf D_{\rm s}.
\label{eq:solid_kinematic_eqn_stacked}
\end{equation}
Note that the discrete time differentiation operator $\mathcal D_{\rm s,\Delta t}$ changes at each time step; in the remainder, we omit this explicit dependence to shorten  notation.

With the definitions and notation established above, we can now formulate the global fully discrete monolithic system, which couples all IRK  stages simultaneously. At each time step, the fully discrete FSI problem consists of solving 
for the stacked unknowns $(\mathbf U_{\rm f},\,\mathbf p_{\rm f},\,\mathbf D_{\rm s})$:
\begin{equation}\label{eq:FSI_IRK_monolithic}
\left\{
\begin{aligned}
&\mathcal P_{\rm f,\Gamma}\,\mathcal R_{\rm f}(\mathbf U_{\rm f},\,\mathbf p_{\rm f},\,
\mathcal P_{\rm s,\Gamma}
\mathbf D_{\rm s})+\mathcal P_{\rm s,\Gamma}\,\mathcal R_{\rm s}(\mathbf D_{\rm s})= 0,  \\[2pt]
&
\mathcal P_{\rm f,\Gamma}\,\mathbf U_{\rm f}=\mathcal P_{\rm s,\Gamma}\, \mathbf U_{\rm s}, 
\\[2pt]
&\mathcal P_{\rm f,\mathrm{in}}\,\mathcal R_{\rm f}(\mathbf U_{\rm f},\,\mathbf p_{\rm f},\,
\mathcal P_{\rm s,\Gamma}
\mathbf D_{\rm s})=0, \\
&\mathcal B_{\rm f}\,\mathbf U_{\rm f}=0,\\[2pt]
 & \mathcal P_{\rm s,\mathrm{in}}\,\mathcal R_{\rm s}(\mathbf D_{\rm s})=0,\\[2pt]
&\mathcal P_{\rm f,\mathrm{dir}}\,\mathbf U_{\rm f}=\mathbf U_{\rm f,\mathrm{dir}},\qquad
 \mathcal P_{\rm s,\mathrm{dir}}\,\mathbf D_{\rm s}=\mathbf D_{\rm s,\mathrm{dir}}.
\end{aligned}
\right.
\end{equation}
Note that \eqref{eq:FSI_IRK_monolithic} shares the same structure with the standard BDF-Newmark implicit time discretization of FSI problems, but it has  $s$-times more equations and unknowns.
Eq. 
\eqref{eq:FSI_IRK_monolithic}$_1$ enforces dynamic equilibrium of normal
stresses at the fluid-structure interface, 
Eq. 
\eqref{eq:FSI_IRK_monolithic}$_2$ 
enforces
kinematic continuity (matching of normal velocities) at the interface, and the
remaining equations impose the fluid and solid governing equations in the
interior together with Dirichlet boundary conditions. 

\subsection{Solution to the  fully discrete FSI system \eqref{eq:FSI_IRK_monolithic}}
\label{sec:opt_SQP}

In order to describe the application of static condensation to \eqref{eq:FSI_IRK_monolithic}, we introduce notation:
\begin{equation}
\label{eq:variables_static_condensation}
\mathbf g := \mathcal P_{\rm s,\Gamma}\mathbf D_{\rm s},
\quad
\mathbf{W}_{\rm f} := {\rm col}\left( 
\mathcal P_{\rm f,in}\mathbf U_{\rm f}  \, , \,
\mathbf{p}_{\rm f}
\right),
\quad
\mathbf{D}_{\rm s}^{\rm in} := 
\mathcal P_{\rm s,in}\mathbf D_{\rm s}.    
\end{equation}
Recalling 
\eqref{eq:solid_kinematic_eqn_stacked} and the kinematic condition
\eqref{eq:FSI_IRK_monolithic}$_2$, we find that 
$\mathbf U_{\rm f}$ is uniquely determined by 
$\mathbf W_{\rm f}$, $\mathbf{g}$ and the solution at the previous time step;
similarly, $\mathbf D_{\rm s}$ is uniquely determined by 
$\mathbf D_{\rm s}^{\rm in}$ and $\mathbf{g}$. 
We can hence introduce the residuals:
\begin{equation}
\label{eq:residuals_static_condensation}
\left\{
\begin{array}{l}
\displaystyle{
\mathcal R_{\Gamma} \left( 
\mathbf W_{\rm f}, 
\mathbf g, 
\mathbf D_{\rm s}^{\rm in}
\right)
:=
\mathcal P_{\rm f,\Gamma}\,\mathcal R_{\rm f}(\mathbf U_{\rm f},\,\mathbf p_{\rm f},\,\mathcal P_{\rm s,\Gamma}\mathbf D_{\rm s})+\mathcal P_{\rm s,\Gamma}\,\mathcal R_{\rm s}(\mathbf D_{\rm s}), 
}
\\[3mm]
\displaystyle{
\mathcal R_{\rm s,in} \left( 
\mathbf D_{\rm s}^{\rm in}, 
\mathbf g
\right)
:=
\left[
\begin{array}{l}
 \mathcal P_{\rm s,\mathrm{in}}\,\mathcal R_{\rm s}(\mathbf D_{\rm s}) \\[3mm]
  \mathcal P_{\rm s,\mathrm{dir}}\,\mathbf D_{\rm s} - \mathbf D_{\rm s,\mathrm{dir}}
  \\
\end{array}
\right],
\quad
\mathcal R_{\rm f,in} \left( 
\mathbf W_{\rm f}, 
\mathbf g
\right)
:=
\left[
\begin{array}{l}
 \mathcal P_{\rm f,\mathrm{in}}\,\mathcal R_{\rm f}(\mathbf U_{\rm f},\,\mathbf p_{\rm f},\,\mathcal P_{\rm s,\Gamma}\mathbf D_{\rm s}) 
 \\[3mm]
 \mathcal B_{\rm f,\Gamma}\, \mathbf U_{\rm f}  \\[3mm]
 \mathcal P_{\rm f,\mathrm{dir}}\,\mathbf U_{\rm f} - \mathbf U_{\rm f,\mathrm{dir}}
  \\
\end{array}
\right].
}
\\
\end{array}
\right.
\end{equation}
In conclusion, we
 rewrite  \eqref{eq:FSI_IRK_monolithic} as
\begin{equation}
\label{eq:FSI_IRK_monolithic_compact}
\left\{
\begin{array}{l}
\displaystyle{
\mathcal R_{\Gamma} \left( 
\mathbf W_{\rm f}, 
\mathbf g, 
\mathbf D_{\rm s}^{\rm in}
\right) = 0,
}
 \\[3mm]
\displaystyle{
\mathcal R_{\rm s,in} \left( 
\mathbf D_{\rm s}^{\rm in}, 
\mathbf g
\right) = 0,
}
\\[3mm]
\displaystyle{
\mathcal R_{\rm f,in} \left( 
\mathbf W_{\rm f}, 
\mathbf g
\right) = 0.
}
  \\
\end{array}
\right.
\end{equation}
If we apply the Newton's method to \eqref{eq:FSI_IRK_monolithic_compact}, we obtain
\begin{equation}
\label{eq:newton_method_intermediate}
\left\{
\begin{array}{l}
\displaystyle{
\mathcal R_{\Gamma}^{(j)} 
+ 
\mathcal{J}_{\Gamma, \rm f}
\delta 
\mathbf{W}_{\rm f}
+ 
\mathcal{J}_{\Gamma, \rm s}
\delta 
\mathbf{D}_{\rm s}^{\rm in}
+ 
\mathcal{J}_{\Gamma, \rm g}
\delta 
\mathbf{g} = 0,
\quad
{\rm where} \;
\mathcal{J}_{\Gamma, \rm f} = 
\dfrac{\partial \mathcal R_{\Gamma}^{(j)}}{\partial \mathbf{W}_{\rm f}},
\quad
\mathcal{J}_{\Gamma, \rm s} =
\dfrac{\partial \mathcal R_{\Gamma}^{(j)}}{\partial \mathbf{D}_{\rm s}^{\rm in}},
\quad
\mathcal{J}_{\Gamma, \rm g} =
\dfrac{\partial \mathcal R_{\Gamma}^{(j)}}{\partial \mathbf{g}};
}
 \\[3mm]
\displaystyle{
\mathcal R_{\rm s,in}^{(j)} 
+ 
\mathcal{J}_{\rm s,s}
\delta 
\mathbf{D}_{\rm s}^{\rm in}
+ 
\mathcal{J}_{\rm s,g}
\delta 
\mathbf{g} = 0,
\quad
{\rm where} \;
\mathcal{J}_{\rm s,s}
=
\dfrac{\partial \mathcal R_{\rm s,in}^{(j)}}{\partial \mathbf{D}_{\rm s}^{\rm in}}, \quad
\mathcal{J}_{\rm s,g} = 
\dfrac{\partial \mathcal R_{\rm s,in}^{(j)}}{\partial \mathbf{g}};
}
\\[3mm]
\displaystyle{
\mathcal R_{\rm f,in}^{(j)} 
+ 
\mathcal{J}_{\rm f,f} 
\delta 
\mathbf{W}_{\rm f}
+ 
\mathcal{J}_{\rm f,g}
\delta 
\mathbf{g} = 0,
\quad
{\rm where} \;
\mathcal{J}_{\rm f,f} = 
\dfrac{\partial \mathcal R_{\rm f,in}^{(j)}}{\partial \mathbf{W}_{\rm f}},
\;
\mathcal{J}_{\rm f,g}
=
\dfrac{\partial \mathcal R_{\rm f,in}^{(j)}}{\partial \mathbf{g}};
}\\
\end{array}
\right.
\end{equation}
and  $j$ denotes the Newton's iteration.
The latter system can be solved using static condensation: first,  we 
exploit the second and the third equations to express 
the local fluid and the solid states in terms of the interface control $\delta 
\mathbf{g}$,
\begin{subequations}
\label{eq:static_condensation_explained}
\begin{equation}
\label{eq:static_condensation_explained_a}
\left\{
\begin{array}{ll}
\displaystyle{
\delta 
\mathbf{D}_{\rm s}^{\rm in}
=
\mathbf{D}_{\star}^{(j)}
+
\mathcal{J}_{\rm g \to s}^{(j)}
\delta 
\mathbf{g}
}
&
\displaystyle{
{\rm where}
\;
\mathbf{D}_{\star}^{(j)}
=
-
\left(
  \mathcal J_{\rm s, s}^{(j)}
\right)^{-1}
\mathcal R_{\rm s, in}^{(j)},
\quad
\mathcal{J}_{\rm g \to s}^{(j)}
=
-
\left(
  \mathcal J_{\rm s, s}^{(j)}
\right)^{-1}
\mathcal J_{\rm s, g}^{(j)};
} \\[3mm]
\displaystyle{
\delta 
\mathbf{W}_{\rm f}
=
\mathbf{W}_{\star}^{(j)}
+
\mathcal{J}_{\rm g \to f}^{(j)}
\delta 
\mathbf{g}
} &
\displaystyle{{\rm where}
\;
\mathbf{W}_{\star}^{(j)}
=
-
\left(
  \mathcal J_{\rm f, f}^{(j)}
\right)^{-1}
\mathcal R_{\rm f, in}^{(j)},
\quad
\mathcal{J}_{\rm g \to f}^{(j)}
=
-
\left(
  \mathcal J_{\rm f, f}^{(j)}
\right)^{-1}
\mathcal J_{\rm f, g}^{(j)} 
}; \\
\end{array}
\right.  
\end{equation}
next, we substitute the expressions for $\delta 
\mathbf{D}_{\rm s}^{\rm in}$ and
$\delta 
\mathbf{W}_{\rm f}$ in the first equation to find  an independent  linear system for the interface control $\delta 
\mathbf{g}$:
\begin{equation}
\label{eq:static_condensation_explained_b}
 \mathcal S_{\rm sc}^{(j)} 
\delta 
\mathbf{g}
= 
\mathcal R_{\rm sc}^{(j)} ,
\qquad
{\rm where} \;
\left\{
\begin{array}{l}
\displaystyle{ \mathcal S_{\rm sc}^{(j)}  = 
\mathcal J_{\Gamma, \rm g}^{(j)}
 -
\mathcal J_{\Gamma, \rm s}^{(j)}
\mathcal J_{\rm g \to s}^{(j)}
 -
\mathcal J_{\Gamma, \rm f}^{(j)}
 \mathcal J_{\rm g \to f}^{(j)}
}, \\[3mm]
\displaystyle{
\mathcal R_{\rm sc}^{(j)} =
- 
\mathcal R_{\Gamma}^{(j)} 
- 
  \mathcal J_{\Gamma, \rm  s}^{(j)}
\mathbf{D}_{\star}^{(j)}
-
  \mathcal J_{\Gamma, \rm  f}^{(j)}
\mathbf{W}_{\star}^{(j)} 
}. \\
\end{array}
\right.
\end{equation}
In conclusion, at each Newton's iteration, we first compute 
$\mathbf{D}_{\star}^{(j)},
\mathcal{J}_{\rm g \to s}^{(j)}$ and 
$\mathbf{W}_{\star}^{(j)},
\mathcal{J}_{\rm g \to f}^{(j)}$ using 
\eqref{eq:static_condensation_explained_a};
next, we compute $\delta \mathbf{g}$ using 
 \eqref{eq:static_condensation_explained_b}; 
 finally, we compute 
 $\delta 
\mathbf{D}_{\rm s}^{\rm in}$ and
$\delta 
\mathbf{W}_{\rm f}$ using again
\eqref{eq:static_condensation_explained_a}.
\end{subequations}

We can interpret the solution to \eqref{eq:FSI_IRK_monolithic_compact} as the solution to the constrained minimization statement:
\begin{equation}
\label{eq:FSI_IRK_optimization_based}
\min_{ 
\mathbf W_{\rm f}, 
\mathbf g, 
\mathbf D_{\rm s}^{\rm in}
 }
 \,
 \big\|
 \mathcal R_{\Gamma} \left( 
\mathbf W_{\rm f}, 
\mathbf g, 
\mathbf D_{\rm s}^{\rm in}
\right) 
  \big\|_2
\quad
{\rm subject \, to} \;
\left\{
\begin{array}{l}
\displaystyle{
\mathcal R_{\rm s,in} \left( 
\mathbf D_{\rm s}^{\rm in}, 
\mathbf g
\right) = 0,
}
\\[3mm]
\displaystyle{
\mathcal R_{\rm f,in} \left( 
\mathbf W_{\rm f}, 
\mathbf g
\right) = 0.
}
  \\
\end{array}
\right.    
\end{equation}
Interestingly, if we apply the same 
sequential quadratic programming (SQP) method as in \cite{taddei2024non,taddei2025optimization} to \eqref{eq:FSI_IRK_optimization_based}, at each iteration we find the  linear problem 
\eqref{eq:newton_method_intermediate}. Note, however, that while \eqref{eq:FSI_IRK_monolithic_compact} critically relies on the assumption that the fluid and the solid meshes are conforming at the interface, the minimization statement \eqref{eq:FSI_IRK_optimization_based} can be readily generalized to non-conforming grids.

\section{Model order reduction}
\label{sec:rom}

In this section, we construct a reduced-order model (ROM) for the FSI problem that preserves the semi-discrete energy balance and the interface coupling structure of the full-order  formulation. We first introduce suitable reduced spaces constructed from high-fidelity snapshots using proper orthogonal decomposition (POD).
The construction of the reduced spaces follows the bubble-port decomposition approach as in the static-condensation reduced basis element (scRBE) method \cite{huynh2013static,eftang2013port,ebrahimi2024hyperreduced}. 
We then derive the semi-discrete ROM and prove an energy balance analogous to the full-order model. Finally, we extend the ROM formulation to the fully discrete IRK setting and discuss the resolution of the ROM.

\subsection{Reduced spaces}
\label{sec:ROM_spaces}

High-fidelity snapshots are collected at each physical time step $k=0,\dots,k_{\max}$ and for each IRK stage $i=1,\dots,s$, i.e., we consider the dataset
\[
  \bigl\{\,
    \mathbf u_{\rm f}^{(k+c_i)},\;
    p_{\rm f}^{(k+c_i)},\;
    \mathbf d_{\rm s}^{(k+c_i)}
  \bigr\}_{k=0,\dots,k_{\max},\; i=1,\dots,s},
\]
where $\mathbf u_{\rm f}^{(k+c_i)}\in\mathbb R^{N_{\rm u}}$ and $p_{\rm f}^{(k+c_i)}\in\mathbb R^{N_{\rm p}}$ denote the discrete fluid velocity and pressure, and $\mathbf d_{\rm s}^{(k+c_i)}\in\mathbb R^{N_{\rm s}}$ the discrete solid displacement.
For long-time integration, we do not store the solution at all time steps to control memory costs.

We rely on POD \cite{sirovich1987turbulence,volkwein2011model} to construct the reduced-order spaces. 
We 
first build the port (interface) space: we collect interface displacement snapshots $  \mathbf d_{{\rm s},\Gamma}^{(k+c_i)} 
  := \mathbf P_{{\rm s},\Gamma}\mathbf d_{\rm s}^{(k+c_i)} \in \mathbb R^{N_{\rm c}}$, 
and apply   POD
based on the $L^2(\widetilde{\Gamma})$-inner product
to obtain the interface modes
$\{\boldsymbol\psi_j\}_{j=1}^{m}\subset\mathbb R^{N_{\rm c}}$. 
Each interface mode $\boldsymbol\psi_j$ is then harmonically extended into both the fluid and solid domains to construct the \emph{port bases}. Specifically, for each $j=1,\dots,m$ we solve a Laplace problem with Dirichlet data $\boldsymbol\psi_j$ on $\Gamma$ and homogeneous Dirichlet conditions on the remaining domain boundaries, for both the fluid and solid domains.  We denote these resulting discrete fields by $\mathbf e_{{\rm s},j}\in\mathbb R^{N_{\rm s}}$ and $\mathbf e_{{\rm f},j}\in\mathbb R^{N_{\rm u}}$, respectively. 
We collect these fields column-wise to obtain the port basis matrices:
\[
  \mathbf E_{\rm s} := [\mathbf e_{{\rm s},1}\;\cdots\;\mathbf e_{{\rm s},m}]
  \in\mathbb R^{N_{\rm s}\times m},
  \qquad
  \mathbf E_{\rm f} := [\mathbf e_{{\rm f},1}\;\cdots\;\mathbf e_{{\rm f},m}]
  \in\mathbb R^{N_{\rm u}\times m}.
\]
By construction, each column of $\mathbf E_{\rm s}$ 
shares the same trace on $\Gamma(t)$
with the corresponding column of 
$\mathbf E_{\rm f}$. 

We next define \emph{bubble spaces} for the fluid velocity and solid displacement. These spaces describe interior fields that satisfy homogeneous conditions at the interface and on the external Dirichlet boundaries. For the fluid velocity, 
we assume that the time-dependent Dirichlet boundary data on $\partial\Omega_{\rm f}(t)\setminus\Gamma(t)$ are prescribed in separable form
$
  \mathbf u_{\rm f,dir}(t) = c_{\rm dir}(t)\,\bar{\mathbf u}_{\rm f}^{\rm dir},
$
where $\bar{\mathbf u}_{\rm f}^{\rm dir}\in\mathbb R^{N_{\rm u}}$ is a fixed spatial lifting function obtained as the harmonic extension of a reference boundary profile.
More specifically, $\bar{\mathbf u}_{\rm f}^{\rm dir}$ solves a Laplace problem with Dirichlet data equal to the reference boundary profile on $\partial\Omega_{\rm f}\setminus\Gamma$ and homogeneous Dirichlet conditions on $\Gamma$. 
In addition, we introduce the fluid harmonic extension operator $\mathcal H_{\rm f,\Gamma}$,   
which maps an interface velocity trace to an interior velocity field with homogeneous Dirichlet data on the remaining boundary. Thus, for each fluid snapshot, we define
\[
  {\mathbf u}_{{\rm f},\Gamma}^{(k+c_i)}
  := \mathcal H_{\rm f,\Gamma}\bigl(\mathbf P_{{\rm f},\Gamma}\mathbf u_{{\rm f}}^{(k+c_i)}\bigr),
\]
and subsequently introduce the lifted fluid snapshot
\[
  \mathbf u_{\rm f,{\rm lift}}^{(k+c_i)}
  := \mathbf u_{\rm f}^{(k+c_i)}
     - c_{\rm dir}(t^{(k+c_i)})\,\bar{\mathbf u}_{\rm f}^{\rm dir}
     - \widetilde{\mathbf u}_{{\rm f},\Gamma}^{(k+c_i)}.
\]
By construction, this lifted snapshot vanishes on both the external Dirichlet boundaries and the interface $\Gamma(t)$. 
We rely on the $H^1(\widetilde{\Omega}_{\rm f})$ POD to devise the fluid velocity
bubble basis
\[
  \mathbf Z_{\rm f}\in\mathbb R^{N_{\rm u}\times n_{\rm u}}.
\]
The number of modes $n_{\rm u}$ is determined using an energy criterion.
Following the same procedure, we construct the lifted solid snapshots and then we obtain
 the solid bubble basis
\[
  \mathbf Z_{\rm s}\in\mathbb R^{N_{\rm s}\times n_{\rm s}},
\]
using  POD based on the 
$H^1(\widetilde{\Omega}_{\rm s})$ inner product.
Finally, we define a reduced-order space for the fluid pressure. We store the pressure snapshots 
$\{p_{\rm f}^{(k+c_i)}\}$ and apply POD
based on the $L^2(\widetilde{\Omega}_{\rm s})$ inner product
to obtain the pressure basis
\[
  \mathbf Z_{\rm p}\in\mathbb R^{N_{\rm p}\times n_{\rm p}}.
\]

To preserve inf-sup stability at the reduced level for the fluid, we apply a standard supremizer enrichment \cite{rozza2007stability}. 
For each pressure mode, we compute an associated velocity supremizer by solving a suitable auxiliary problem, and then we augment the velocity bubble basis $\mathbf Z_{\rm f}$ with the resulting 
$n_{\rm p}$ supremizer vectors. This enrichment step practically ensures a bounded inf-sup constant at the reduced-order level.

With the above reduced-order spaces defined, we  now state the continuous-time reduced approximation ansatz that will be used to construct a semi-discrete ROM. The reduced approximations for the solid displacement, fluid velocity, and fluid pressure read: 
\begin{subequations}\label{eq:RB_continuous}
\begin{align}
\widehat{\mathbf{d}}_{\rm s}(t) &=
      \mathbf Z_{\rm s}\,\alpha_{\rm s}(t)
     +\mathbf E_{\rm s}\,\beta(t),
\\
\widehat{\mathbf{u}}_{\rm f}(t) &=
      \bar{\mathbf u}_{\rm f}^{\rm dir}\,c_{\rm dir}(t)
     +\mathbf Z_{\rm f}\,\alpha_{\rm f}(t)
     +\mathbf E_{\rm f}\,\dot{\beta}(t),
\\
\widehat{\mathbf{p}}_{\rm f}(t) &=
      \mathbf Z_{\rm p}\,\alpha_{\rm p}(t),
\end{align}
\end{subequations}
where the time-dependent generalized coordinates are
$\alpha_{\rm f}(t)\in\mathbb R^{n_{\rm u}}$,  
$\alpha_{\rm s}(t)\in\mathbb R^{n_{\rm s}}$,  
$\alpha_{\rm p}(t)\in\mathbb R^{n_{\rm p}}$, and  
$\beta(t)\in\mathbb R^{m}$. 
The term $\mathbf E_{\rm s}\beta(t)$ explicitly describes the displacement of the fluid-structure interface, whereas $\mathbf E_{\rm f}\dot\beta(t)$ represents the fluid velocity field induced by this interface motion.
Since $\mathbf E_{\rm s}$ and $\mathbf E_{\rm f}$ share the same trace on $\Gamma(t)$, this choice automatically enforces kinematic compatibility at the interface. 

\begin{remark}
Since the ALE mapping is the solution to a linear parameter- and time-independent pseudo-elasticity problem,  $\Phi_{\rm f}$ is an affine function of the port modes. Therefore, its construction  and online evaluation is trivial. We notice that 
if either  
nonlinear elasticity models  for the ALE morphing need to be employed or the problem features geometric parameters, we should treat the mapping as an additional state variable and devise an offline/online strategy accordingly.
\end{remark}

\subsection{Semi-discrete ROM}
\label{sec:semi_discrete_ROM}

We now insert the reduced approximation \eqref{eq:RB_continuous} into the semi-discrete FSI equations \eqref{eq:FSI_semidiscrete} and project onto the port and bubble spaces to derive the semi-discrete ROM. 
We define the reduced fluid and solid residuals as
\[
  \widehat{\mathbf  R}_{\rm f}(\cdot):=
  \frac{d}{dt} \left(
     \mathbf M_{\rm f}\, {\widehat{\mathbf u}}_{\rm f} \right)
    +\mathbf R_{\rm f}'(\cdot),
\qquad
  \widehat{\mathbf R}_{\rm s}(\cdot):=
     \mathbf M_{\rm s}\,\ddot{\widehat{\mathbf d}}_{\rm s}
    +\mathbf R_{\rm s}'(\cdot). 
\]
Galerkin projection of the semi-discrete FSI equations onto the port and bubble spaces yields the ROM:
\begin{equation}
\left\{
\begin{aligned}
& \mathbf E_{\rm s}^{\top}\,\widehat{\mathbf  R}_{\rm s}
 +\mathbf E_{\rm f}^{\top}\,\widehat{\mathbf  R}_{\rm f}=0, \\
& \mathbf Z_{\rm f}^{\top}\,\widehat{\mathbf  R}_{\rm f}=0,
\qquad
  \mathbf Z_{\rm p}^{\top}\,\mathbf B_{\rm f}\,\widehat{\mathbf u}_{\rm f}=0,
\\
& \mathbf Z_{\rm s}^{\top}\,\widehat{\mathbf  R}_{\rm s}=0. 
\end{aligned}
\right.
\label{eq:semi_discrete_ROM}
\end{equation}
In this formulation, the first equation corresponds to the dynamic equilibrium condition at the interface, while 
the second and the third equations enforce the fluid and solid momentum balances, respectively, within their corresponding bubble spaces. The 
kinematic continuity at the interface is automatically satisfied by construction.

In view of the stability analysis, we provide the variational formulation associated with \eqref{eq:semi_discrete_ROM}. Towards this end, we introduce the spaces
$$
\left\{
\begin{array}{l}
\displaystyle{
\widehat{V}_{\rm f}(t) \, = \,
\left\{
\tilde{u}\circ \Phi_{\rm f}^{-1}(t) \, : \,
\tilde{u} = \bar{u}_{\rm f}^{\rm dir} c_{\rm dir}(t)
+ 
\sum_{i=1}^{n_{\rm u}} \alpha_i \tilde{\zeta}_{{\rm f},i}
+ 
\sum_{j=1}^{m} \beta_j \tilde{e}_{{\rm f},j},
\; \;
\boldsymbol{\beta}\in \mathbb{R}^m, \;\;
\boldsymbol{\alpha}\in \mathbb{R}^{n_{\rm u}}
\right\},
} \\[3mm]
\displaystyle{
\widehat{V}_{\rm f, 0}(t) \, = \,
\left\{
\tilde{u}\circ \Phi_{\rm f}^{-1}(t) \, : \,
\tilde{u}  \in {\rm span} \{ \tilde{\zeta}_{i, \rm f} \}_{i=1}^{n_{\rm u}} \cup 
{\rm span} \{ \tilde{e}_{{\rm f},j} \}_{j=1}^{m} 
\right\}, 
\quad
\widehat{Q}_{\rm f}(t) \, = \,
\left\{
\tilde{q}\circ \Phi_{\rm f}^{-1}(t) \, : \,
\tilde{q}  \in  {\rm span} \{ \tilde{\zeta}_{{\rm p}, i} \}_{i=1}^{n_{\rm p}}
\right\},
}
\\[3mm]
\displaystyle{
\widehat{V}_{\rm s}  \, = \,
 {\rm span} \{ \tilde{\zeta}_{{\rm s},i} \}_{i=1}^{n_{\rm s}} \cup 
{\rm span} \{ \tilde{e}_{{\rm s},j} \}_{j=1}^m
}, \\
\end{array}
\right.
$$
where $\{ \tilde{\zeta}_{{\rm f},i} \}_{i=1}^{n_{\rm u}}$,
$\{ \tilde{\zeta}_{{ \rm p}, i} \}_{i=1}^{n_{\rm p}}$,
$\{ \tilde{\zeta}_{{\rm s},i} \}_{i=1}^{n_{\rm s}}$,
$\{ \tilde{e}_{{\rm f}, j} \}_{j=1}^{m}$,
$\{ \tilde{e}_{{\rm s}, j} \}_{j=1}^{m}$, 
are the FE fields associated with the POD bases introduced in section \ref{sec:ROM_spaces}.
We denote by  $\widehat{u}_{\rm f} \in \widehat{V}_{\rm f}(t)$,
$\widehat{p}_{\rm f} \in \widehat{Q}_{\rm f}(t)$, 
 $\widehat{d}_{\rm s} \in \widehat{V}_{\rm s}$ the reduced-order estimates associated with the FE vectors $\widehat{\mathbf{u}}_{\rm f} \in \mathbb{R}^{N_{\rm u}}$,
$\widehat{\mathbf{p}}_{\rm f} \in \mathbb{R}^{N_{\rm p}}$,
$\widehat{\mathbf{d}}_{\rm s} \in \mathbb{R}^{N_{\rm s}}$ in \eqref{eq:RB_continuous}, respectively. For any $t>0$, the reduced-order estimates satisfy:
\begin{equation}
\label{eq:weak_formulation_rom}
\left\{
\begin{array}{l}
\displaystyle{
\frac{d}{dt}\int_{ \Omega_{\rm f}(t)   }
\rho_{\rm f} 
 \widehat{u}_{\rm f} 
 \cdot  v \, dx
+
R_{\rm f}'( \widehat{u}_{\rm f}, \widehat{p}_{\rm f}, \bs\omega_{\rm f}, v) 
+
\int_{ \widetilde{\Omega}_{\rm s}   }
\widetilde\rho_{\rm s} 
\frac{\partial^2 \widehat{d}_{\rm s} }{\partial t^2}
\cdot w \, dx
+
R_{\rm s}'(\widehat{d}_{\rm s},  w)
 = 0
};
\\[3mm]
\displaystyle{
b_{\rm f}(\widehat{u}_{\rm f}, q)
 = 0
};
\\[3mm]
\displaystyle{
\Phi_{\rm f}  \big|_{\widetilde{\Gamma}} = \texttt{id} + \widehat{d}_{\rm s}\big|_{\widetilde{\Gamma}},
\quad
\widehat{u}_{\rm f} \circ \Phi_{\rm f} \big|_{\widetilde{\Gamma}}
=
\frac{\partial \widehat{d}_{\rm s}}{\partial t}
\big|_{\widetilde{\Gamma}},
\quad
\bs\omega_{\rm f} =
\dfrac{\partial \Phi_{\rm f}}{\partial t}
},
  \\
\end{array}
\right.
\end{equation}
for any $v\in \widehat{V}_{\rm f,0}(t)$, $q\in \widehat{Q}_{\rm f}(t)$ and $w\in \widehat{V}_{\rm s}$ such that
$v\circ \Phi_{\rm f}(t) \big|_{\widetilde{\Gamma}} = w\big|_{\widetilde{\Gamma}}$. Note that \eqref{eq:weak_formulation_rom} is the Galerkin projection of \eqref{eq:weak_formulation_a}. Next, we define the kinetic energies in the fluid and the solid domain, and the solid potential:
$$
\widehat{T}_{\rm f}(t) : =
\dfrac{1}{2} \int_{\Omega_{\rm f}(t)} \rho_{\rm f} \big| \widehat{u}_{\rm f}(t) \big|^2 \, dx,
\quad
\widehat{T}_{\rm s}(t) : =
\dfrac{1}{2} \int_{\widetilde{\Omega}_{\rm s}} \rho_{\rm s} \big|  \partial_t \widehat{d}_{\rm s}(t) \big|^2 \, dx,
\quad
\widehat{\mathcal W}_{\rm s}(t) : =
\int_{\widetilde{\Omega}_{\rm s}} W(\widehat{\mathbf{E}}_{\rm s} ) \, dx,
$$
where $\widehat{\mathbf{E}}_{\rm s}$ is the Green-Lagrange tensor associated with the predicted displacement and $W$ is the density of elastic energy associated with the prescribed constitutive model  --- i.e., $\widetilde{\boldsymbol{\Sigma}}_{\rm s}  = \dfrac{\partial W}{\partial \mathbf{E}_{\rm s} } ( \widehat{\mathbf{E}}_{\rm s}   )$. Finally, we define the viscous dissipation:
$$
\mathcal{D}(\widehat{u}_{\rm f})
=
2 \mu_{\rm f} \int_{\Omega_{\rm f}(t)} 
\boldsymbol{\varepsilon}_\mathrm{f}(\widehat{u}_\mathrm{f}):\boldsymbol{\varepsilon}_\mathrm{f}(\widehat{u}_\mathrm{f}) \, dx.
$$

We can now state the following result that states that the
semi-discrete ROM solution satisfies the same 
energy balance as the full-order model.

\begin{lemma}[ROM semi-discrete energy balance]\label{lem:ROM_energy}
Assume that the reduced FSI system~\eqref{eq:semi_discrete_ROM} is isolated, i.e., 
$\widehat{\mathbf u}_{\mathrm{f}}=0$ on $\partial\Omega_{\mathrm{f}}(t)\setminus\Gamma(t)$ and 
$\widetilde{\boldsymbol\Pi}_{\mathrm{s}}\widetilde{\mathbf{n}}_{\mathrm{s}}=0$ on 
$\partial\widetilde{\Omega}_{\mathrm{s}}\setminus\widetilde{\Gamma}$.
Then the reduced-order solution satisfies the energy balance
\begin{equation}
\label{eq:ROM_energy_balance}
  \frac{{\rm d}}{{\rm d}t}\bigl(
  \widehat{\mathcal T}_{\rm f}+ \widehat{\mathcal T}_{\rm s}+ \widehat{\mathcal W}_{\rm s}\bigr)
  \;+\;{\mathcal D}_{\rm f}(\widehat{u}_{\rm f}) =0 .
\end{equation}
\end{lemma}

\begin{proof}
The proof is an immediate consequence of  Galerkin orthogonality and the ALE transport formula (see, \cite[Lemma 1]{taddei2025optimization}). We omit the details.
\end{proof}

\subsection{Fully discrete ROM}
\label{sec:fully_discrete_ROM}

The reduced bases introduced in Section~\ref{sec:ROM_spaces} are purely spatial and do not depend on the IRK stage index. For an $s$-stage integrator we therefore augment them with an identity Kronecker factor:
\[
  \mathcal Z_{\rm f}:=I_s\otimes \mathbf Z_{\rm f},\quad
  \mathcal Z_{\rm s}:=I_s\otimes \mathbf Z_{\rm s},\quad
  \mathcal Z_{\rm p}:=I_s\otimes \mathbf Z_{\rm p},\quad
  \mathcal E_{\rm f}:=I_s\otimes \mathbf E_{\rm f},\quad
  \mathcal E_{\rm s}:=I_s\otimes \mathbf E_{\rm s},
\]
and denote the stage-stacked reduced coordinates by
\[
  \boldsymbol{\alpha}_{\rm f}\in\mathbb R^{s n_{\rm u}},\quad
  \boldsymbol{\alpha}_{\rm s}\in\mathbb R^{s n_{\rm s}},\quad
  \boldsymbol{\alpha}_{\rm p}\in\mathbb R^{s n_{\rm p}},\quad
  \boldsymbol{\beta}\in\mathbb R^{sm}.
\]
Recalling the stacked full-order unknowns
$(\mathbf U_{\rm f},\mathbf p_{\rm f},\mathbf D_{\rm s})$ from Section~\ref{sec:full_discrete}, we adopt the following fully discrete ROM ansatz
\begin{subequations}\label{eq:RB_ansatz}
\begin{align}
\widehat{\mathbf D}_{\rm s} &=
      \mathcal Z_{\rm s}\,\boldsymbol{\alpha}_{\rm s}
     +\mathcal E_{\rm s}\,\boldsymbol{\beta},
\\[4pt]
\widehat{\mathbf U}_{\rm f} &=
      \bigl(\mathbf 1_s\otimes\bar{ \mathbf u}_{\rm f}^{\rm dir}\bigr)
      c_{\rm dir}
     +\mathcal Z_{\rm f}\,\boldsymbol{\alpha}_{\rm f}
     +\mathcal E_{\rm f}\,
        \mathcal D_{{\rm s},\Delta t}\boldsymbol{\beta},
\\[4pt]
\widehat{\mathbf p}_{\rm f} &=
      \mathcal Z_{\rm p}\,\boldsymbol{\alpha}_{\rm p}, 
\end{align}
\end{subequations}
where $\mathcal D_{{\rm s},\Delta t}$ is the discrete time differentiation operator introduced in \eqref{eq:solid_kinematic_eqn_stacked}. Note that the port contribution in the fluid velocity and solid displacement automatically enforces the discrete kinematic condition at the interface at each IRK stage.

If we insert the ansatz \eqref{eq:RB_ansatz} into the full residual blocks of
formulation \eqref{eq:FSI_IRK_monolithic} and 
we perform Galerkin projection with the stage-extended port and bubble bases, we obtain the $s$-stage fully-discrete ROM:
\[
\left\{
\begin{aligned}
& \mathcal E_{\rm s}^{\top}\mathcal R_{\rm s}(\widehat{\mathbf D}_{\rm s})
 +\mathcal E_{\rm f}^{\top}\mathcal R_{\rm f}
      (\widehat{\mathbf U}_{\rm f},
       \widehat{\mathbf p}_{\rm f},
       \widehat{\mathbf D}_{\rm s})=0,
\\[4pt]
& \mathcal Z_{\rm f}^{\top}\mathcal R_{\rm f}(\widehat{\mathbf U}_{\rm f},
 \widehat{\mathbf p}_{\rm f},\widehat{\mathbf D}_{\rm s})=0,
\qquad
  \mathcal Z_{\rm p}^{\top}
  \mathcal B_{\rm f}
  \widehat{\mathbf U}_{\rm f}=0,
\\[4pt]
& \mathcal Z_{\rm s}^{\top}\mathcal R_{\rm s}(\widehat{\mathbf D}_{\rm s})=0.
\end{aligned}
\right.
\]
The structure is identical to that of the semi-discrete ROM \eqref{eq:semi_discrete_ROM}, but now all quantities are stage-stacked and evaluated at the IRK stage times.

To solve the fully discrete ROM, we employ the same  procedure described in Section \ref{sec:opt_SQP}, but now applied within the reduced-order setting. In this context, the control variable is the interface vector of port coordinates, $\boldsymbol{\beta}$. At each Newton iteration, the reduced interface residual and its Jacobian are computed by leveraging the HF residual operators $\mathcal R_{\rm f}$ and $\mathcal R_{\rm s}$ and the associated reduced-order bases $\mathcal Z_{\bullet}$ and $\mathcal E_{\bullet}$. Each Newton iteration thus involves the solution to the linearized reduced fluid and solid subproblems and the solution to
a reduced Schur complement system. 
This procedure closely follows the methodology described in \cite{huynh2013static,eftang2013port,ebrahimi2024hyperreduced,
craig1968coupling,seshu1997substructuring}. 
We omit the details.

\begin{remark}
The fully discrete ROM introduced above can equivalently be obtained by first deriving the semi-discrete ROM \eqref{eq:semi_discrete_ROM} and then applying the IRK time integration scheme to the reduced system, using the same Butcher tableau as for the full-order model.
\end{remark}

\begin{remark}
In this work,  the velocity basis is enriched through supremizer modes 
 to maintain the inf-sup stability of the ROM. This enrichment procedure is fundamentally different from the reduced-order basis (ROB) enrichment introduced in \cite{taddei2024non,taddei2025optimization}, which aims instead to restore full rank of the reduced Jacobian. 
As extensively discussed in the scRBE literature,  the port-bubble decomposition and Galerkin projection ensure the well-posedness of the ROM for coercive problems, while the flux-based formulation of \cite{taddei2024non,taddei2025optimization} require a specialized enrichment strategy.
 We emphasize, however, that this theoretical result strictly applies to strongly coercive cases and does not directly extend to the noncoercive (saddle-point) structure of the  incompressible FSI problems. Nevertheless, numerical experiments presented in Section \ref{sec:numerics} demonstrate that, even for the nonlinear FSI cases considered in this work, the proposed bubble-port construction remains stable and accurate without requiring the ROB enrichment of \cite{taddei2024non,taddei2025optimization}. 
\end{remark}

\section{Numerical results}

\label{sec:numerics}

\subsection{Vertical beam}

We consider an incompressible flow that interacts with a slender elastic beam oriented perperdicularly to the free-stream flow direction. The beam is clamped at its bottom end, while its upper part remains free to deform under the fluid loading. The geometric configuration can be found in \cite{taddei2025optimization}. 
No-slip wall boundary conditions are prescribed on the top and bottom boundaries of the fluid domain. At the inlet, a parabolic velocity profile is imposed with peak velocity ${u}_{\infty}=1$, while homogeneous Neumann boundary conditions are applied at the outlet.
The material parameters are chosen as
$$
\rho_{\rm f} = 1, \quad \widetilde{\rho}_{\rm s}=1.1, \quad E=10^3, \quad
\nu=0.3, \quad \mu_{\rm f}= 0.035.
$$
The coupled system is integrated until the final time $T=6\,$s. 

\subsubsection{HF results}

\begin{figure}[H]
\centering
\subfloat[$t=1$]{\includegraphics[width=0.33\textwidth]{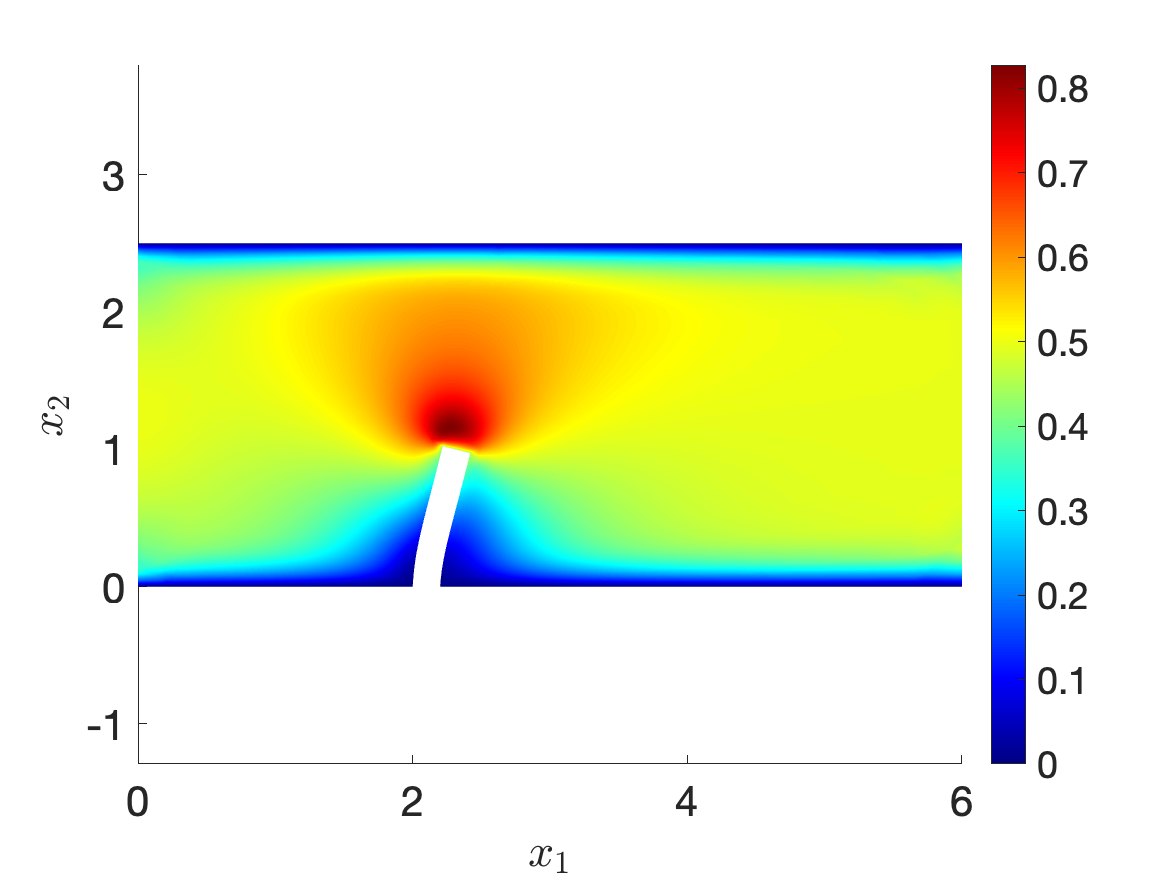}}
~~
\subfloat[$t=2$]{\includegraphics[width=0.33\textwidth]{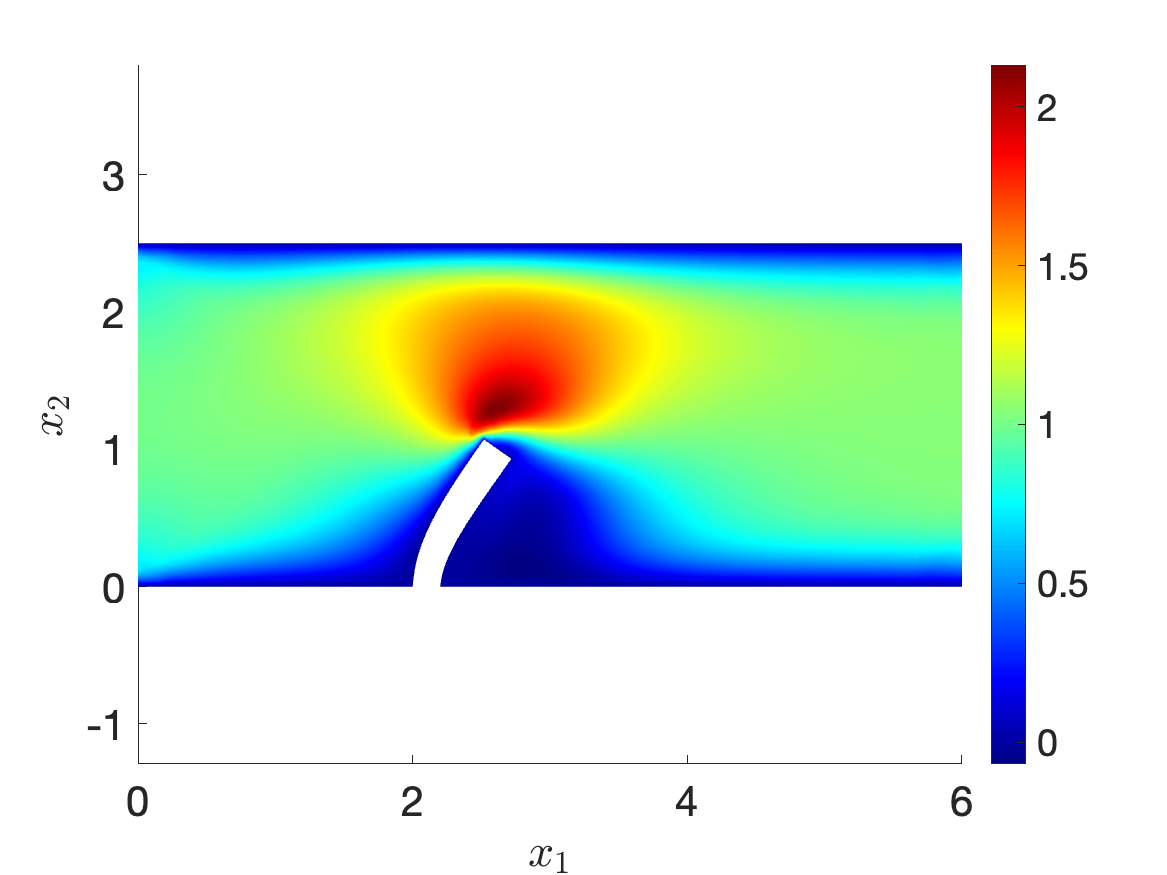}}
~~
\subfloat[$t=6$]{\includegraphics[width=0.33\textwidth]{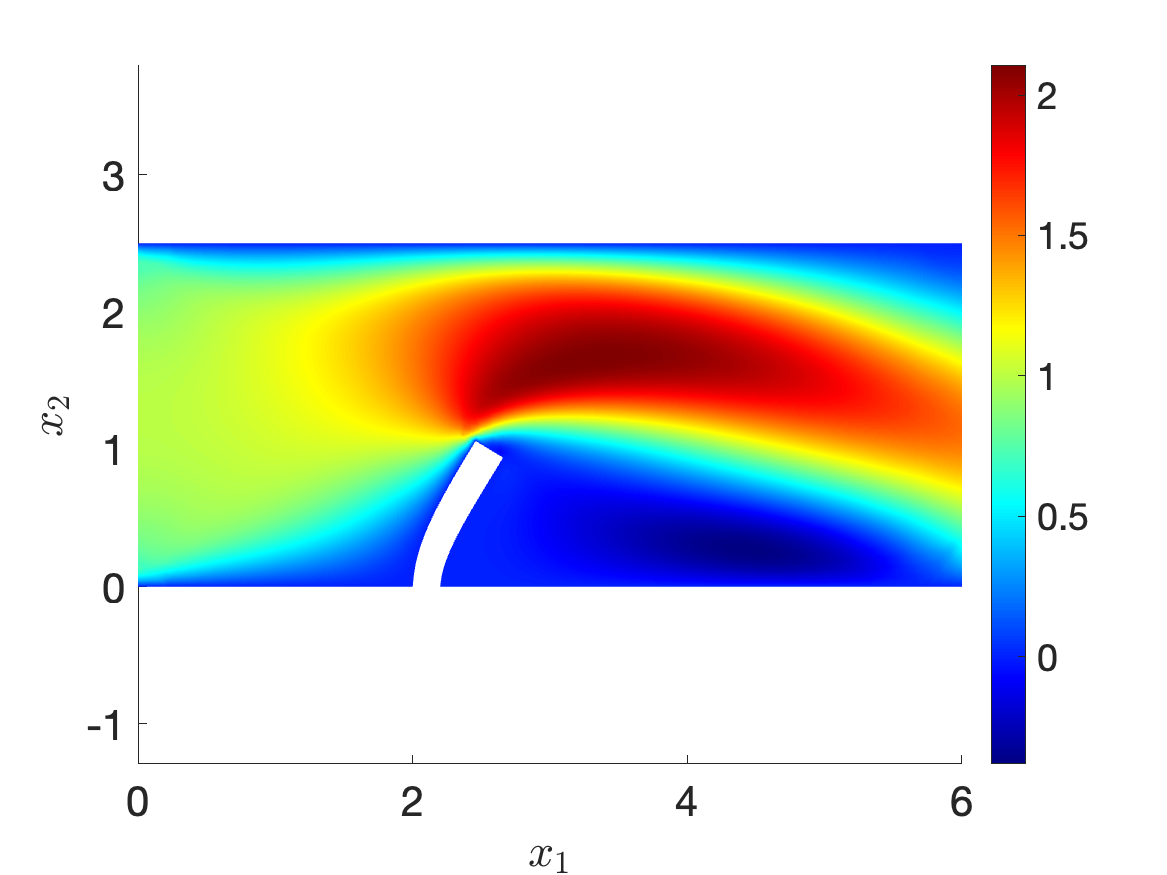}}
\caption{Vertical beam; $x$-velocity field at $3$ time instants ($s=2$, $\Delta t=0.05$).}
\label{fig:vbeam_HF_u_time126_s2}
\end{figure}

Figure \ref{fig:vbeam_HF_u_time126_s2} shows the horizontal fluid velocity field at three time instants for the Radau-IIA scheme with $s=2$ stages and time step $\Delta t=0.05$\,s. The beam is progressively bent by the flow, reaches a maximum deflection, and then recovers slightly towards its undeformed configuration.

\begin{figure}[H]
\centering
\subfloat[$F_x$]{\includegraphics[width=0.32\textwidth]{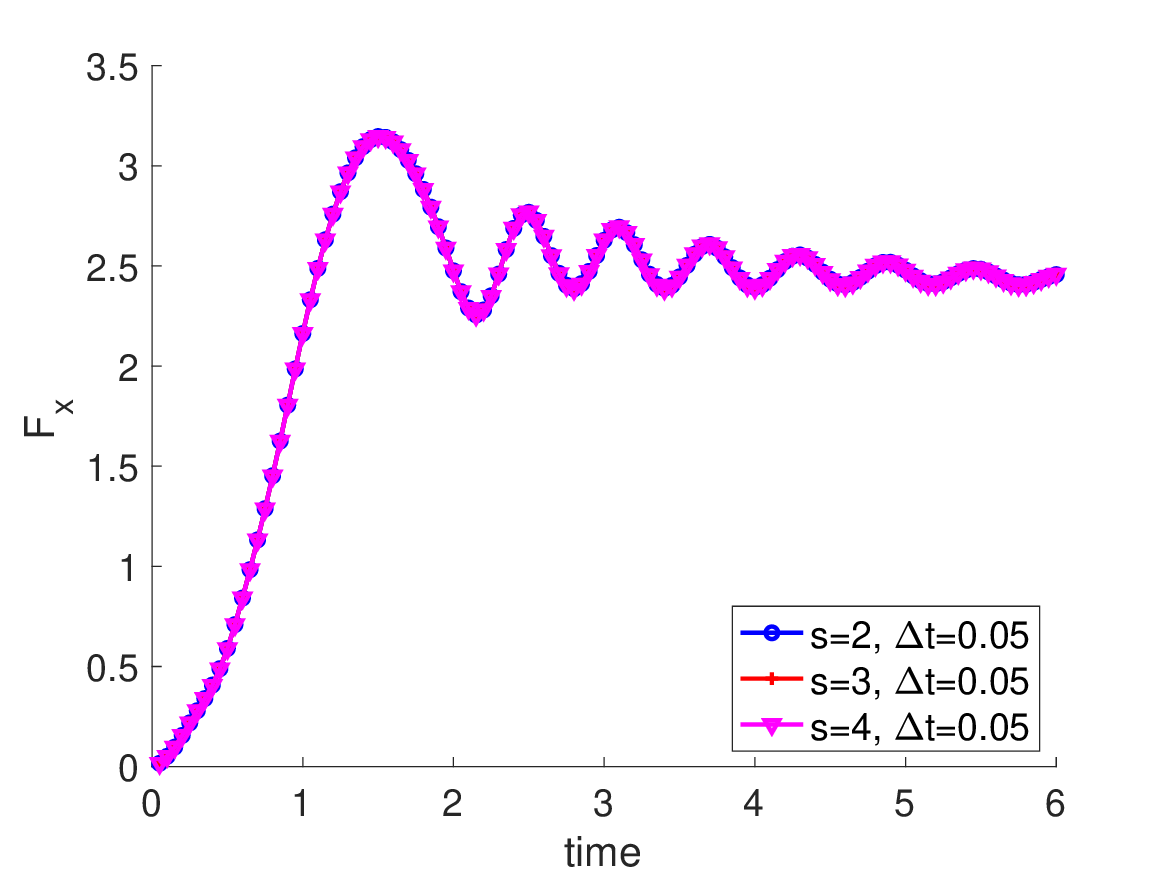}}
~~
\subfloat[$F_y$]{\includegraphics[width=0.32\textwidth]{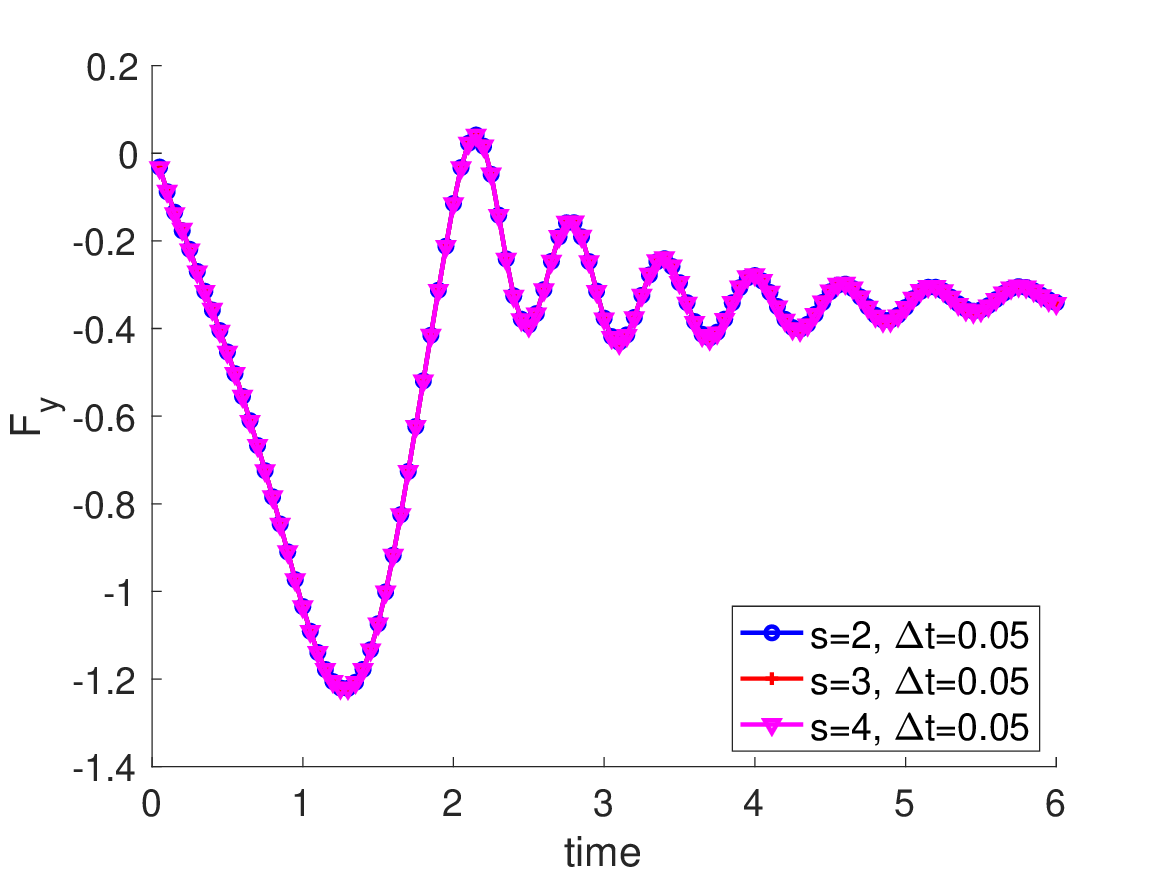}}
\subfloat[$x-$displacement]{\includegraphics[width=0.32\textwidth]{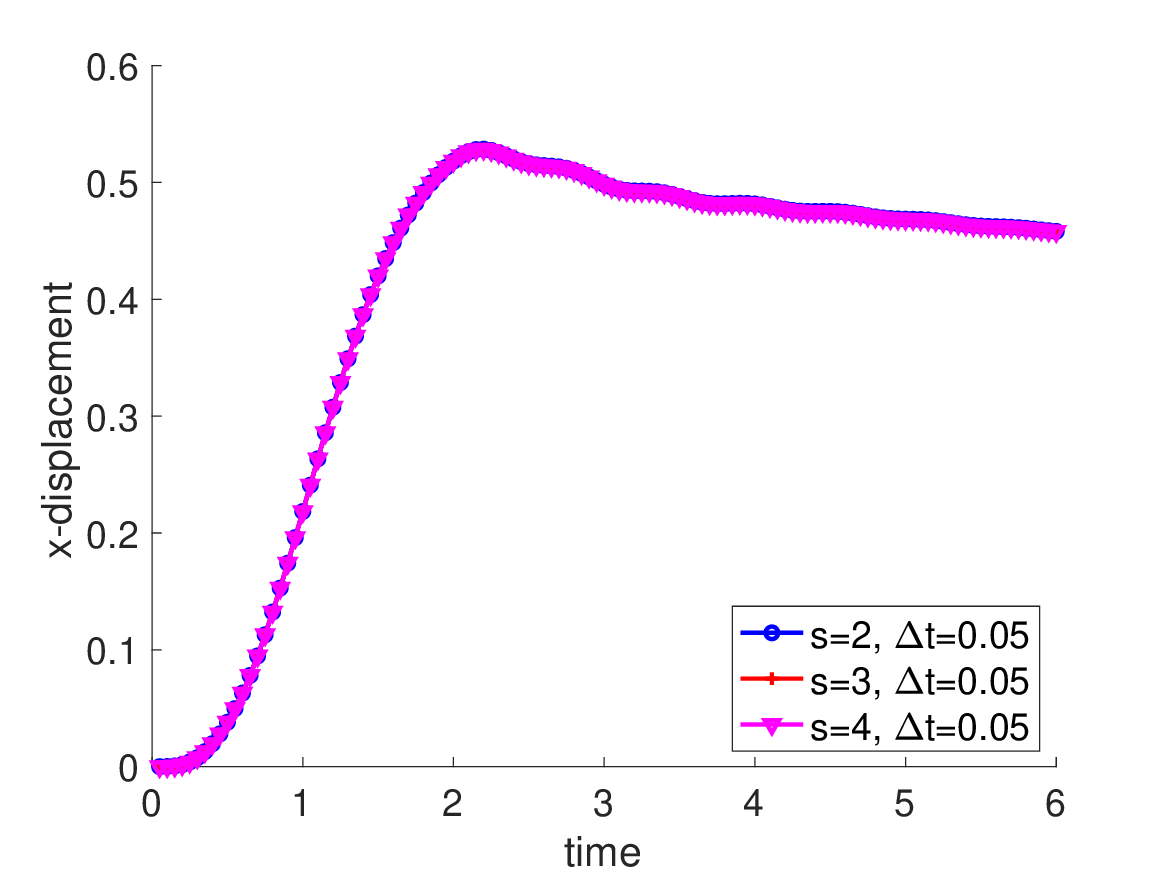}}
\caption{Vertical beam; comparison of fluid forces
and tip $x$-displacement using $s=2,3,4$ and $\Delta t=0.05$.}
\label{fig:hf_vbeam_irk357_0d05_force}
\end{figure}

We next investigate the influence of the IRK order and time step on the high-fidelity solution. We perform simulations with $s=2,3,4$ stages and time steps $\Delta t=0.025, 0.05, 0.1$\,s, and monitor the time histories of drag $F_x$, lift $F_y$, and the horizontal displacement of a probe point located at the beam tip. Results for $\Delta t=0.05$\,s are reported in Figures \ref{fig:hf_vbeam_irk357_0d05_force}. Similar results are observed for the other two time steps and all three IRK schemes produce closely matching results. 
Table \ref{tab:vbeam_H1_disp_error} reports the $H^1$-norm error in the solid displacement for Radau-IIA schemes. 
Here, the solution obtained with $s=4$ and time step 
$\Delta t=0.001$ is used as the reference solution.
For each method, the error increases monotonically with $\Delta t$, and  an accuracy gain is observed when increasing the number of stages. 
By performing a least-squares fit to the error data, 
we estimate effective temporal convergence rates of approximately
$p \approx 2.1$, $p \approx 2.4$, and $p \approx 3.5$ for $s=2,3,4$, respectively. 
The results confirm high-order convergence in time and a
systematic accuracy improvement as $s$ increases. The observed convergence rates remain below the formal Radau-IIA orders $2s-1$: this   can be attributed to the differential algebraic structure of the  FSI system, where temporal errors interact with spatial discretization of the fluid and solid subproblems, the ALE mesh motion, and the iterative interface-coupling.

\begin{table}[H]
  \centering
  \begin{tabular}{c|c|c|c}
    \toprule
    \diagbox{$s$}{$\Delta t$} & 0.025 & 0.05 & 0.10 \\
    \midrule
    2 & $7.09\times 10^{-5}$ & $2.22\times 10^{-4}$ & $1.32\times 10^{-3}$  \\
    3 & $3.81\times 10^{-6}$ & $6.55\times 10^{-5}$ & $1.04\times 10^{-4}$  \\
    4 & $8.52\times 10^{-7}$ & $6.38\times 10^{-6}$ & $1.07\times 10^{-4}$  \\
    \bottomrule
  \end{tabular}
  \caption{Vertical beam; $H^1$-norm of the solid displacement error   for Radau-IIA schemes with $s=2,3,4$ and time steps
  $\Delta t = 0.20, 0.10, 0.05$.}
  \label{tab:vbeam_H1_disp_error}
\end{table}

\subsubsection{ROM results}

\begin{figure}[H]
  \centering
\subfloat[$t=1$]{\includegraphics[width=0.33\textwidth]{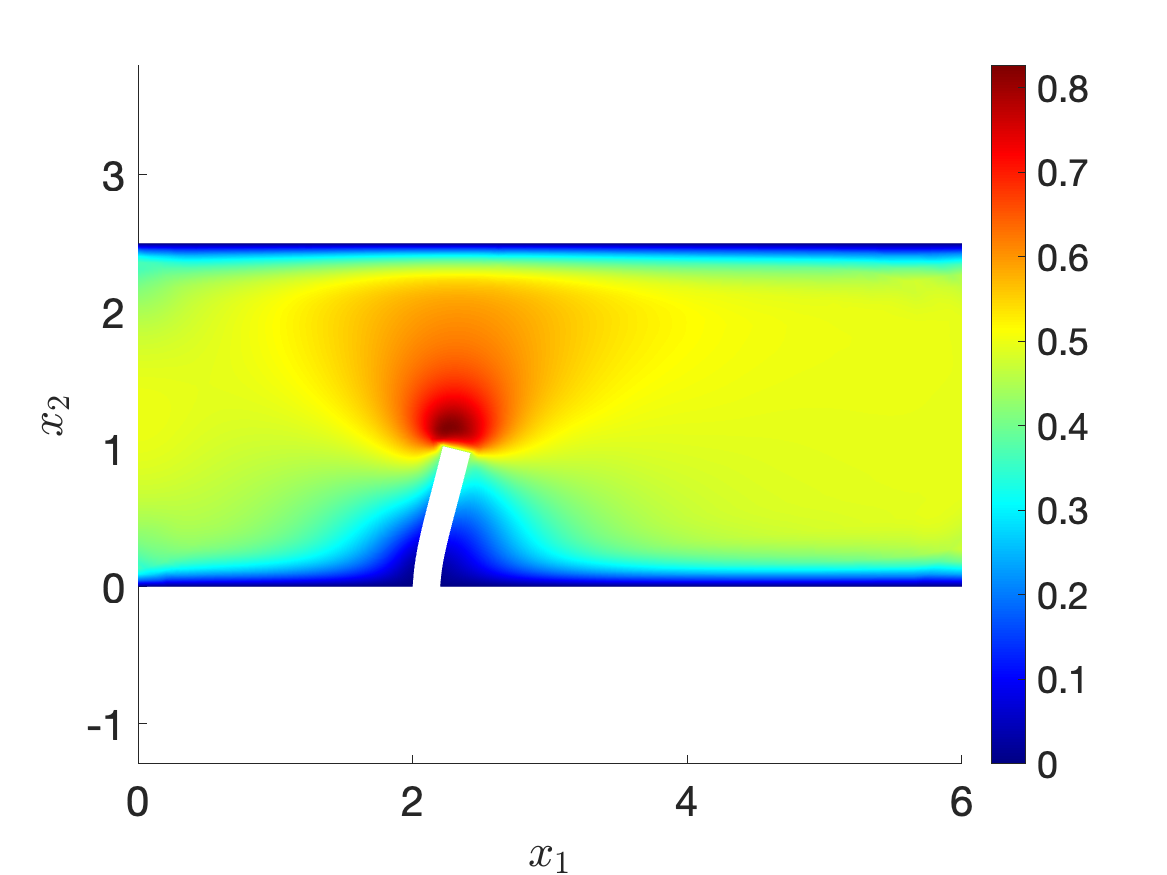}}
~~
\subfloat[$t=2$]{\includegraphics[width=0.33\textwidth]{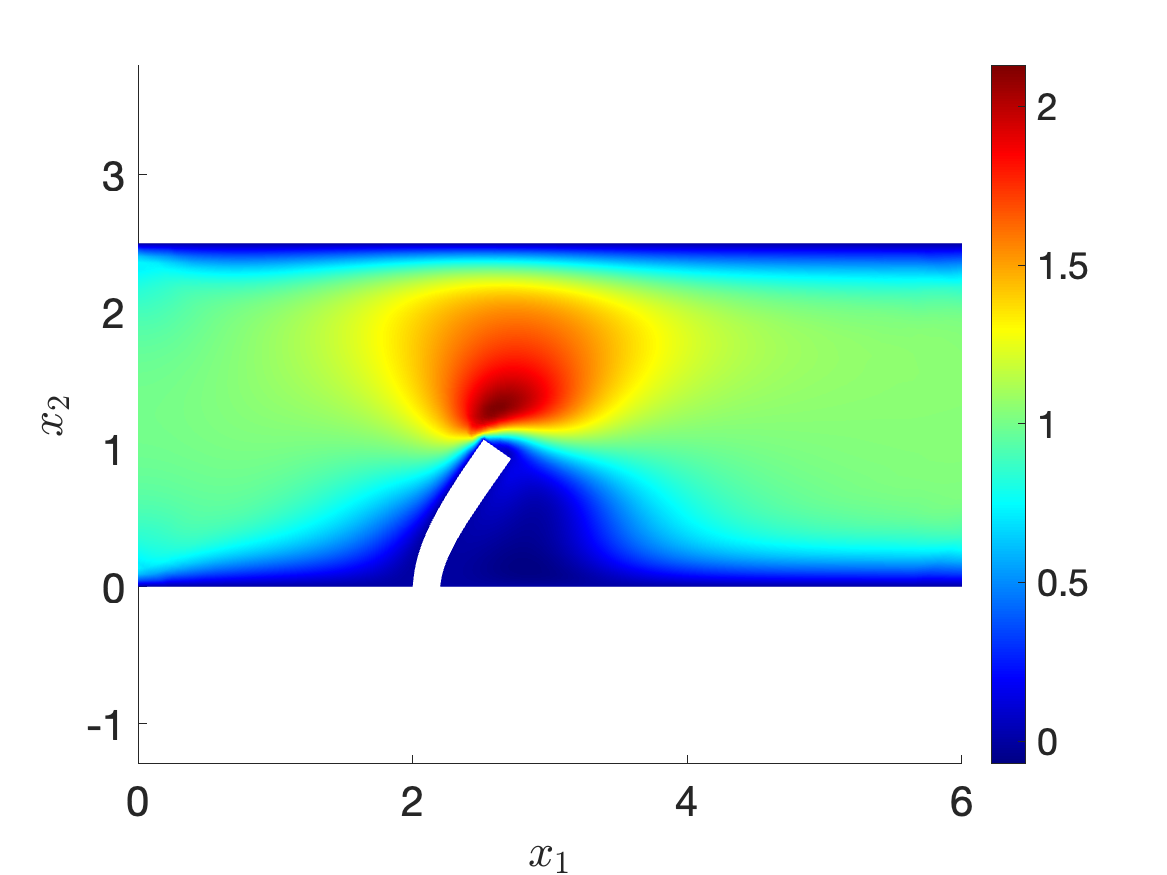}}
~~
\subfloat[$t=6$]{\includegraphics[width=0.33\textwidth]{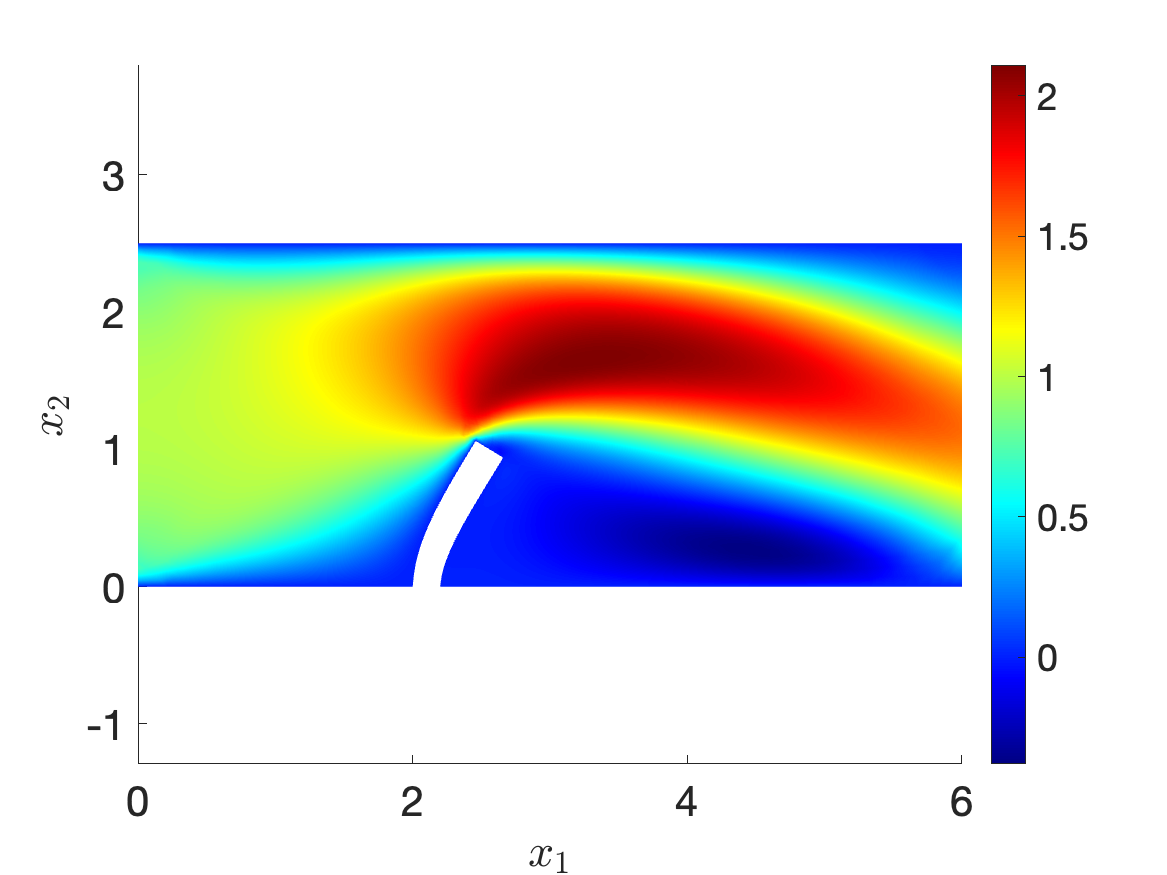}}
\\
\vspace*{-0.4cm}
\subfloat[$t=1$ (error)]{\includegraphics[width=0.33\textwidth]{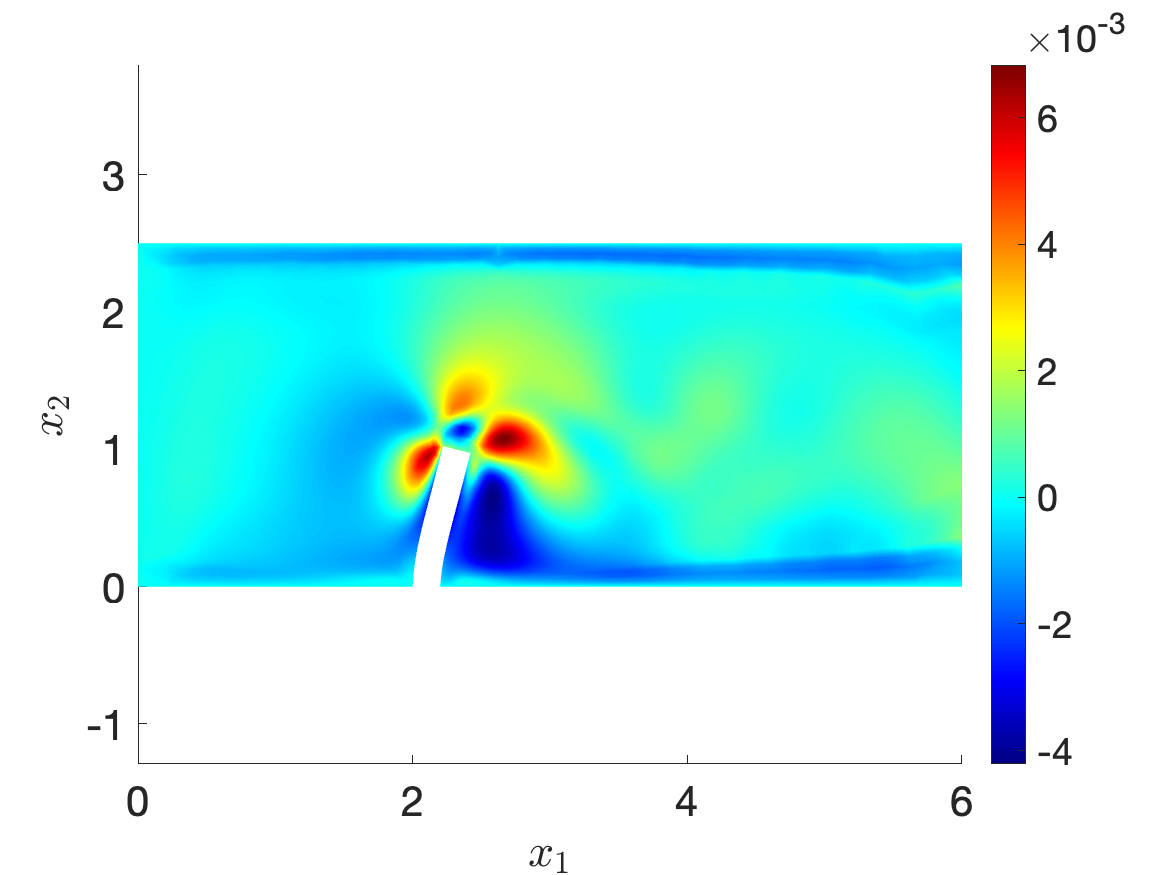}}
~~
\subfloat[$t=2$ (error)]{\includegraphics[width=0.33\textwidth]{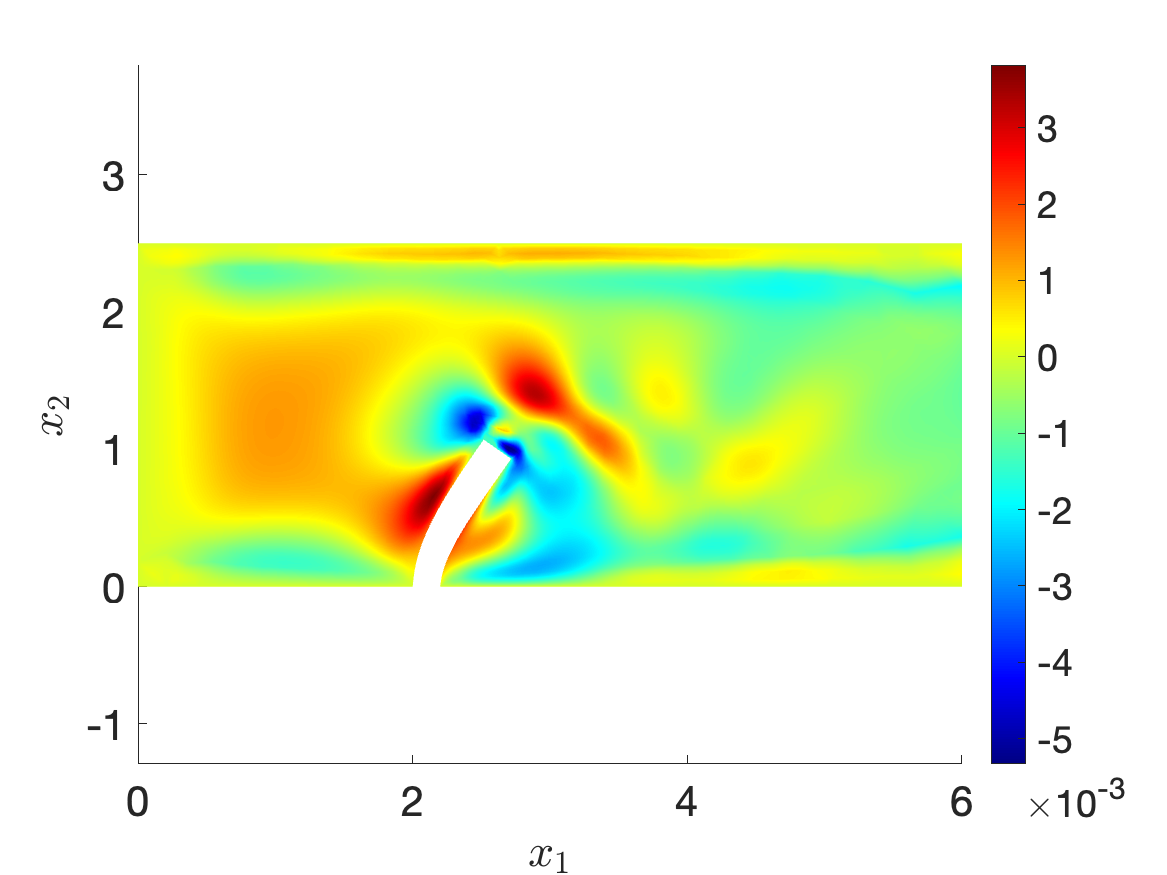}}
~~
\subfloat[$t=6$ (error)]{\includegraphics[width=0.33\textwidth]{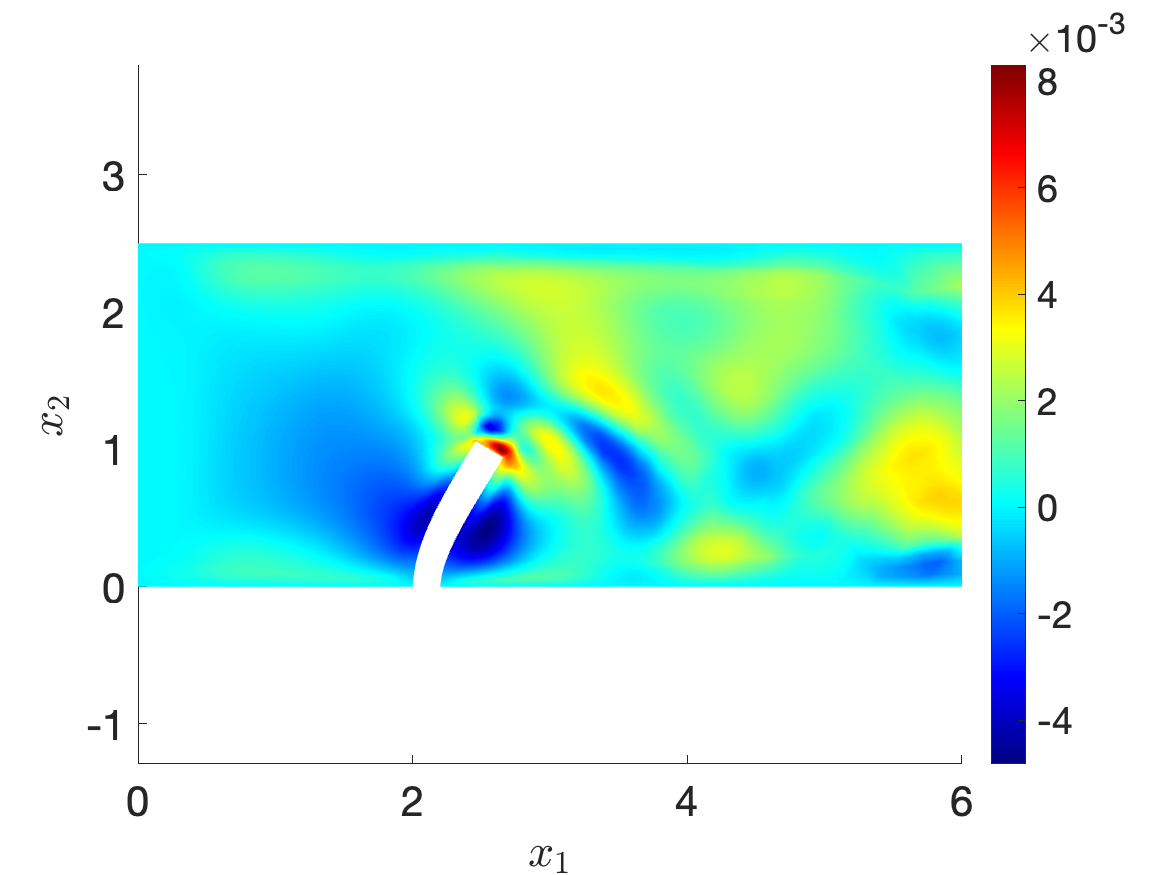}}
\caption{Vertical beam; ROM results obtained using IRK ($s=2,\Delta t=0.05\,$s, $\rm{tol}_{\rm POD} = 10^{-5}$).}
\label{fig:vbeam_ROM_irk3_stab}
\end{figure}

We first consider ROMs built from $s=2$ IRK snapshots with time step $\Delta t=0.05$\,s and POD truncation tolerance $\mathrm{tol}_{\rm POD}=10^{-5}$. Figure \ref{fig:vbeam_ROM_irk3_stab} shows the ROM horizontal fluid velocity field and the corresponding pointwise error at $t=1,2,6$\,s. The ROM accurately captures the main flow and beam deformation features, and the errors remain small throughout the domain.

\begin{figure}[H]
\centering
\subfloat[ROM errors]{\includegraphics[width=0.45\textwidth]{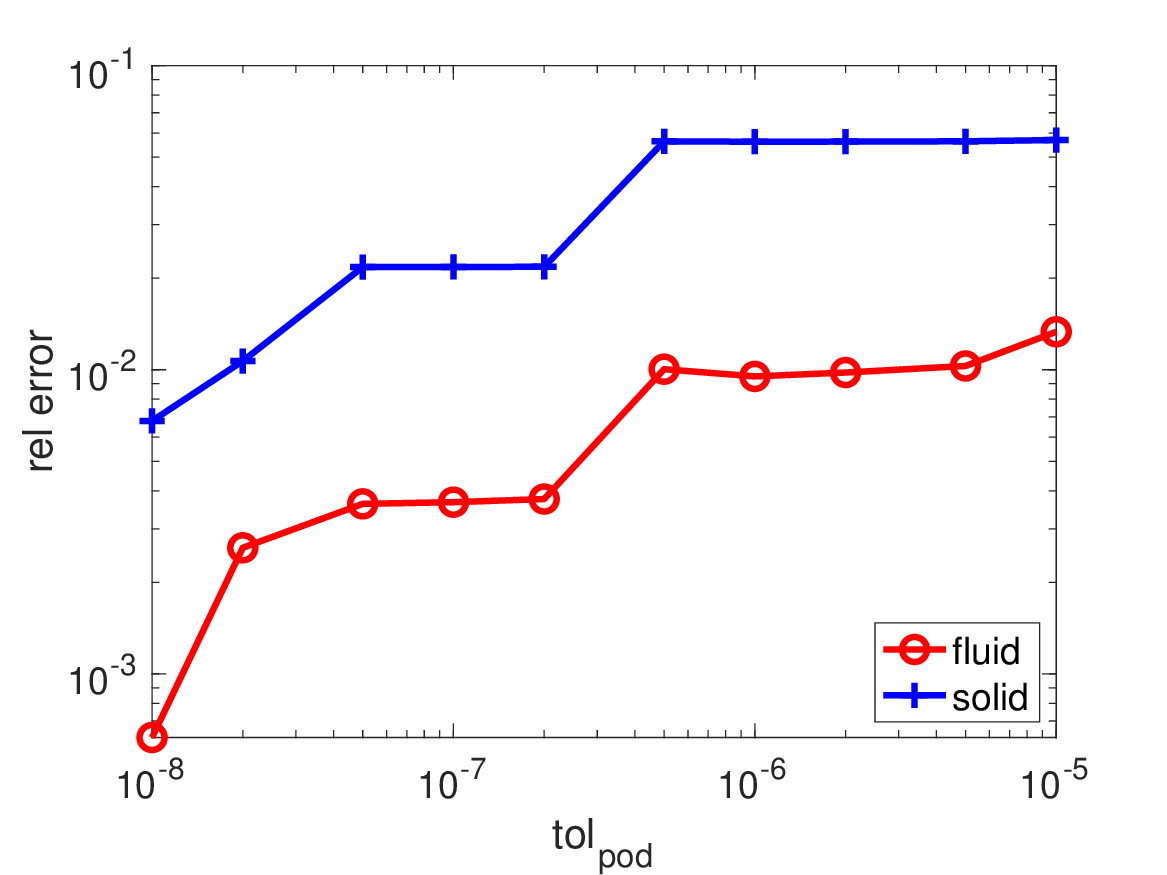}}
~~
\subfloat[number of modes]{\includegraphics[width=0.45\textwidth]{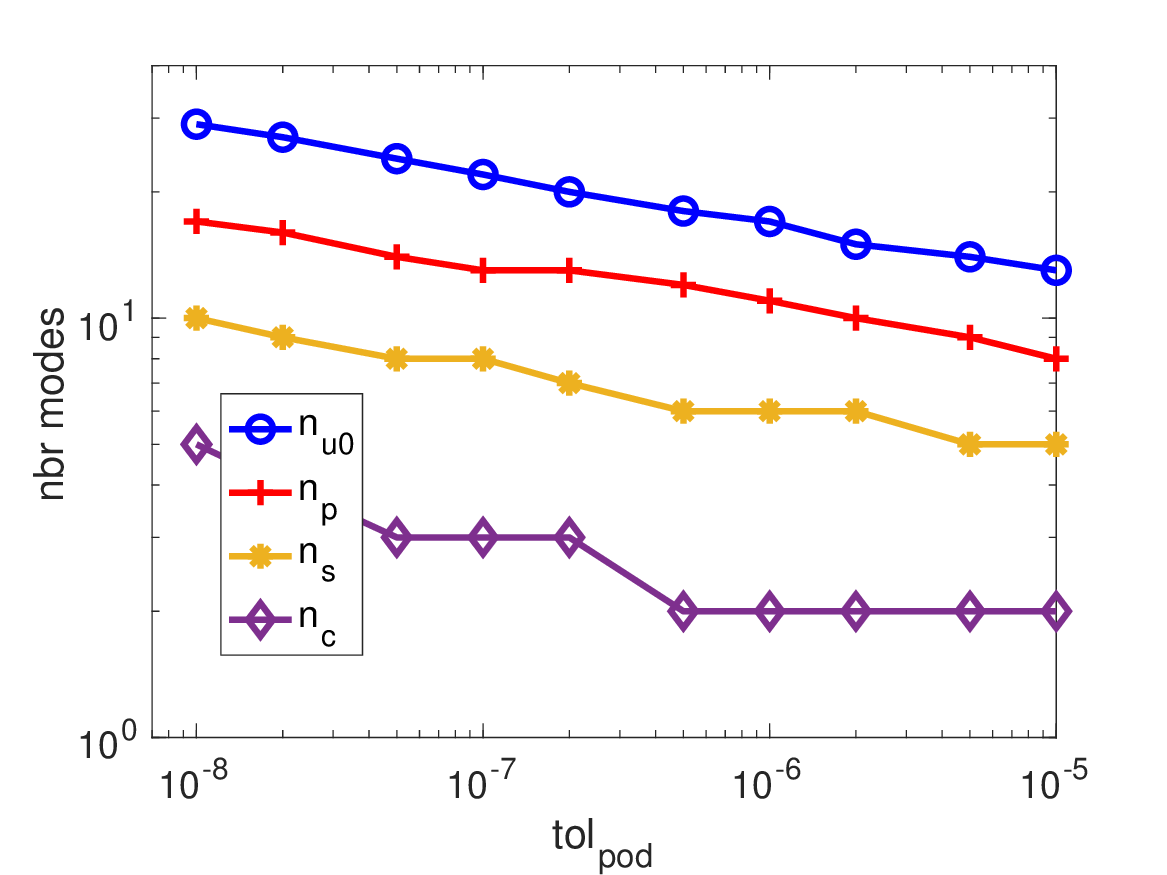}}
\caption{Vertical beam; ROM errors and number of modes in terms of POD tolerance  ($s=2,\Delta t=0.05\,$s).}
\label{fig:vbeam_ROM_err_s2_0d05}
\end{figure}

Figure \ref{fig:vbeam_ROM_err_s2_0d05} reports the $H^1\times L^2$ relative error of the velocity-pressure pair as a function of $\mathrm{tol}_{\rm POD}$, together with the corresponding number of reduced modes. As expected, tighter tolerances lead to smaller errors at the price of larger reduced spaces.

\begin{figure}[H]
\centering
\subfloat[$F_x$]{\includegraphics[width=0.33\textwidth]{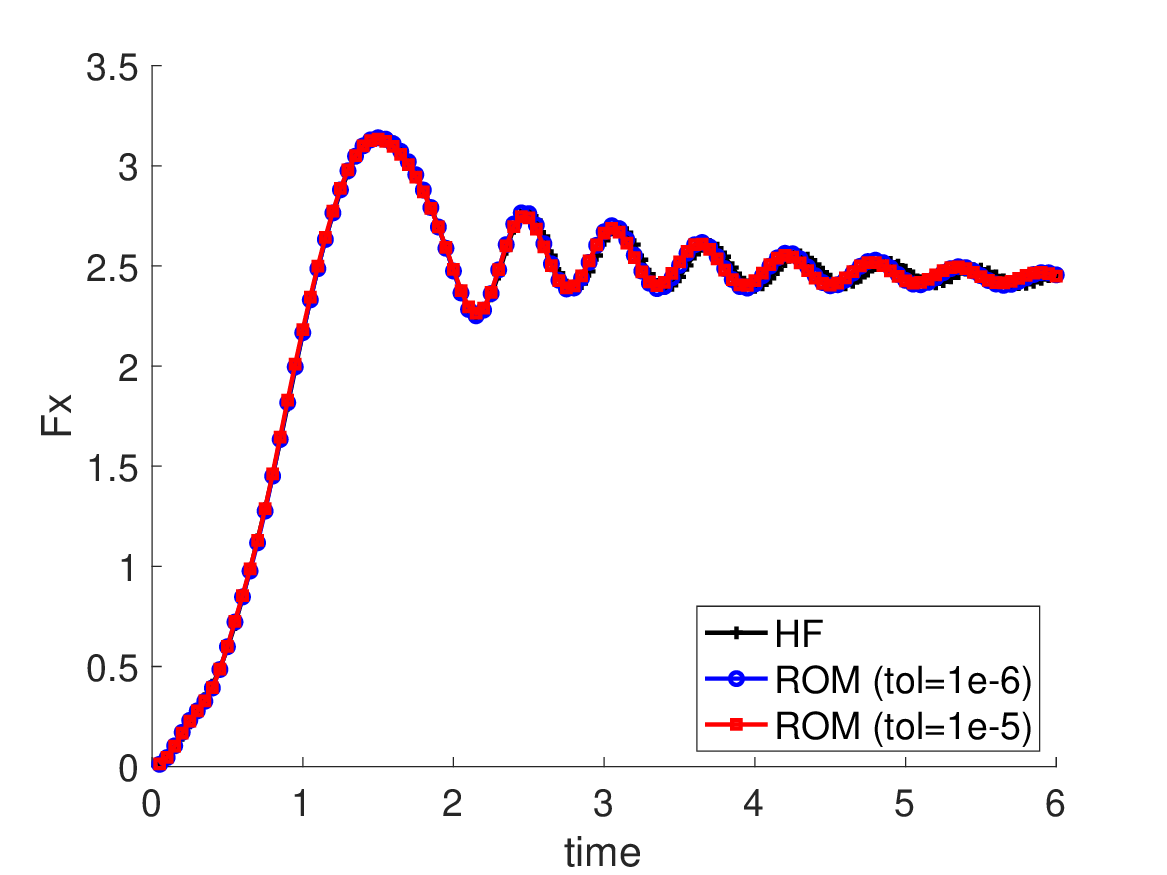}}
~~
\subfloat[$F_y$]{\includegraphics[width=0.33\textwidth]{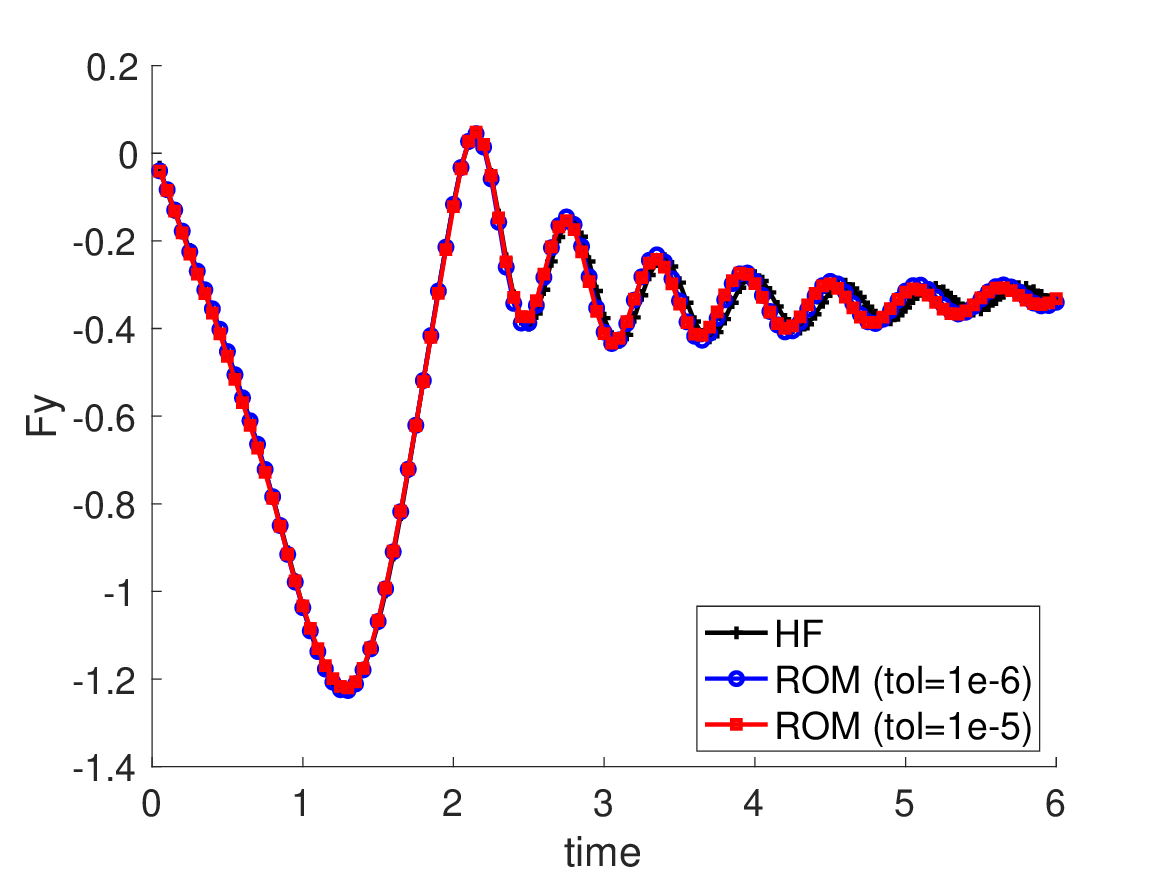}}
~~
\subfloat[x-displacement]{\includegraphics[width=0.33\textwidth]{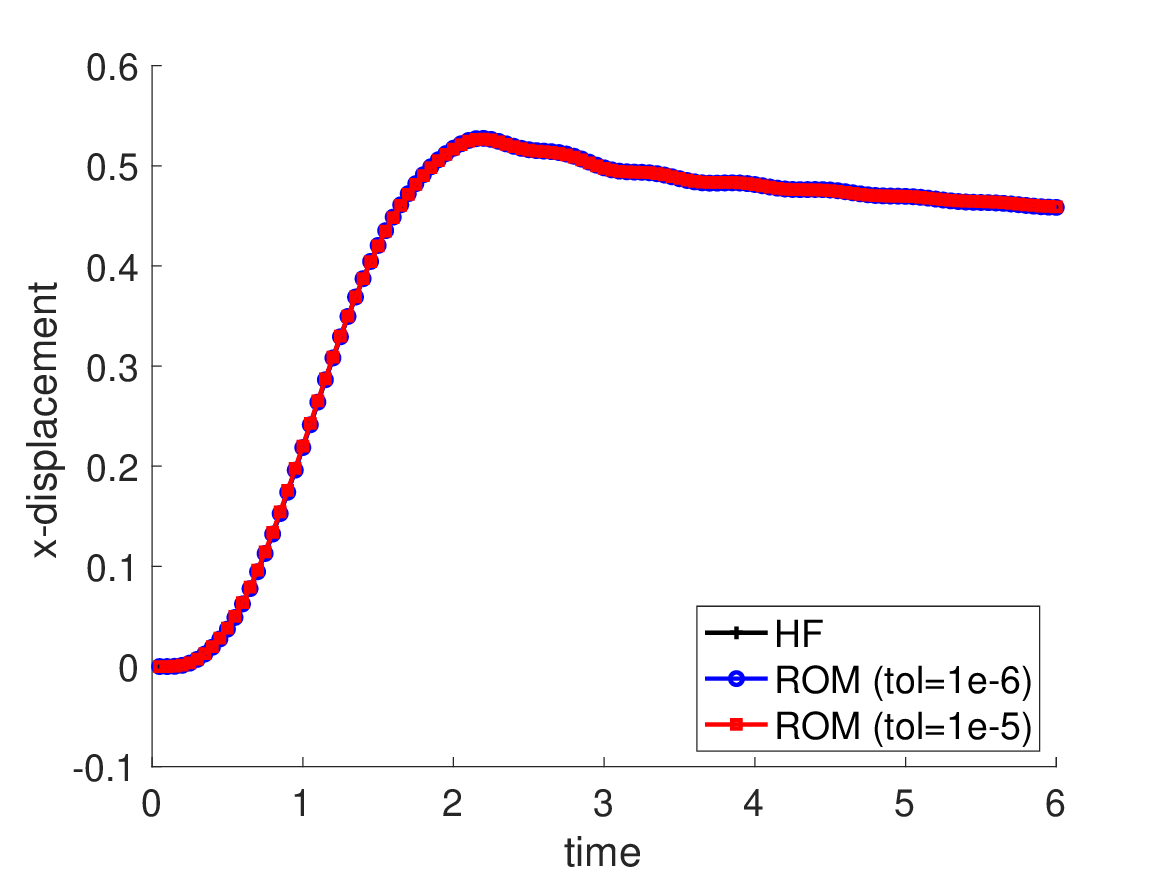}}
\caption{Vertical beam; comparison between HF and two ROM results ($s=2,\Delta t=0.05\,$s).}
\label{fig:vbeam_compr_s2_0d05}
\end{figure}

In Figure \ref{fig:vbeam_compr_s2_0d05}, we compare high-fidelity and ROM predictions for drag, lift, and tip $x$-displacement for two representative ROMs built with $\mathrm{tol}_{\rm POD}=10^{-6}$ and $10^{-5}$. Both ROMs accurately reproduce the force signals and the $x$-displacement.

We also repeat the ROM analysis for IRK with $s=3$ and $\Delta t=0.1$\,s. The corresponding velocity fields, error plots, and error versus tolerance curves confirm the same qualitative behavior observed for $s=2$: ROMs remain stable and accurate, and the use of larger  POD tolerances reduces the number of retained POD modes at the cost of increased approximation errors. These additional results are reported in \ref{app:vbeam_rom_s3}.

\begin{figure}[H]
\centering
\subfloat[$F_x$]{\includegraphics[width=0.33\textwidth]{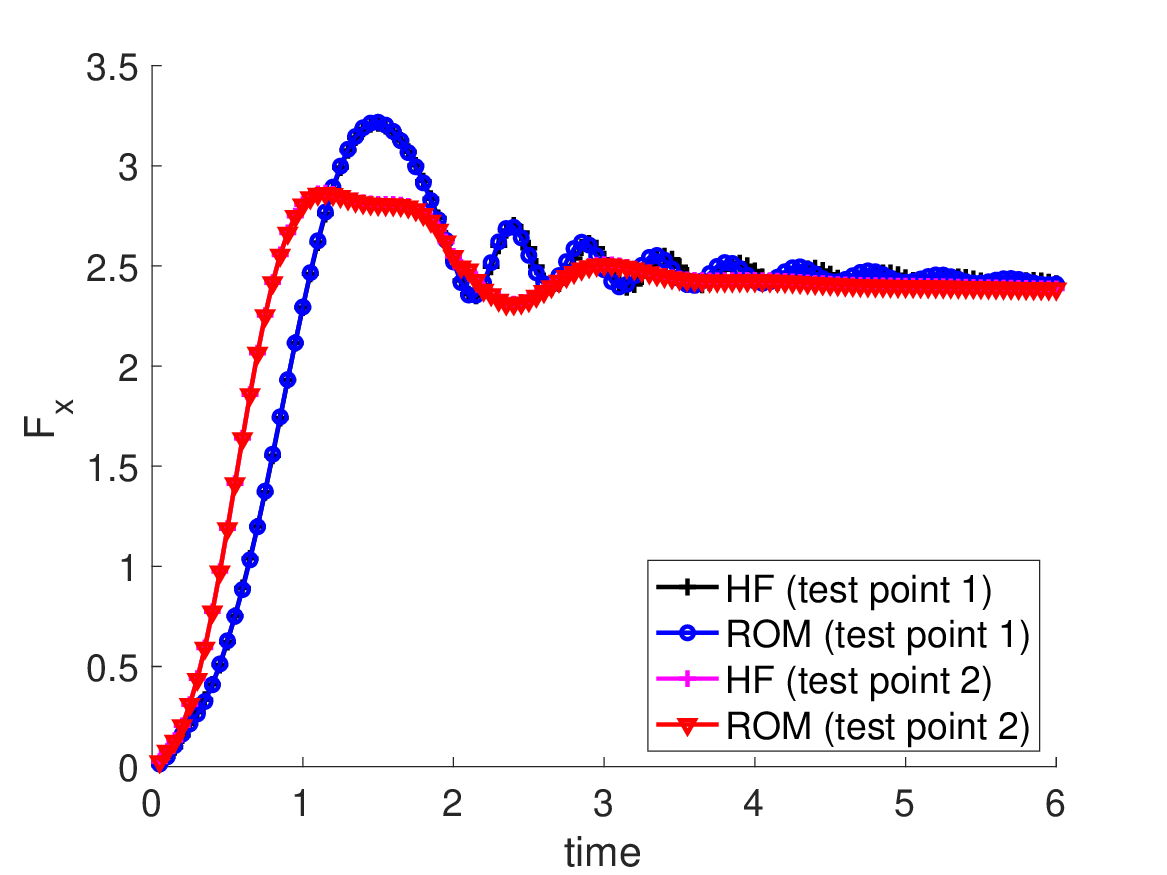}}
~~
\subfloat[$F_y$]{\includegraphics[width=0.33\textwidth]{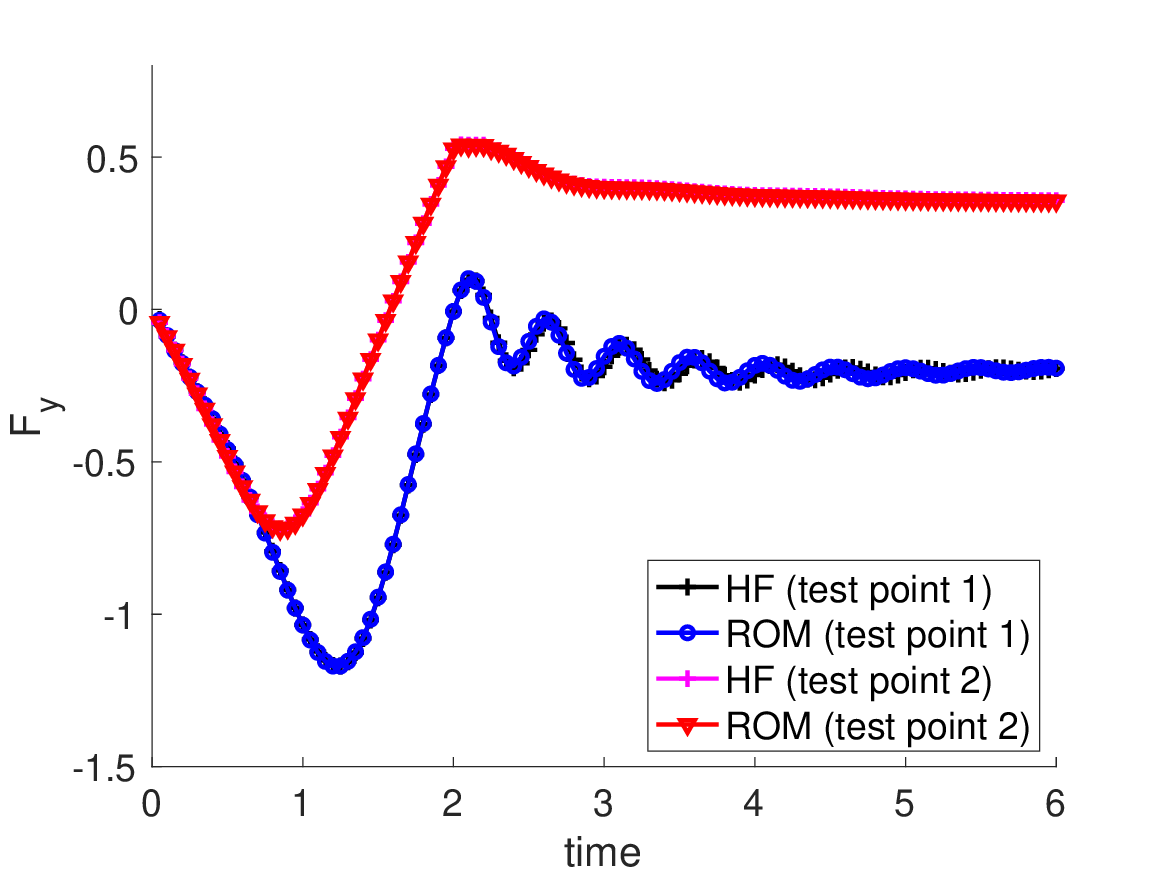}}
~~
\subfloat[x-displacement]{\includegraphics[width=0.33\textwidth]{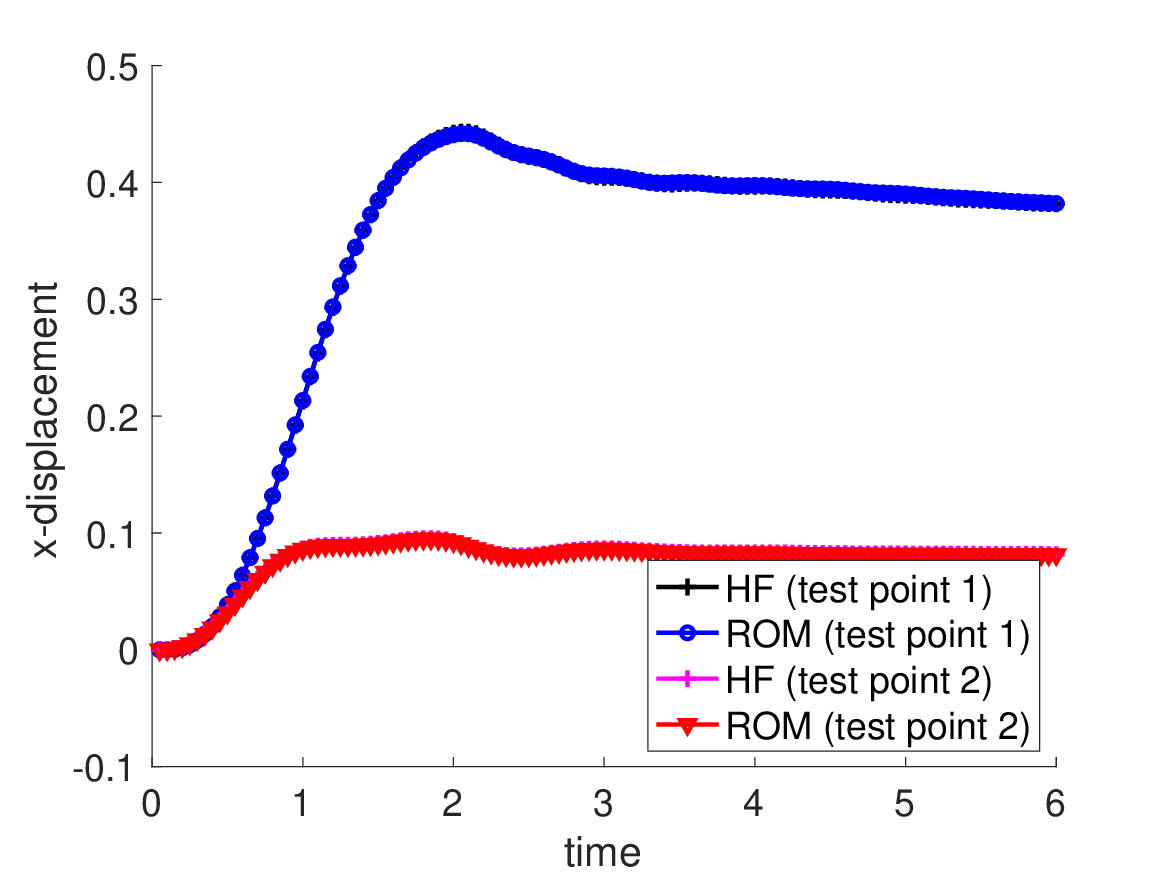}}
\caption{Vertical beam; drag and lift forces, and tip $x$-displacement at two representative test points.}
\label{fig:vbeam_param}
\end{figure}

\textbf{Parametric case.} We consider the Young's modulus $E_{\rm s}$ and the mass ratio $m^{*}:=\frac{\widetilde{\rho}_{\rm s}}{\rho_{\rm f}}$ as parameters, with ranges $E_{\rm s}/10^3\in [1, 10]$ and $m^{*}\in [0.5,2]$, while keeping $\rho_{\rm f}$ fixed. We consider $s=2$ stages and time step $\Delta t=0.05$ for the IRK time integrator. 
To build and assess the ROM, we generate $20$ training points and $5$ additional test points. In Figure \ref{fig:vbeam_param}, we show predictions of the lift and drag forces, as well as the beam tip displacement in the $x$-direction, obtained using two different ROM tolerances. The predictions demonstrate good accuracy and robustness of the reduced-order model across the considered parametric ranges.

\subsection{Turek}

We now consider the widely-used Turek FSI3  benchmark \cite{TurekHron2006}, which consists of an incompressible flow past a cylinder with an attached elastic beam. The cylinder is rigid and stationary, and the elastic beam is fixed at its left edge (attached to the cylinder) while the remainder is free to deform under fluid-induced forces.   
We follow the standard geometric configuration and parameter choices described in 
\cite{TurekHron2006,taddei2025optimization}. 
We initialize the FSI simulation by first holding the beam fixed in its undeformed configuration for $T_0=30\,$s to allow the vortex shedding pattern to develop. Then, at $t=T_0$,  we release the beam and integrate the coupled FSI system forward in time.

\subsubsection{HF results}

\begin{figure}[H]
\centering
\subfloat[$t=32$]{\includegraphics[width=0.33\textwidth]{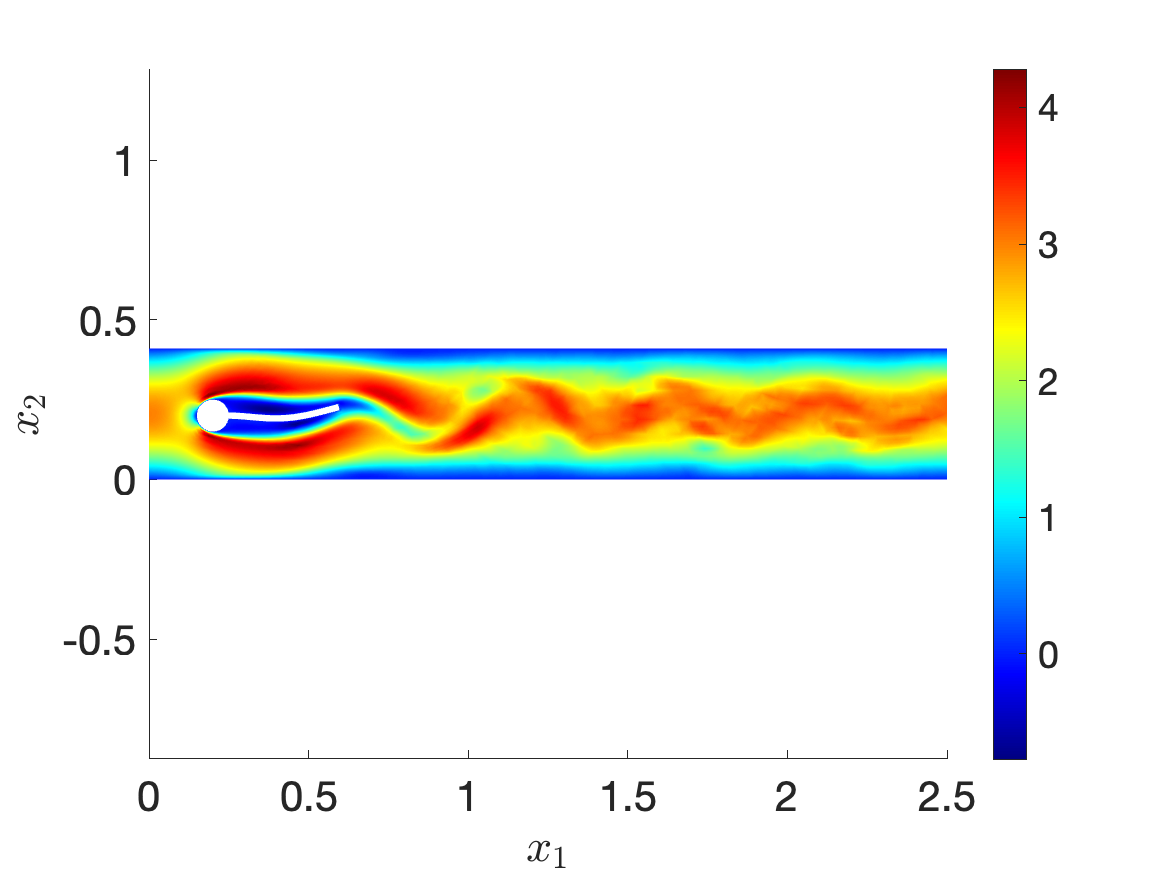}}
~~
\subfloat[$t=35$]{\includegraphics[width=0.33\textwidth]{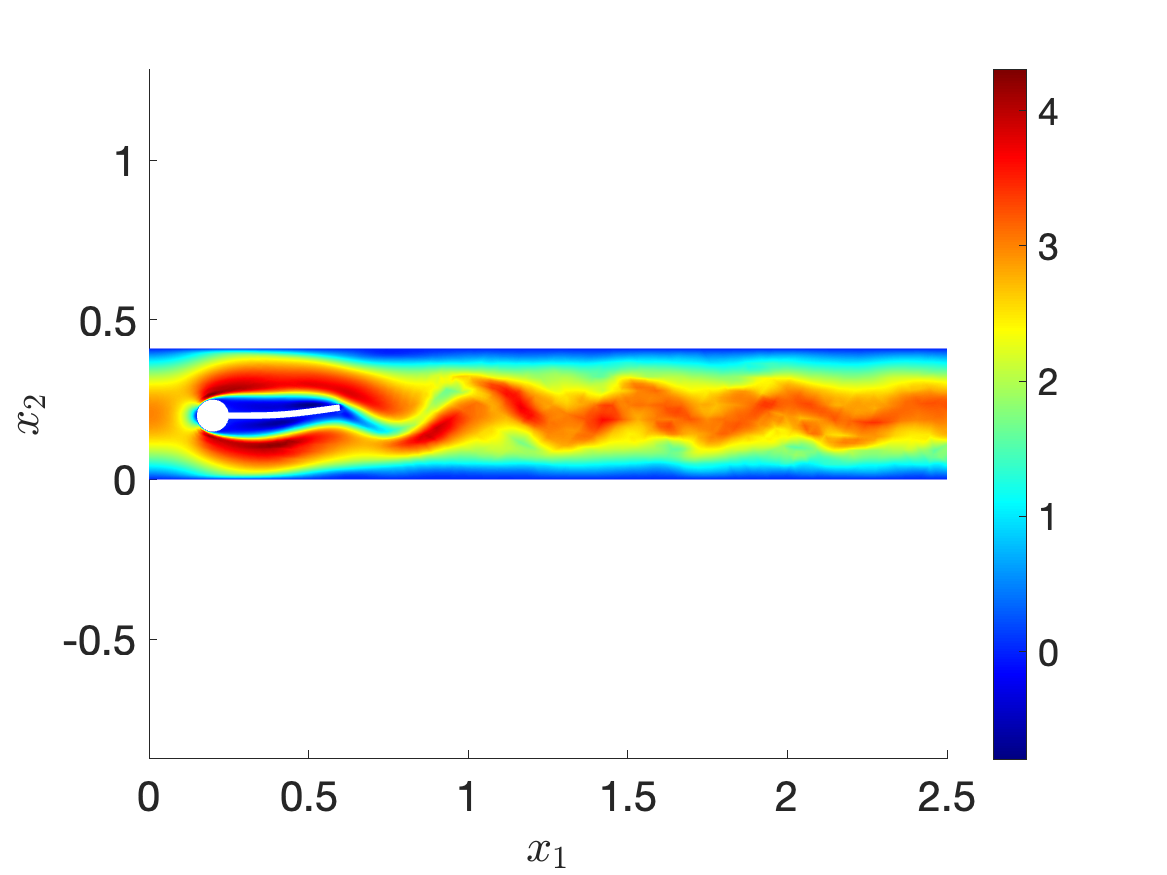}}
~~
\subfloat[$t=38$]{\includegraphics[width=0.33\textwidth]{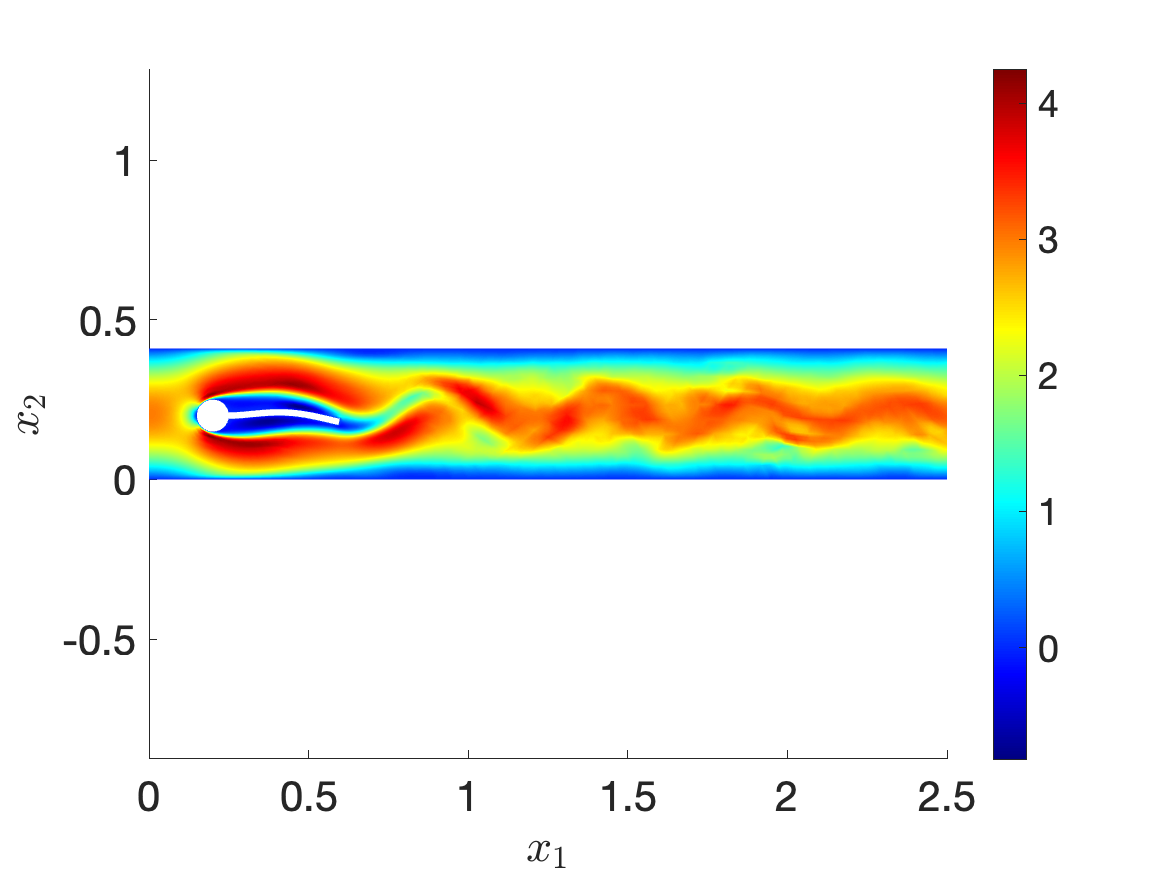}}
\caption{Turek; x-velocity field at 3 time instants.}
\label{fig:turek_HF_u_time258_s2}
\end{figure}

Figure \ref{fig:turek_HF_u_time258_s2} shows the horizontal velocity field at three time instants  $t=32,35,38\,$s obtained using the Radau-IIA IRK scheme with $s=2$ stages ($3$rd-order accuracy, hereafter referred as IRK3) and time step $\Delta t=0.01\,$s. The plots illustrate the formation of a periodic vortex street downstream of the cylinder and the associated oscillatory motion of the beam.

\begin{figure}[H]
\centering
\subfloat[$F_x$]{\includegraphics[width=0.33\textwidth]{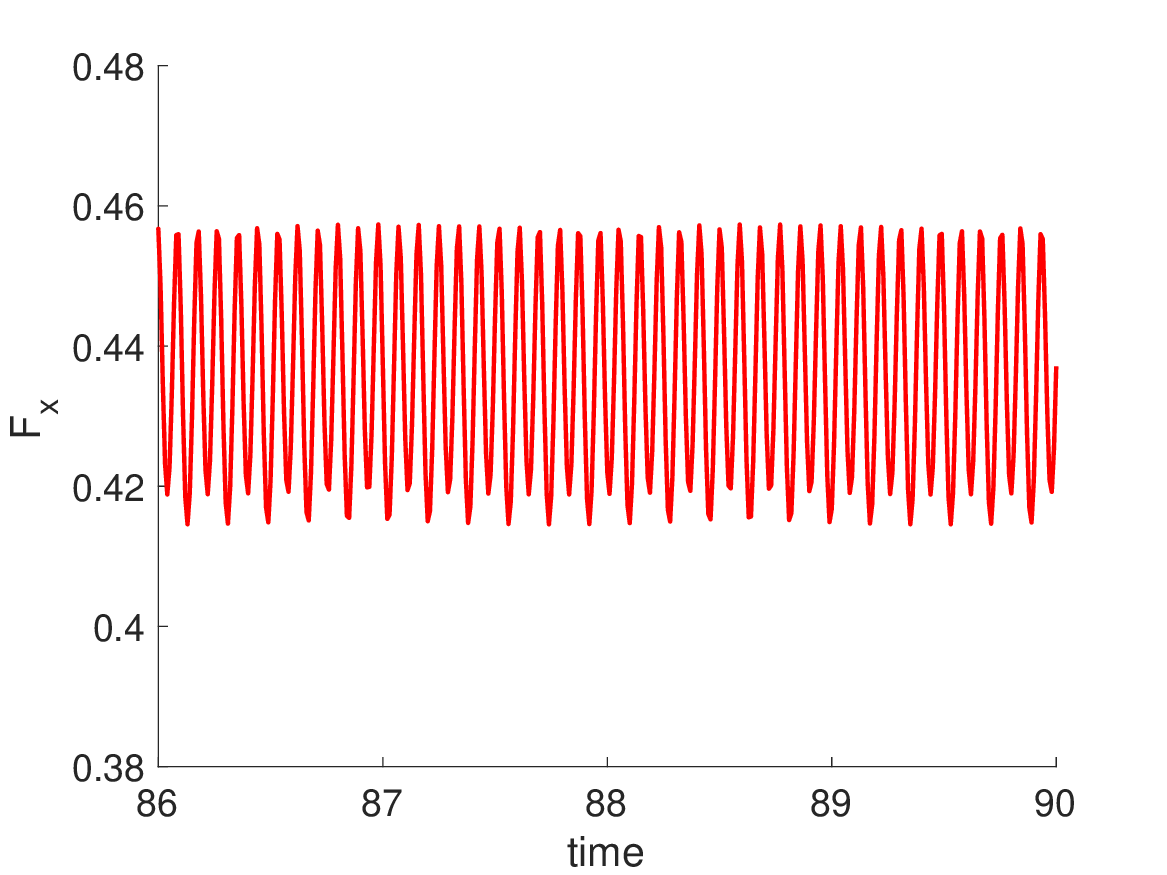}}
~~
\subfloat[$F_y$]{\includegraphics[width=0.33\textwidth]{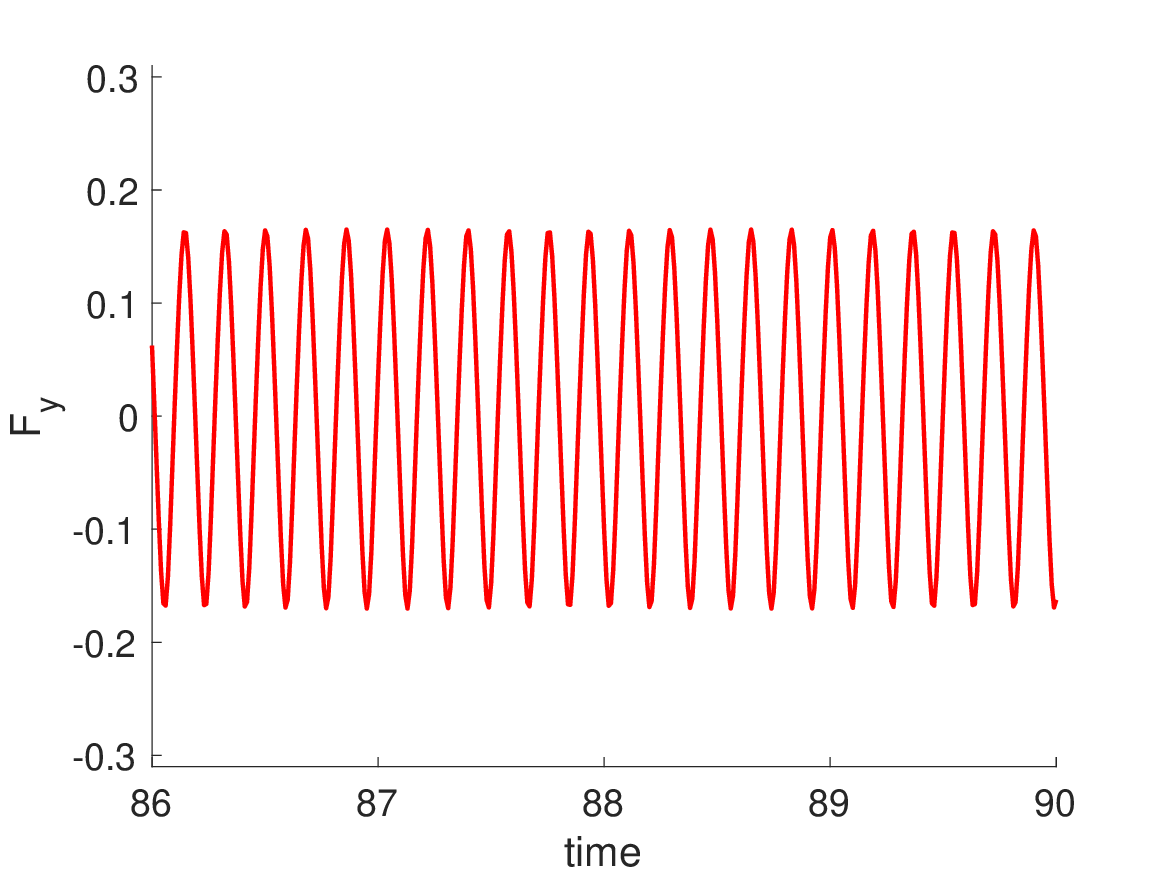}}
~~
\subfloat[Total energy]{\includegraphics[width=0.33\textwidth]{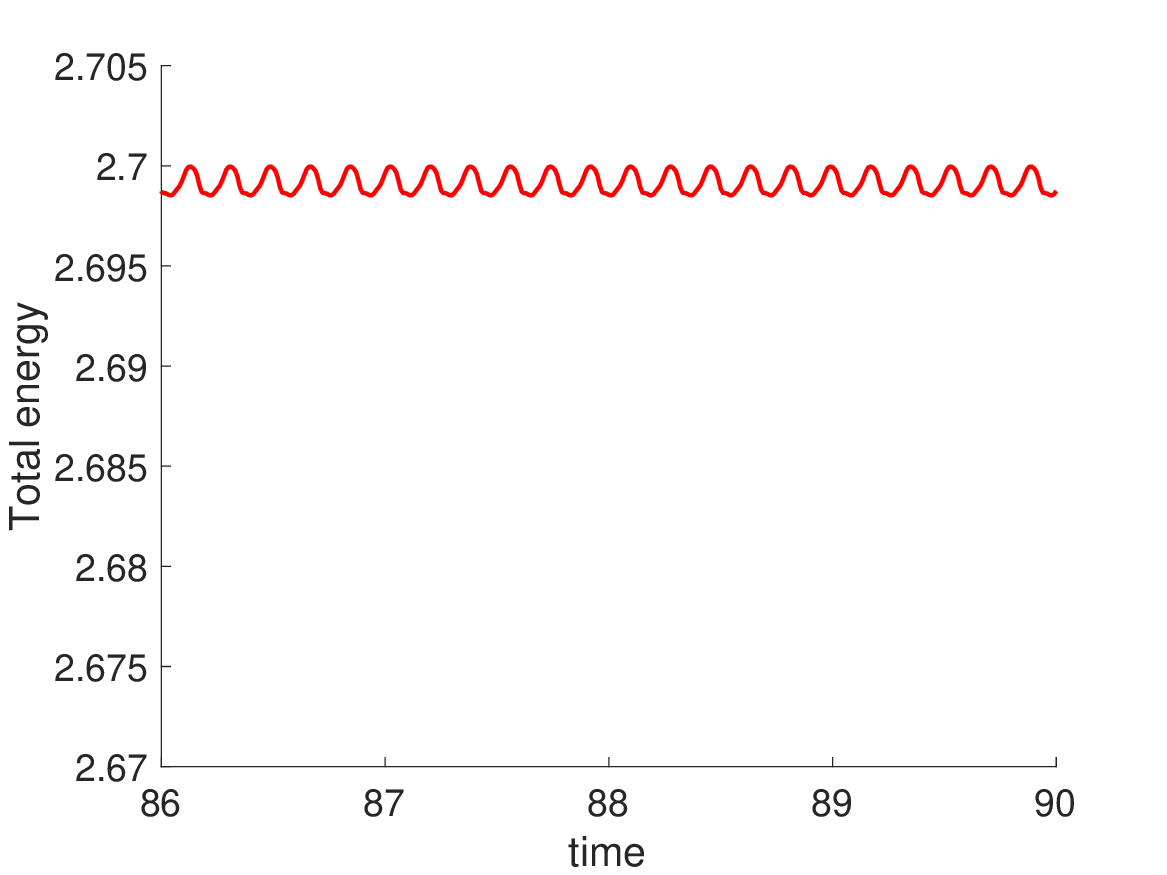}}
\caption{Turek; stable long-time HF results obtained using IRK3.}
\label{fig:turek_HF_irk3_stab}
\end{figure}

In Figure \ref{fig:turek_HF_irk3_stab} we report the long-time histories of the drag and lift forces acting on the cylinder-beam system, together with the total energy for IRK3. After an initial transient, all three signals converge to a periodic regime. The total energy remains bounded and exhibits small periodic variations, which confirms the long-time stability of the IRK discretization.

\begin{figure}[H]
\centering
\subfloat[$F_x$]{\includegraphics[width=0.33\textwidth]{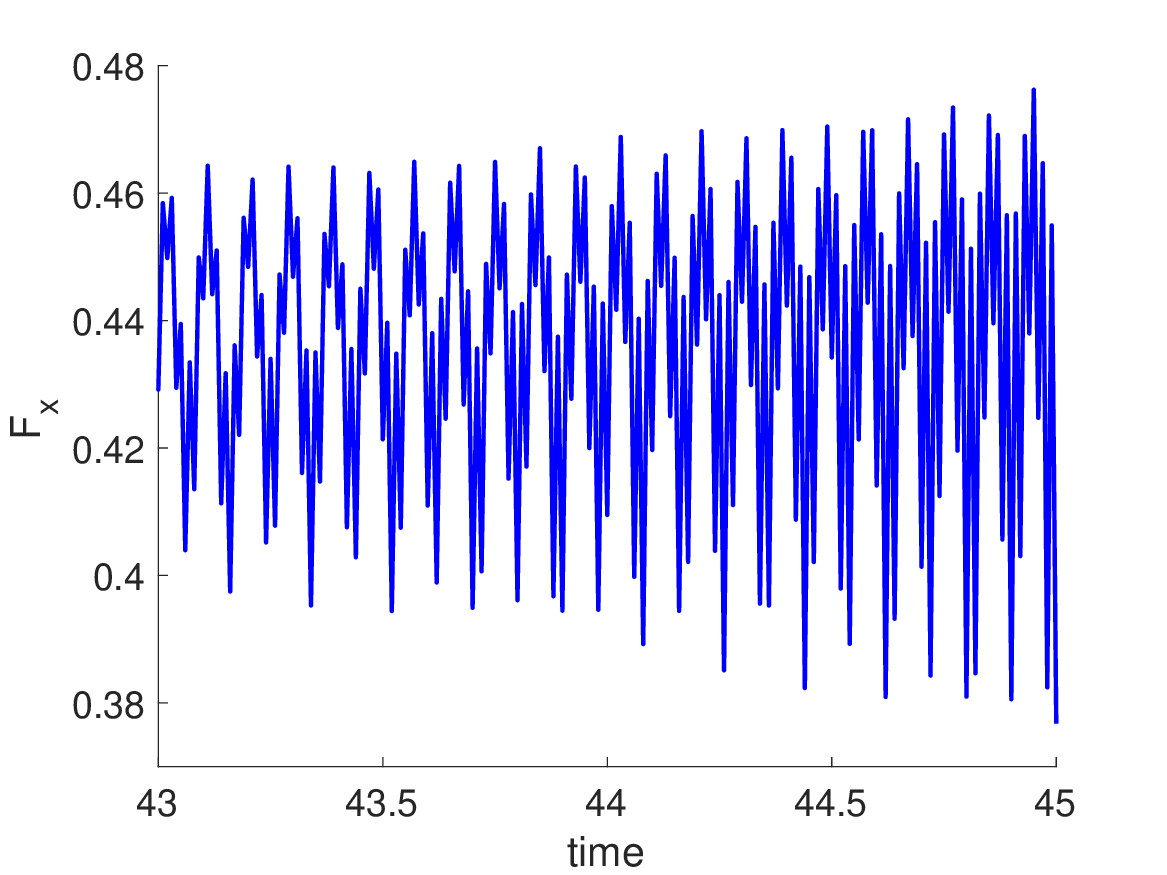}}
~~
\subfloat[$F_y$]{\includegraphics[width=0.33\textwidth]{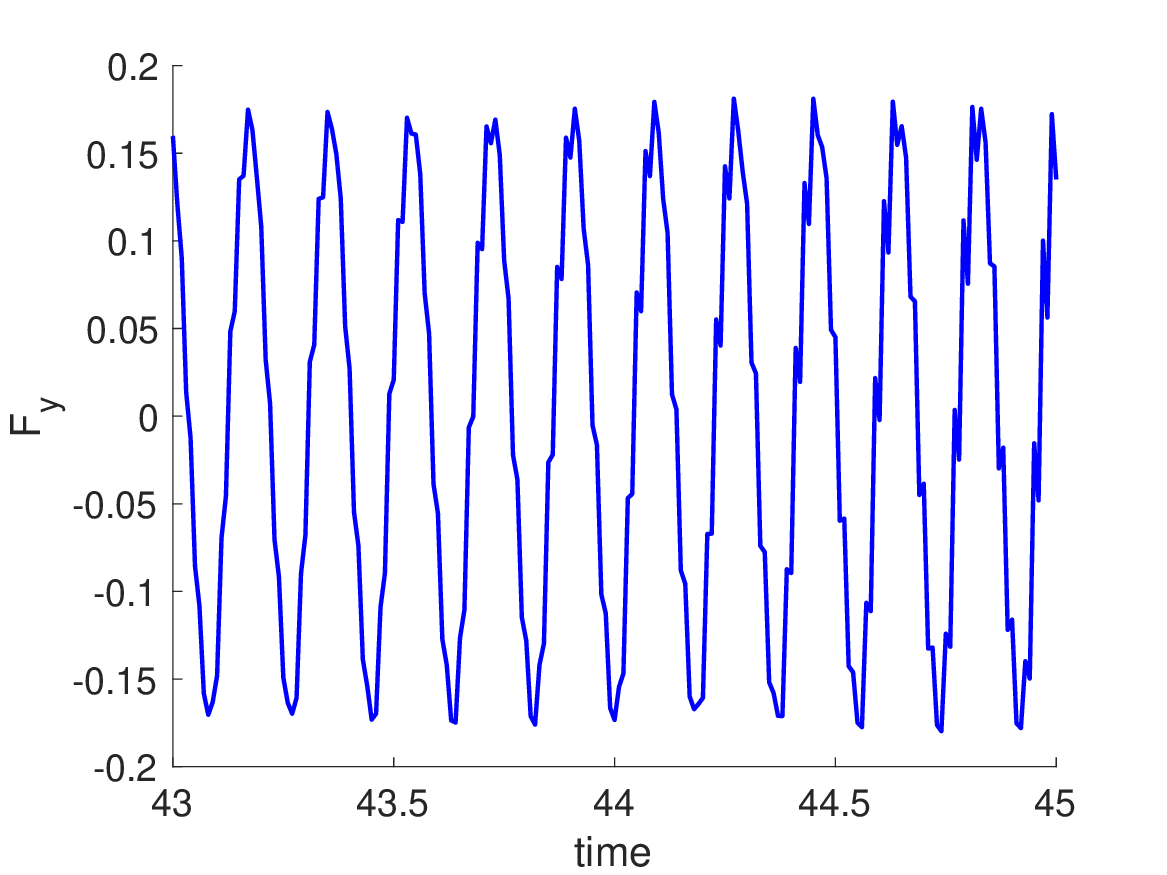}}
~~
\subfloat[Total energy]{\includegraphics[width=0.33\textwidth]{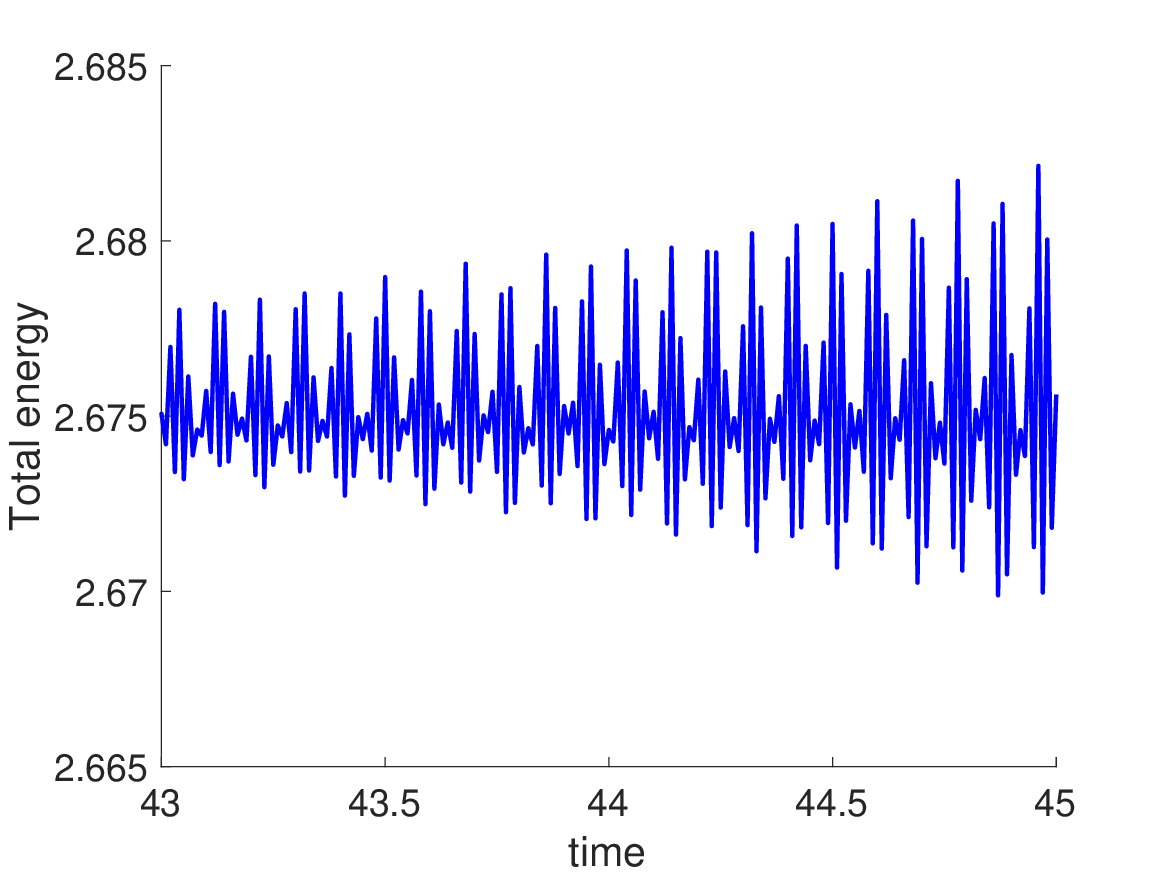}}
\caption{Turek; HF instability using BDF2 + Newmark.}
\label{fig:turek_HF_bdf2_instability}
\end{figure}

For comparison, we also consider the time integrator with BDF2 for the fluid and the Newmark scheme ($\beta=\frac14$, $\gamma=\frac12$)  for the solid. The corresponding long time histories of drag, lift, and total energy, shown in Figure \ref{fig:turek_HF_bdf2_instability}, reveal an instability: the force amplitudes and the total energy grow in time: this indicates that this time integration pair is not robust in this FSI setting.

\begin{figure}[H]
\centering
\subfloat[$F_x$]{\includegraphics[width=0.33\textwidth]{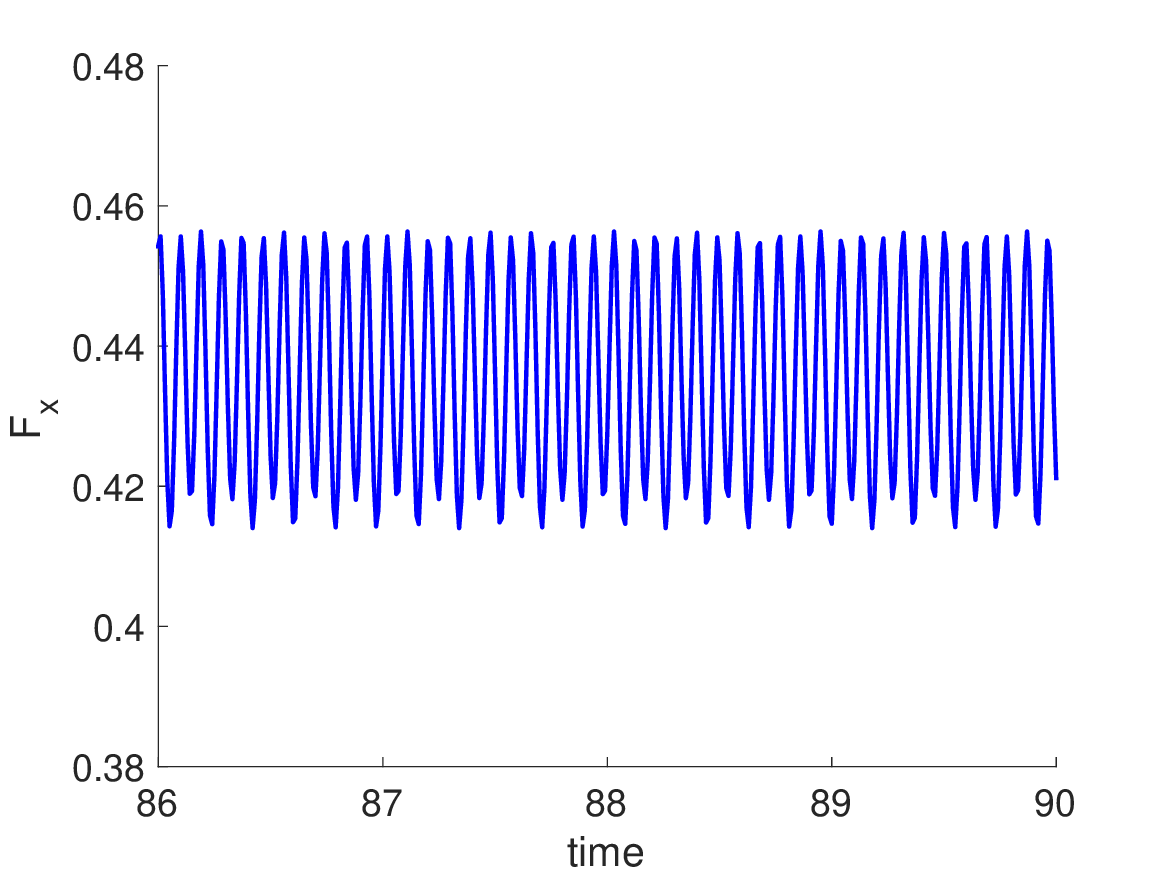}}
~~
\subfloat[$F_y$]{\includegraphics[width=0.33\textwidth]{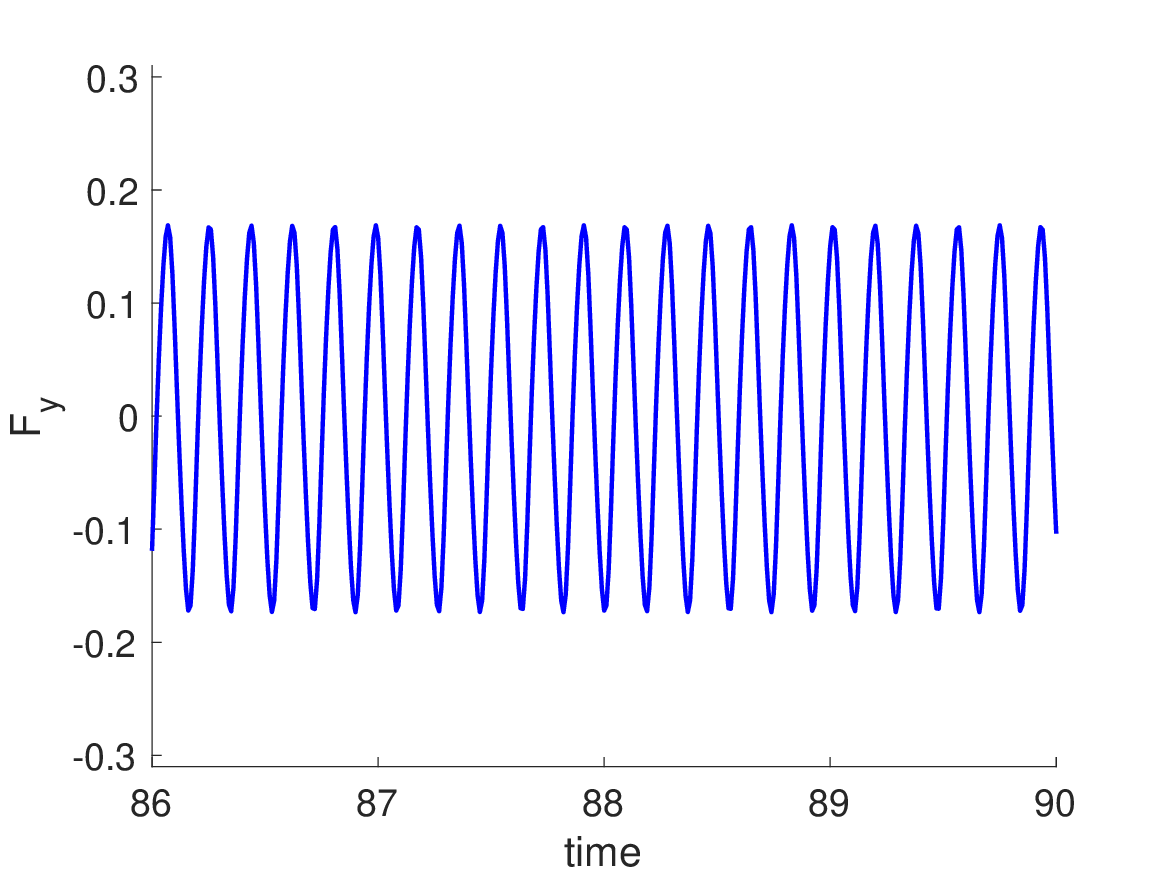}}
~~
\subfloat[Total energy]{\includegraphics[width=0.33\textwidth]{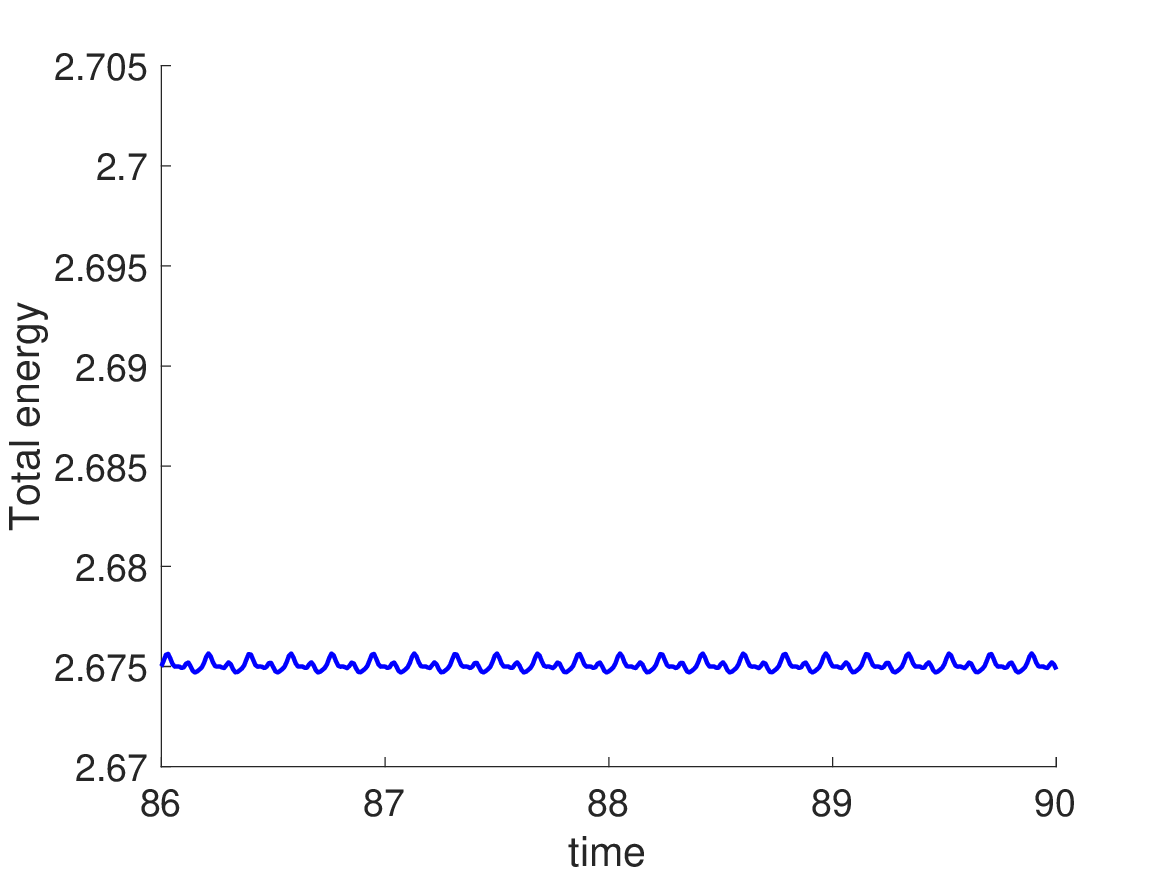}}
\caption{Turek; stable long-time HF results obtained using Newmark method with damping.}
\label{fig:turek_HF_bdf2_stab}
\end{figure}

To mitigate this problem, we introduce numerical damping in the Newmark method and choose $\beta=\frac13$ and $\gamma=\frac35$ ($2\beta \geq \gamma \geq \frac12$ to ensure unconditional stability). 
With these parameters, the instability observed with the undamped Newmark scheme disappears,
 as shown in Figure \ref{fig:turek_HF_bdf2_stab}. 
Note, however, that the predicted total energy is significantly lower than the one predicted by the IRK scheme: this  is consistent with the observation    that IRK schemes are less dissipative than damped BDF-Newmark methods.

\begin{figure}[H]
\centering
\subfloat[$F_x$]{\includegraphics[width=0.45\textwidth]{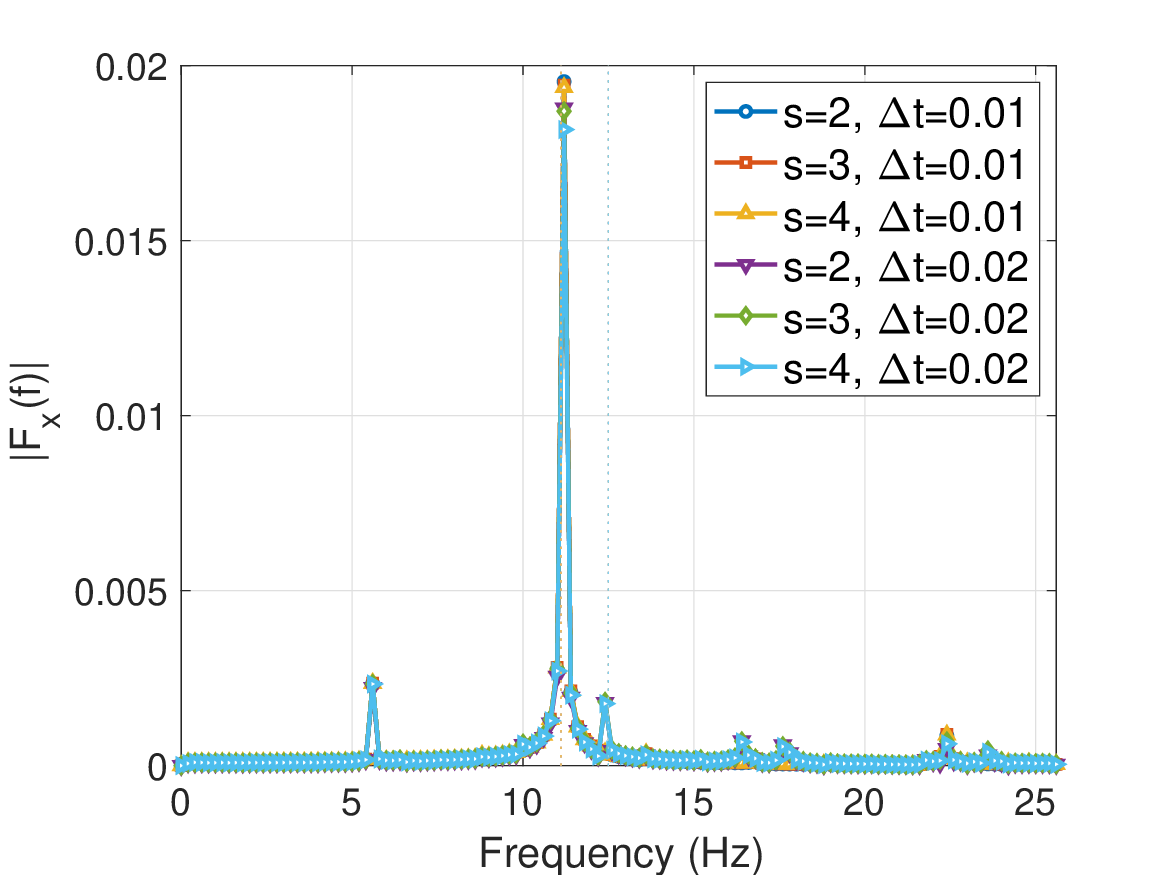}}
~~
\subfloat[$F_y$]{\includegraphics[width=0.45\textwidth]{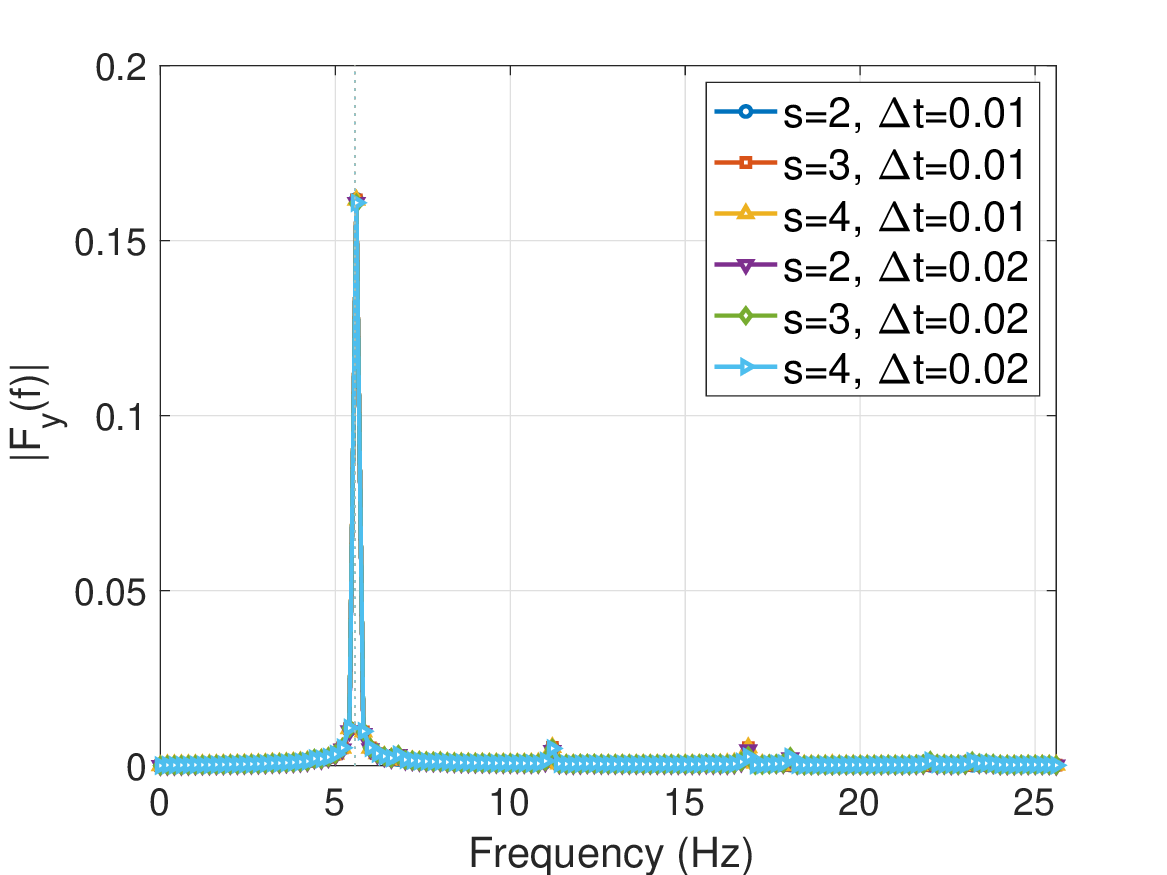}}
\caption{Turek; comparison of spectrum for 6 runs.}
\label{fig:hf_turek_irk357_spectrum_force}
\end{figure}

We further assess the influence of the IRK order and time step by performing simulations with $s=2,3,4$ and $\Delta t = 0.01, 0.02$\,s, and comparing the resulting drag and lift histories.
For each fixed time step, the force histories obtained with $s=2,s=3,s=4$ align well. 
All runs remain stable over the considered time interval, and the differences in amplitude and phase are very small. 
The corresponding time histories are reported in  \ref{app:turek_hf}. 
We further quantify these observations in the frequency domain. Figure \ref{fig:hf_turek_irk357_spectrum_force} shows the single-sided amplitude spectra of the drag ($F_x$) and lift ($F_y$) forces for six runs: \\
(1) $s=2,\,\Delta t=0.01$;
(2) $s=3,\,\Delta t=0.01$;
(3) $s=4,\,\Delta t=0.01$;
(4) $s=2,\,\Delta t=0.02$;
(5) $s=3,\,\Delta t=0.02$;
(6) $s=4,\,\Delta t=0.02$.
The dominant frequency and the overall spectral shape are essentially identical across all runs. 
Table \ref{tab:hf_turek_amp} summarizes the drag ($F_x$) and lift ($F_y$) amplitudes (mean $\pm$ half peak-to-peak) for all six runs and confirms the robustness of the IRK schemes with respect to $s$ and $\Delta t$.

\begin{table}[H]
\centering
\caption{Drag ($F_x$) and lift ($F_y$) amplitudes for the six runs (reported as mean $\pm$ half peak-to-peak).}
\label{tab:hf_turek_amp}
\begin{tabular}{cccccc}
\toprule
Run & $s$ & $\Delta t$ & $F_x$ amplitude & $F_y$ amplitude \\
\midrule
1 & 2 & 0.01 & $0.436 \pm 0.0198$ & $0.0029 \pm 0.1668$ \\
2 & 3 & 0.01 & $0.436 \pm 0.0198$ & $0.0029 \pm 0.1676$ \\
3 & 4 & 0.01 & $0.436 \pm 0.0197$ & $0.0029 \pm 0.1676$ \\
4 & 2 & 0.02 & $0.436 \pm 0.0190$ & $0.0026 \pm 0.1619$ \\
5 & 3 & 0.02 & $0.436 \pm 0.0190$ & $0.0030 \pm 0.1647$ \\
6 & 4 & 0.02 & $0.436 \pm 0.0184$ & $0.0030 \pm 0.1646$ \\
\bottomrule
\end{tabular}
\end{table}

\subsubsection{ROM results}
\begin{figure}[H]
  \centering
\subfloat[$t=32$]{\includegraphics[width=0.33\textwidth]{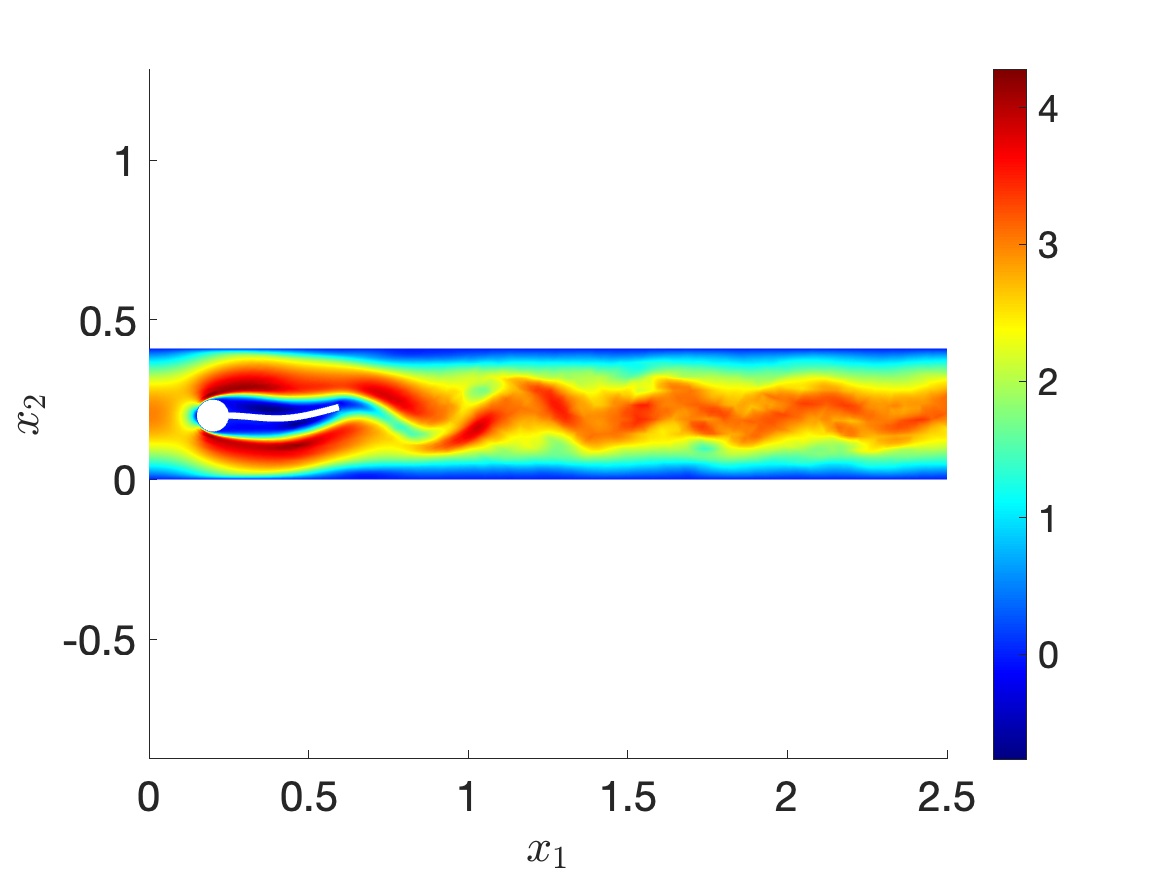}}
~~
\subfloat[$t=35$]{\includegraphics[width=0.33\textwidth]{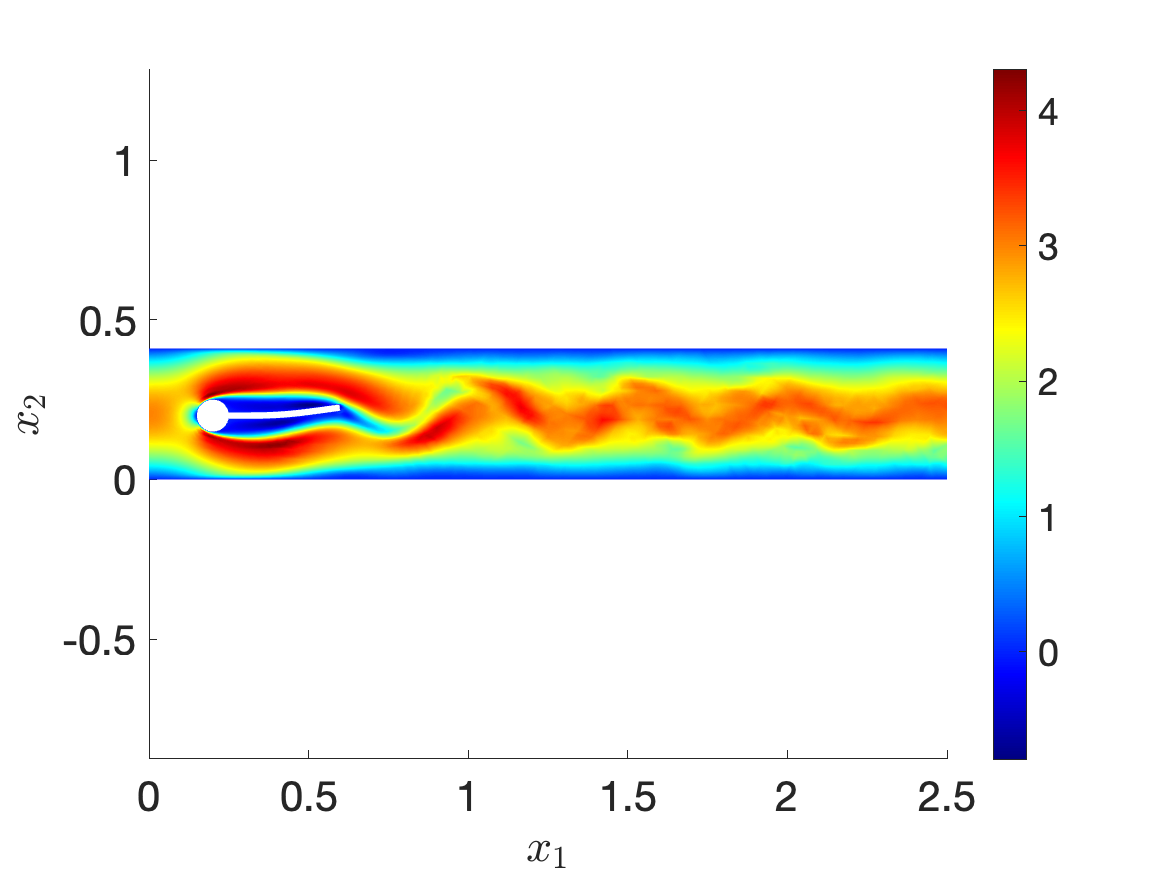}}
~~
\subfloat[$t=38$]{\includegraphics[width=0.33\textwidth]{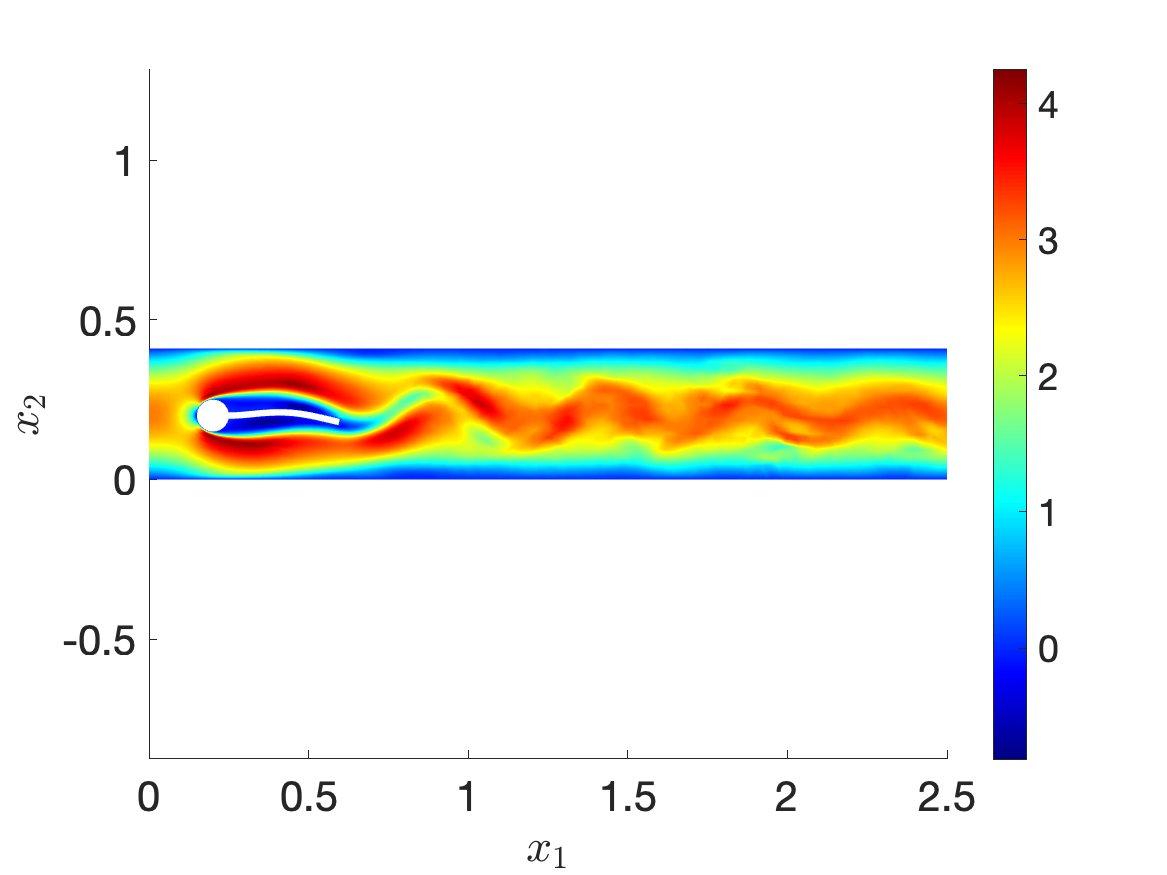}}
\\
\vspace*{-0.4cm}
\subfloat[$t=32$ (error)]{\includegraphics[width=0.33\textwidth]{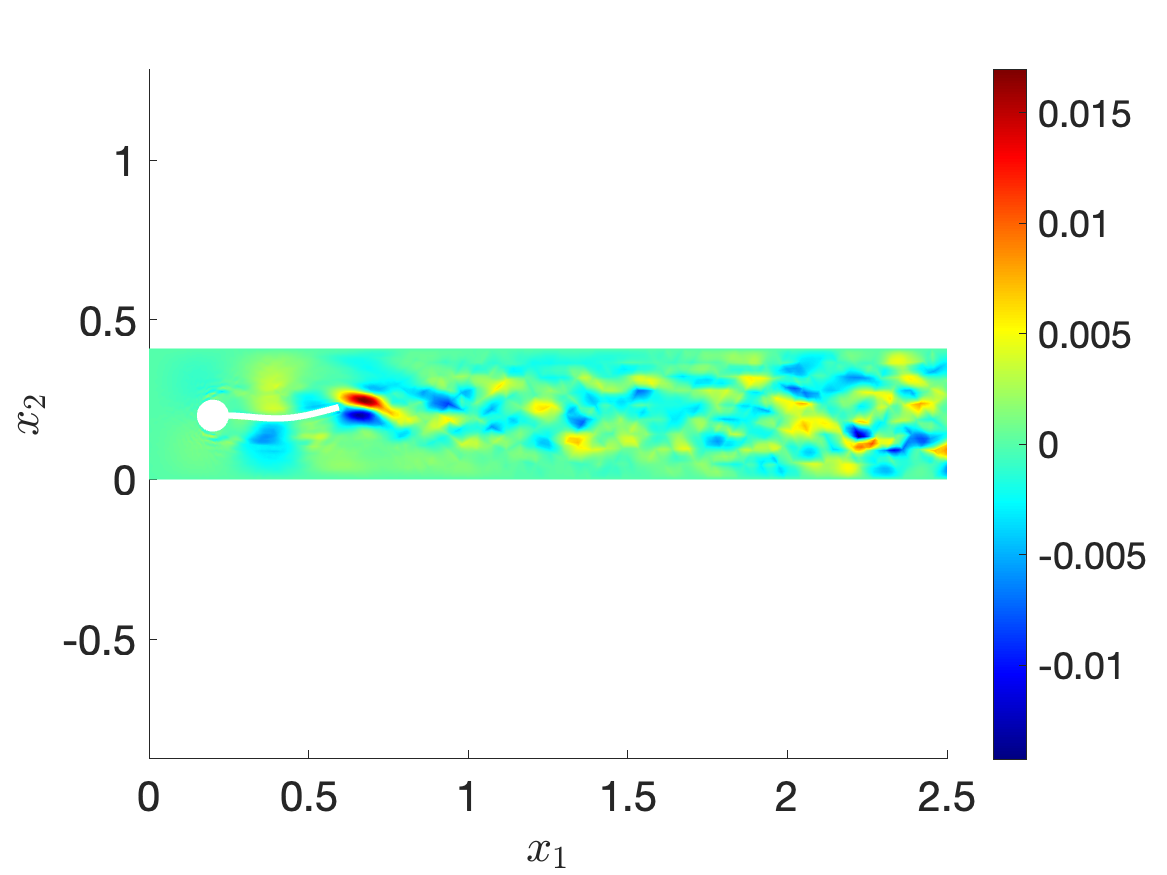}}
~~
\subfloat[$t=35$ (error)]{\includegraphics[width=0.33\textwidth]{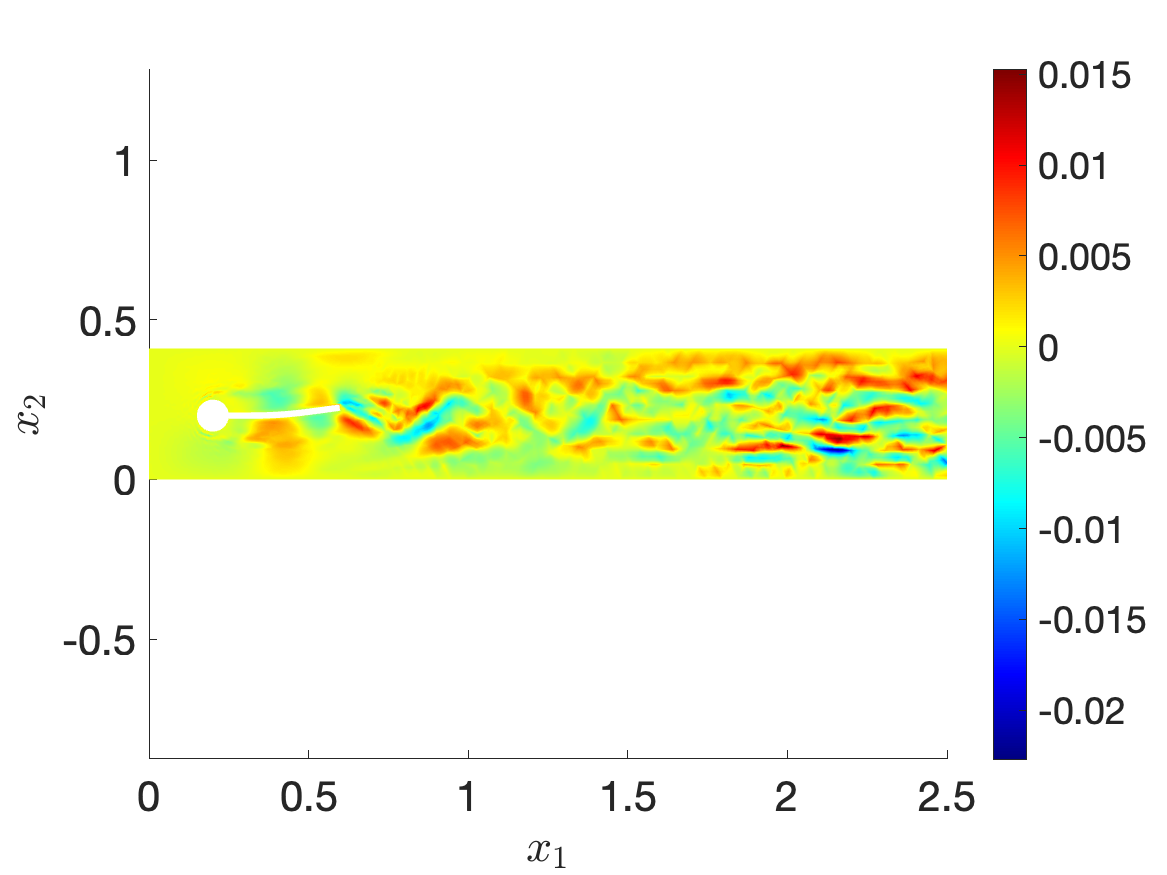}}
~~
\subfloat[$t=38$ (error)]{\includegraphics[width=0.33\textwidth]{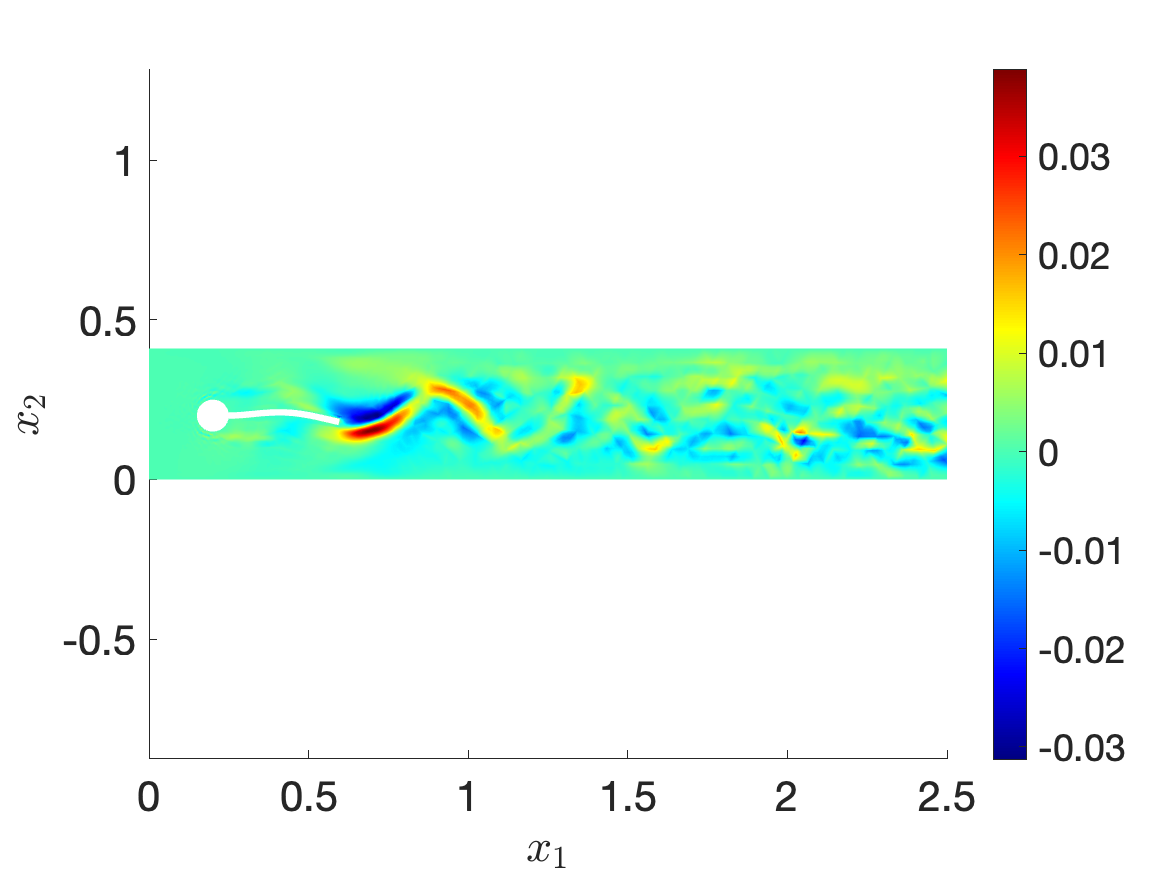}}
\caption{Turek; stable ROM results obtained using IRK ($s=2,\Delta t=0.01\,$s).}
\label{fig:turek_ROM_irk3_stab}
\end{figure}

We next build ROMs from IRK snapshots and assess their performance. 
We train with snapshots from time instants $t=30-33\,$s, then we perform ROM calculation for $t=30-39\,$s. 
We first consider $s=2,\Delta t=0.01\,$s for the time integrator. 
We show the horizontal velocity field computed by the ROM at $t=32,35,38\,$s for a POD tolerance  tol$_{\rm POD}=10^{-6}$ in Figure \ref{fig:turek_ROM_irk3_stab}. The errors remain  small compared to the magnitude of the flow, which indicates that the ROM correctly reproduces the main flow features and the beam motion. 

\begin{figure}[H]
\centering
\subfloat[ROM errors]{\includegraphics[width=0.45\textwidth]{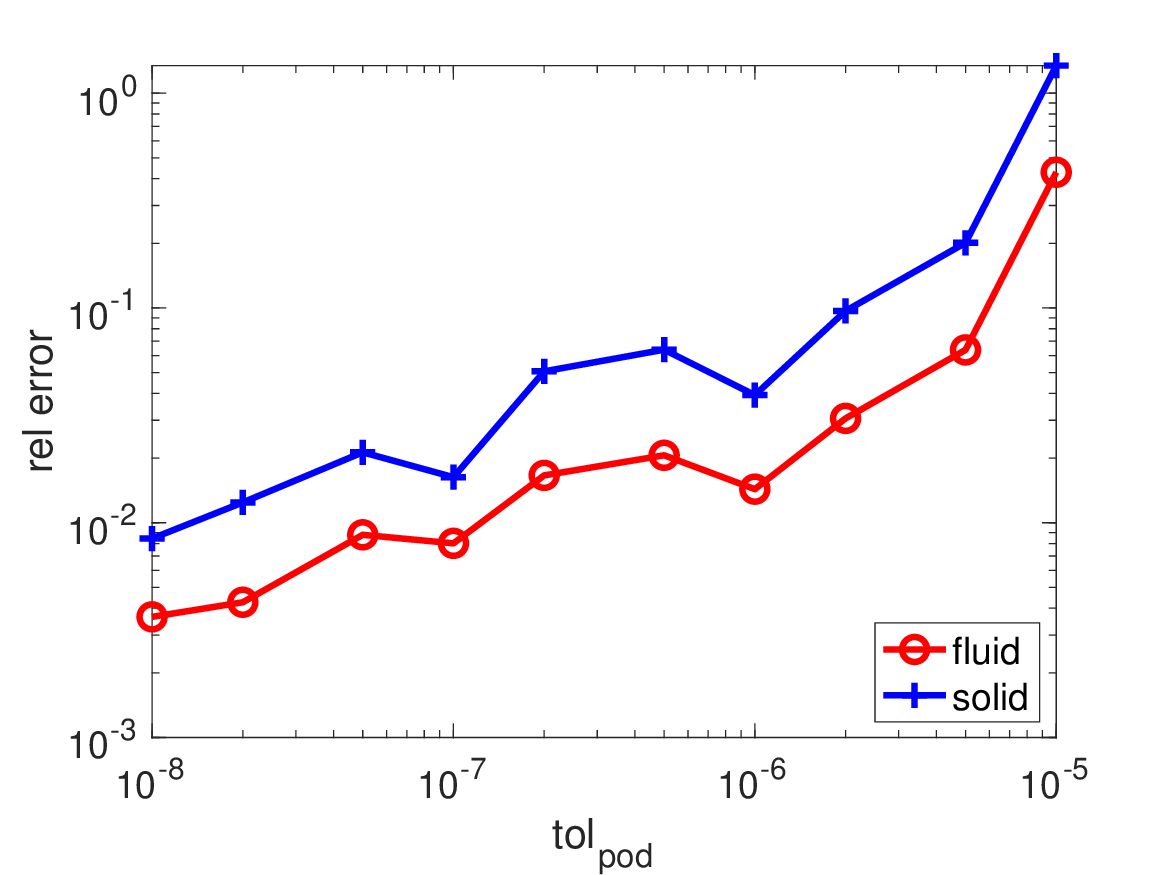}}
~~
\subfloat[number of modes]{\includegraphics[width=0.45\textwidth]{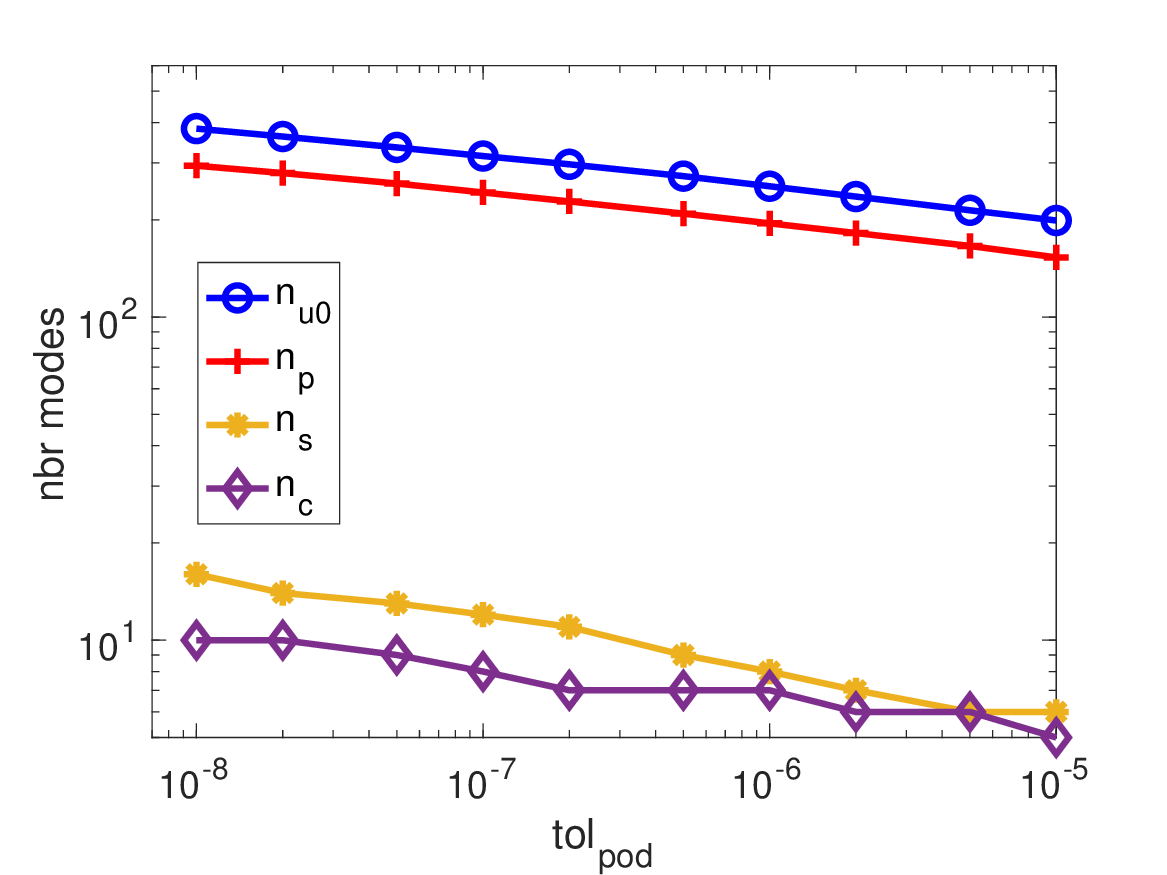}}
\caption{Turek; ROM errors and number of modes in terms of POD tolerance  ($s=2,\Delta t=0.01\,$s).}
\label{fig:turek_ROM_err_s2_0d01}
\end{figure}

We show $H^1\times L^2$ relative error of the velocity-pressure pair as function of the POD tolerance, together with the corresponding number of reduced modes in Figure \ref{fig:turek_ROM_err_s2_0d01}. 
As expected, the use of tighter tolerances  tol$_{\rm POD}$ leads to smaller errors but larger reduced spaces. The fluid velocity basis is enriched with supremizer modes  to ensure the inf-sup stability, and thus the total number of velocity modes for the fluid is $n_{\rm u}=n_{\rm{u}0}+n_{\rm p}$. \\

\begin{figure}[H]
\centering
\subfloat[$F_x$]{\includegraphics[width=0.33\textwidth]{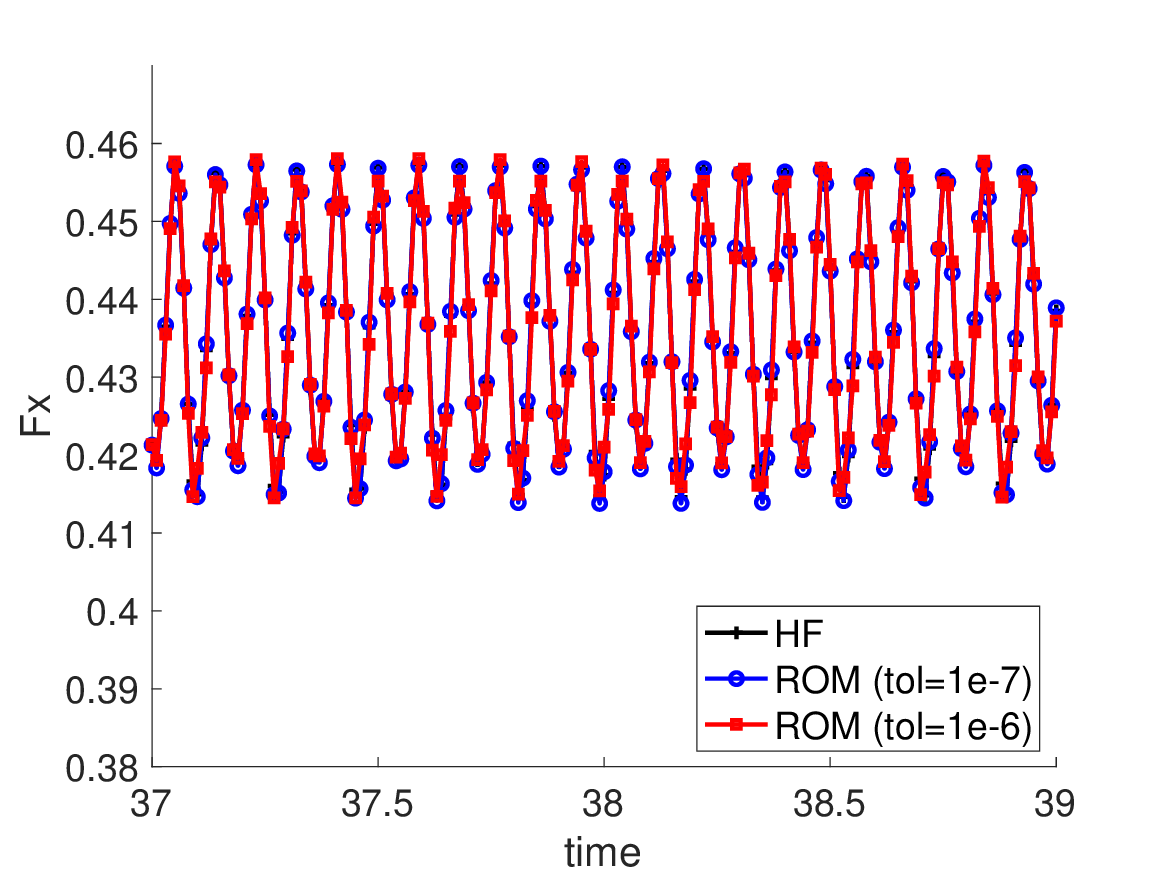}}
~~
\subfloat[$F_y$]{\includegraphics[width=0.33\textwidth]{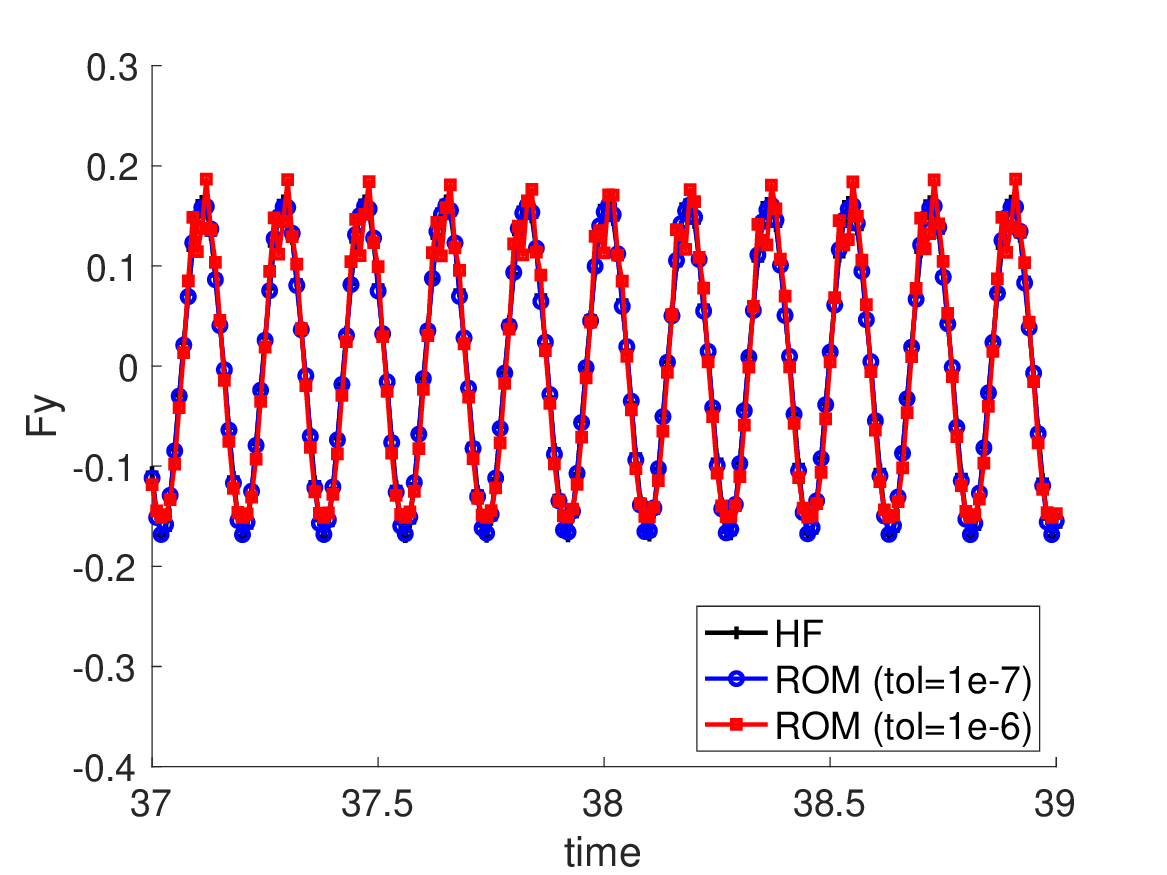}}
~~
\subfloat[Total energy]{\includegraphics[width=0.33\textwidth]{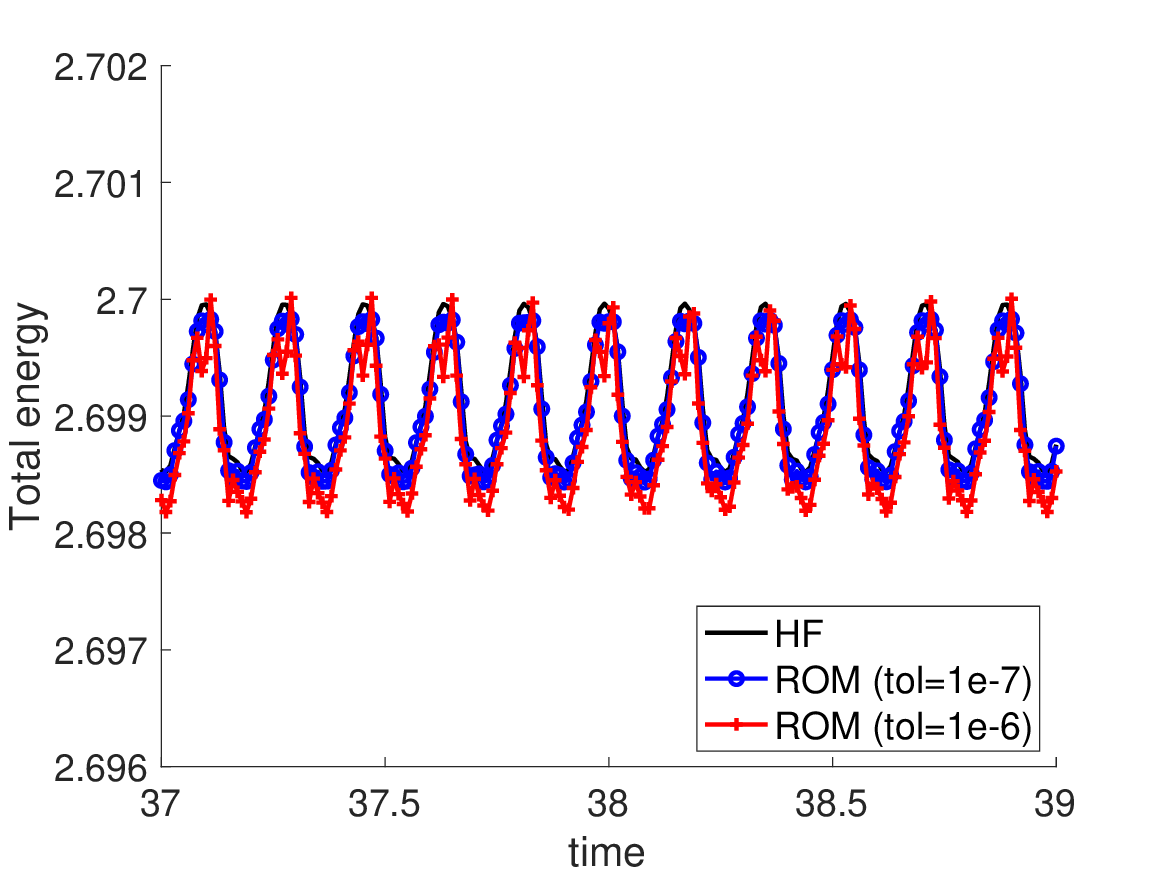}}
\caption{Turek; comparison between HF and two ROM results ($s=2,\Delta t=0.01\,$s).}
\label{fig:turek_compr_s2_0d01}
\end{figure}

In Figure \ref{fig:turek_compr_s2_0d01}, we compare high-fidelity and ROM predictions for drag, lift, and total energy for two representative ROMs built with $\mathrm{tol}_{\rm POD}=10^{-7}$ and $10^{-6}$. Both ROMs can  reproduce the force signals and total energy; the ROM with the smaller tolerance provides a slightly better phase match.

\begin{figure}[H]
\centering
\subfloat[$F_x$]{\includegraphics[width=0.45\textwidth]{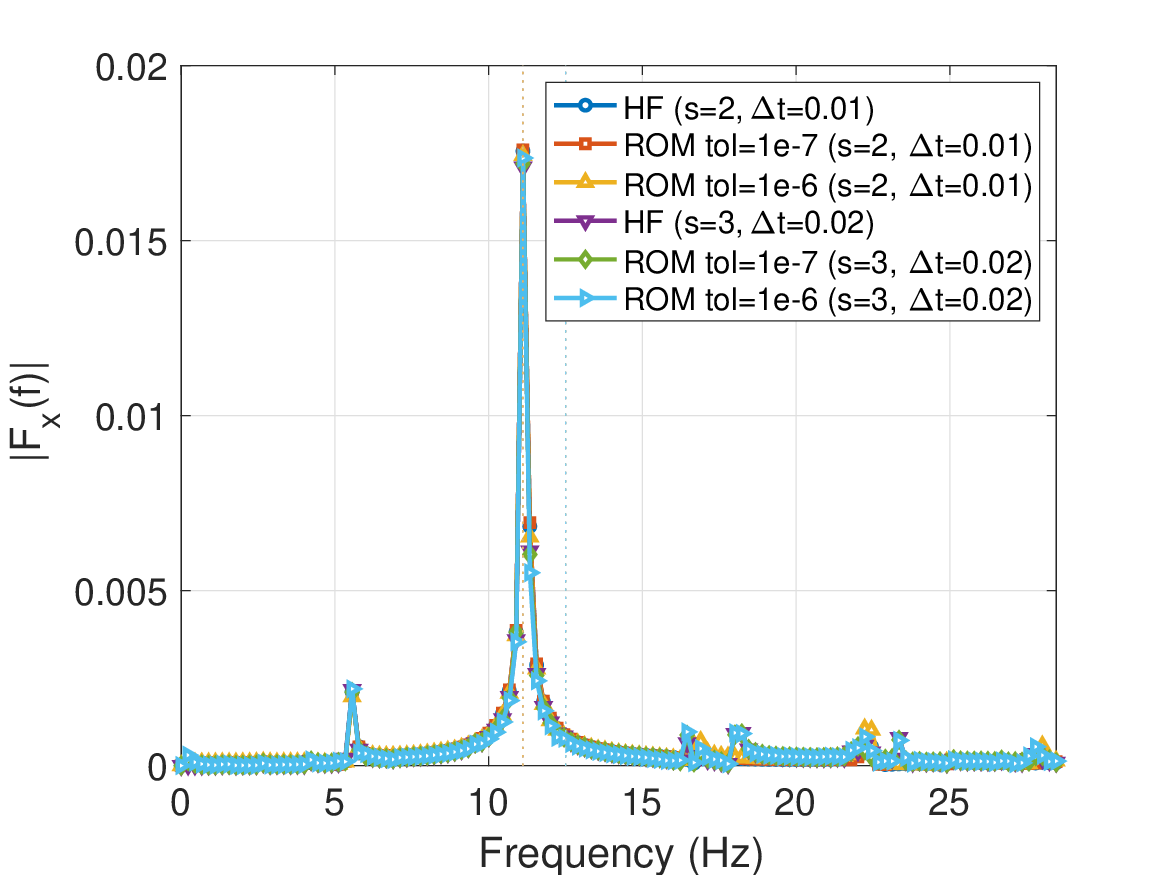}}
~~
\subfloat[$F_y$]{\includegraphics[width=0.45\textwidth]{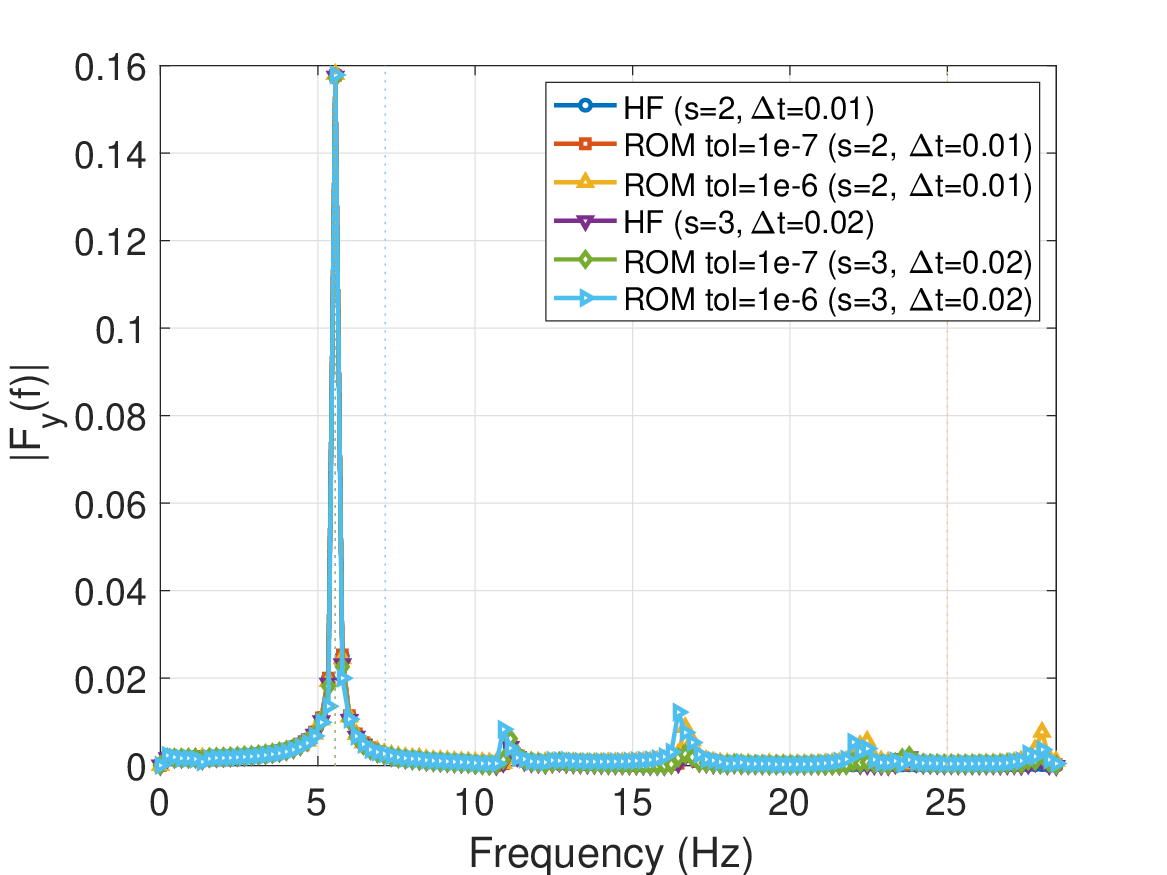}}
\caption{Turek; spectra of drag and lift forces obtained by ROMs for different $s$, $\Delta t$, and $\text{tol}_{\text{POD}}$.}
\label{fig:ROM_turek_irk35_spectrum_force}
\end{figure}

We also repeat the analysis for IRK with $s=3$ and $\Delta t = 0.02$\,s. The corresponding velocity fields, error plots, and error versus tolerance curves confirm the same qualitative behavior: ROMs remain stable and accurate, and the dimension accuracy trade-off with respect to $\mathrm{tol}_{\rm POD}$ is similar to the $s=2$ case. These additional results are reported in \ref{app:turek_rom_s3}. 
Figure \ref{fig:ROM_turek_irk35_spectrum_force} shows the single-sided amplitude spectra of the drag ($F_x$) and lift ($F_y$) forces computed by the ROMs for different numbers of IRK stages $s$, time steps $\Delta t$, and POD tolerances $\mathrm{tol}_{\mathrm{POD}}$. All ROM spectra closely match the dominant frequency and overall spectral distribution of the high-fidelity reference, indicating that the proposed ROM framework accurately captures the frequency content of fluid forces over a range of numerical settings.

\begin{figure}[H]
\centering
\subfloat[$F_x$]{\includegraphics[width=0.44\textwidth]{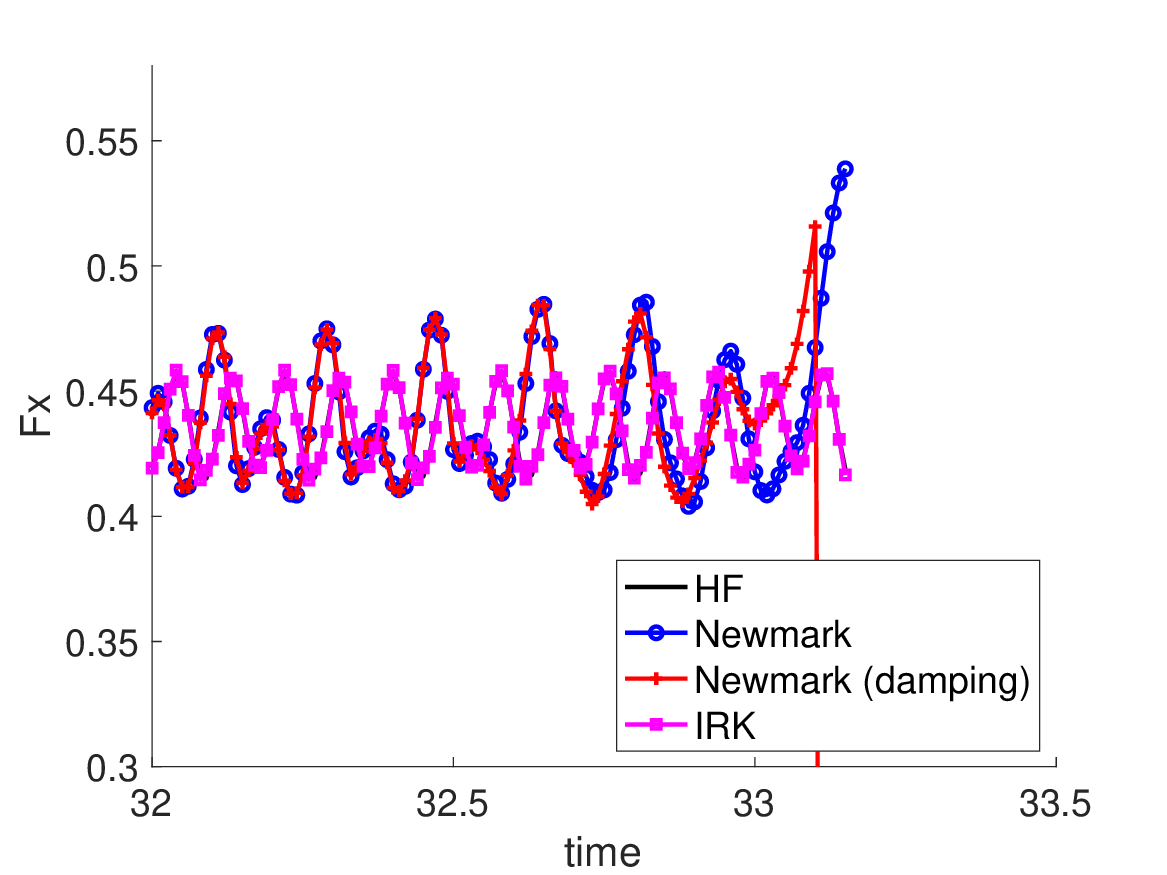}}
~~
\subfloat[$F_y$]{\includegraphics[width=0.44\textwidth]{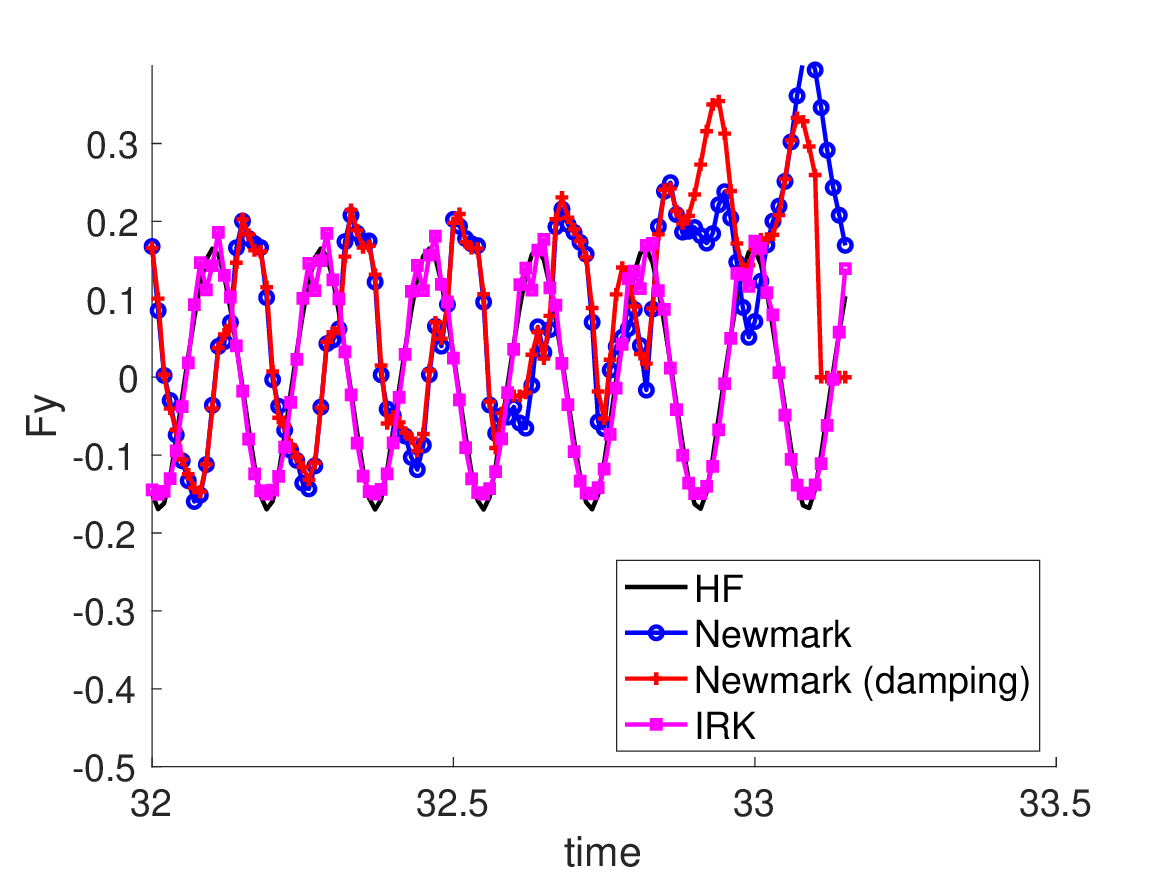}}
\caption{Turek; comparison between BDF2 + Newmark and IRK ROM results ($s=2,\Delta t=0.01\,$s, tol$_{\rm POD} = 2 \cdot 10^{-7}$).}
\label{fig:turek_compr_s2_0d01_bdf2_irk3}
\end{figure}

A key advantage of the present method with respect to the flux-based strategy of \cite{taddei2025optimization}
is its long time robustness at the reduced level. To illustrate this, we construct a POD basis from IRK snapshots ($s=2$) with $\mathrm{tol}_{\rm POD}=2\times 10^{-7}$ and compare the effect of advancing this basis in time with different time integrators. When the ROM is integrated with IRK ($s=2$), the solution remains stable and accurately follows the high-fidelity reference. In contrast,
the use of the BDF2-Newmark scheme for time integration  renders the ROM unstable, irrespective of numerical (Newmark) damping. This behaviour is illustrated in Figure \ref{fig:turek_compr_s2_0d01_bdf2_irk3} for the drag and lift forces.

\begin{figure}[H]
\centering
\subfloat[ROM errors]{\includegraphics[width=0.45\textwidth]{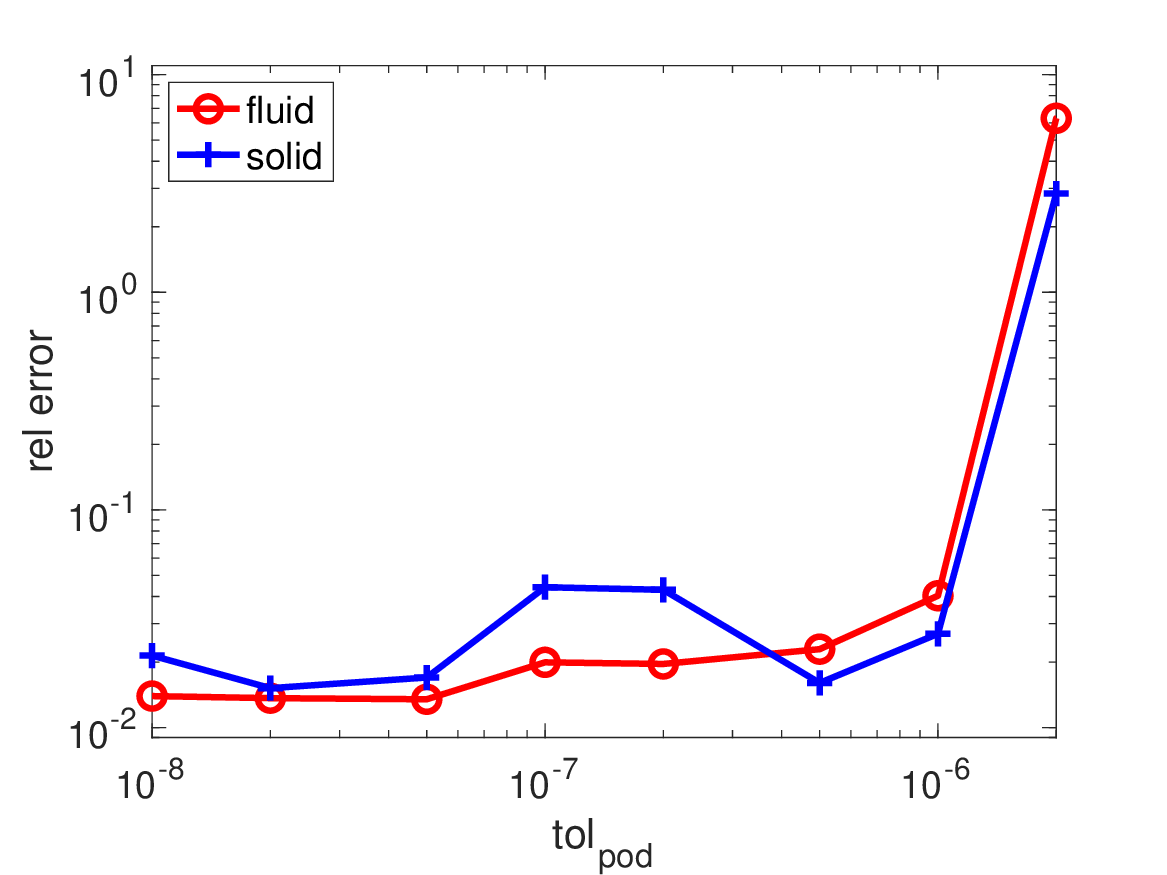}}
~~
\subfloat[number of modes]{\includegraphics[width=0.45\textwidth]{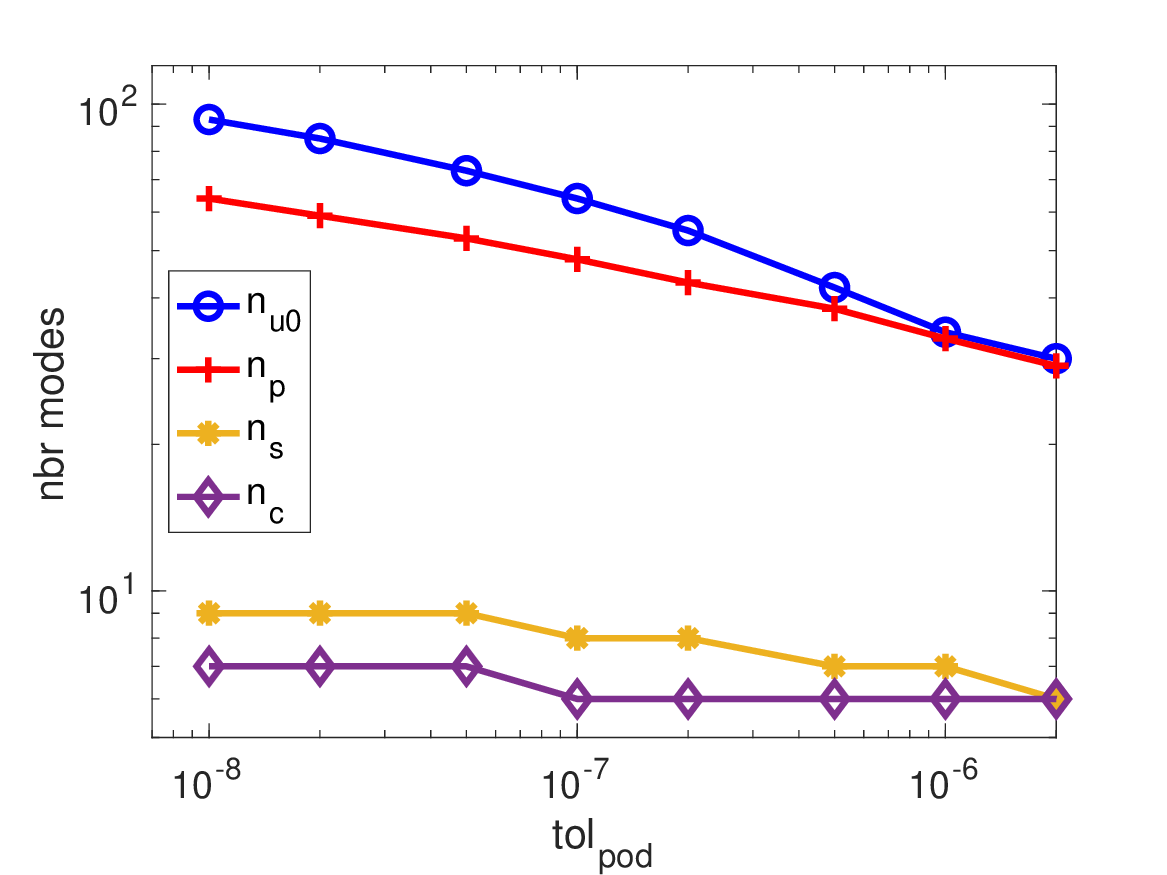}}
\caption{Turek; ROM errors and number of modes in terms of POD tolerance using periodic snapshots  ($s=2,\Delta t=0.01\,$s).}
\label{fig:turek_ROM_err_periodic}
\end{figure}

We also investigate ROMs trained on periodic snapshots only. In this case, we use snapshots from time instants $t=33-35\,$s, as training set. The snapshots are in a fully developed periodic regime. We then perform the ROM calculation for $T=33-39\,$s. Figure \ref{fig:turek_ROM_err_periodic} shows the $H^1 \times L^2$ relative error of the velocity–pressure pair and the number of reduced modes as functions of $\mathrm{tol}_{\rm POD}$.

\begin{figure}[H]
\centering
\subfloat[$F_x$]{\includegraphics[width=0.45\textwidth]{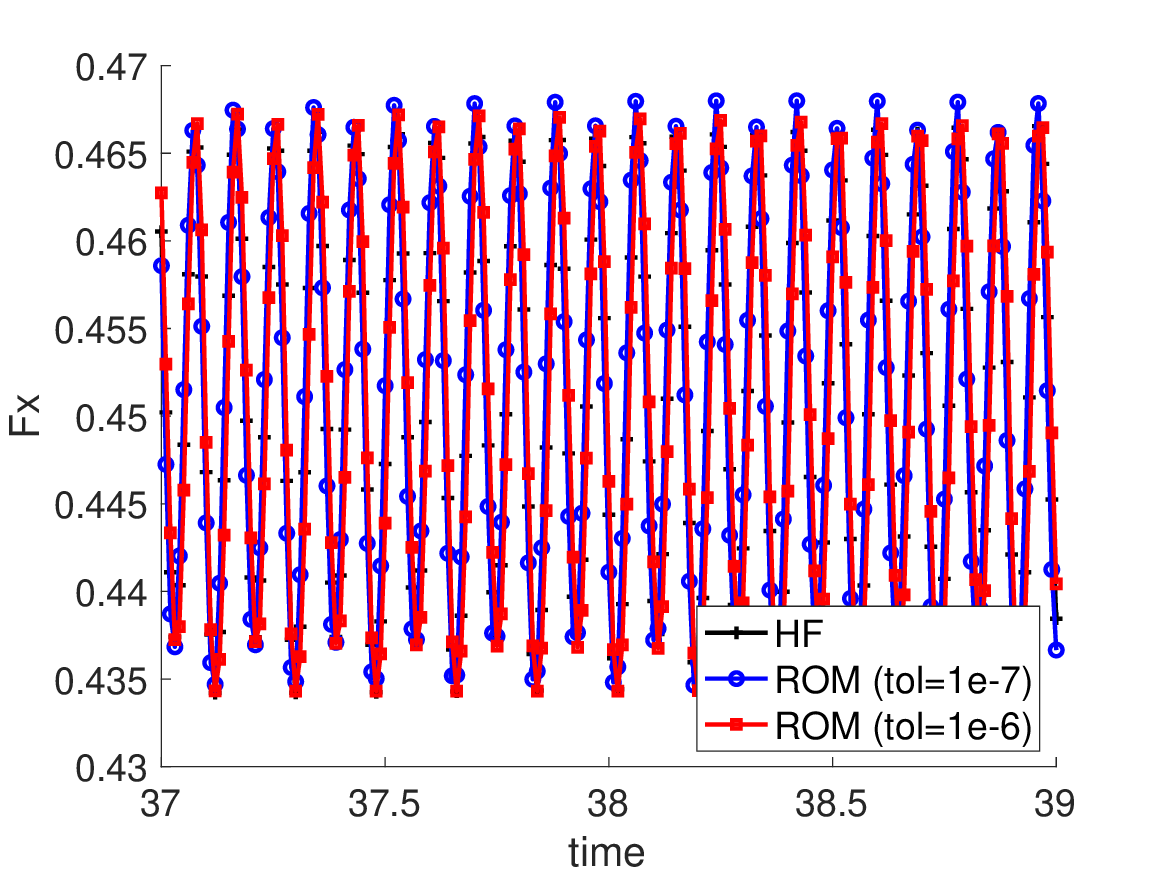}}
~~
\subfloat[$F_y$]{\includegraphics[width=0.45\textwidth]{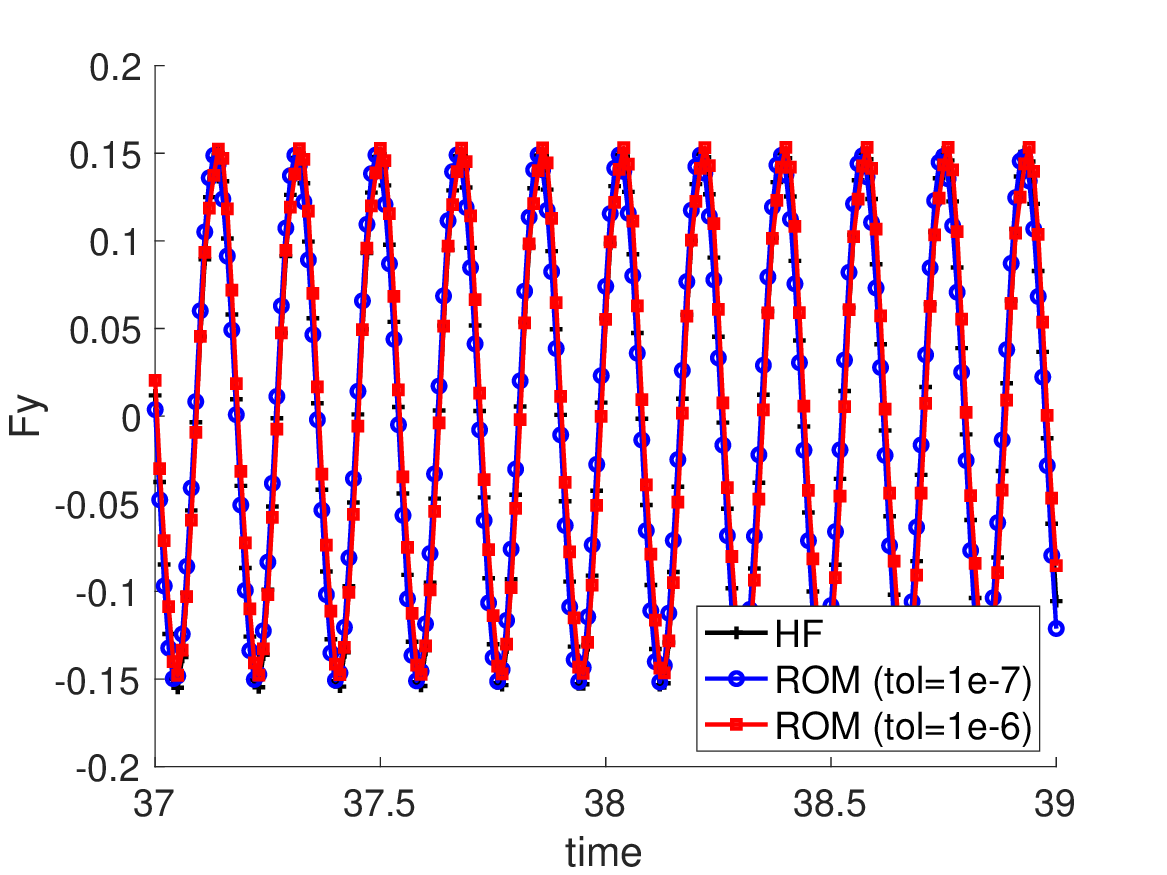}}
\\
\subfloat[$d_x$]{\includegraphics[width=0.45\textwidth]{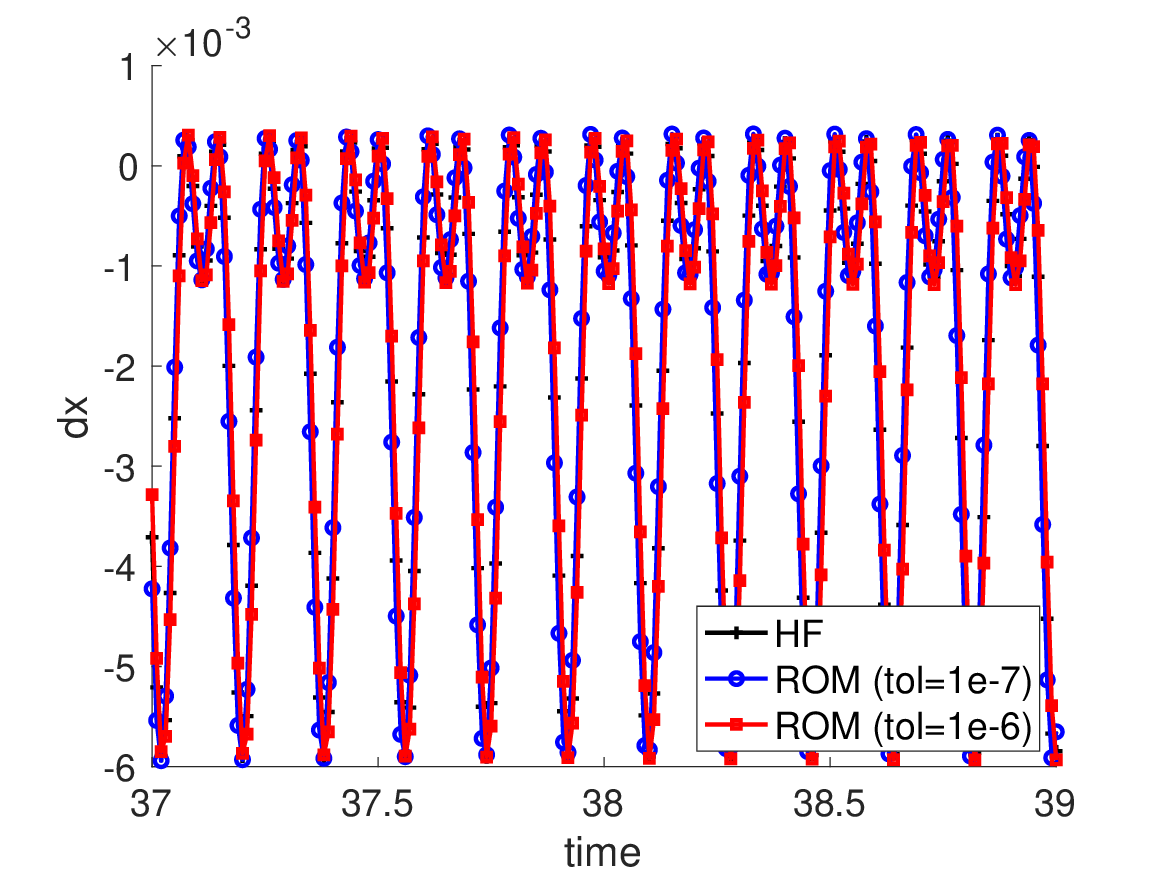}}
~~
\subfloat[$d_y$]{\includegraphics[width=0.45\textwidth]{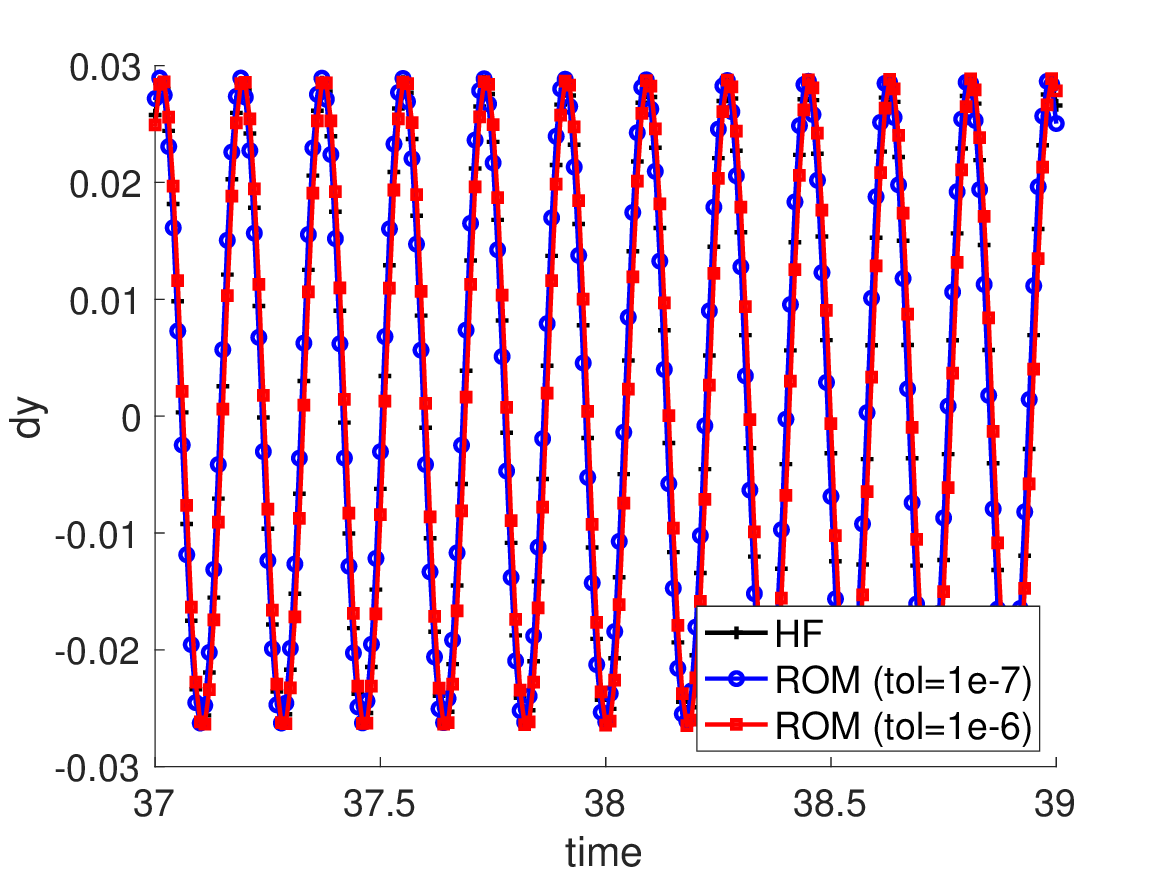}}
\caption{Turek; comparison between HF and two ROM results (test point $1$). }
\label{fig:ROM_turek_irk3_compr_param1}
\end{figure}

\begin{figure}[H]
\centering
\subfloat[$F_x$]{\includegraphics[width=0.45\textwidth]{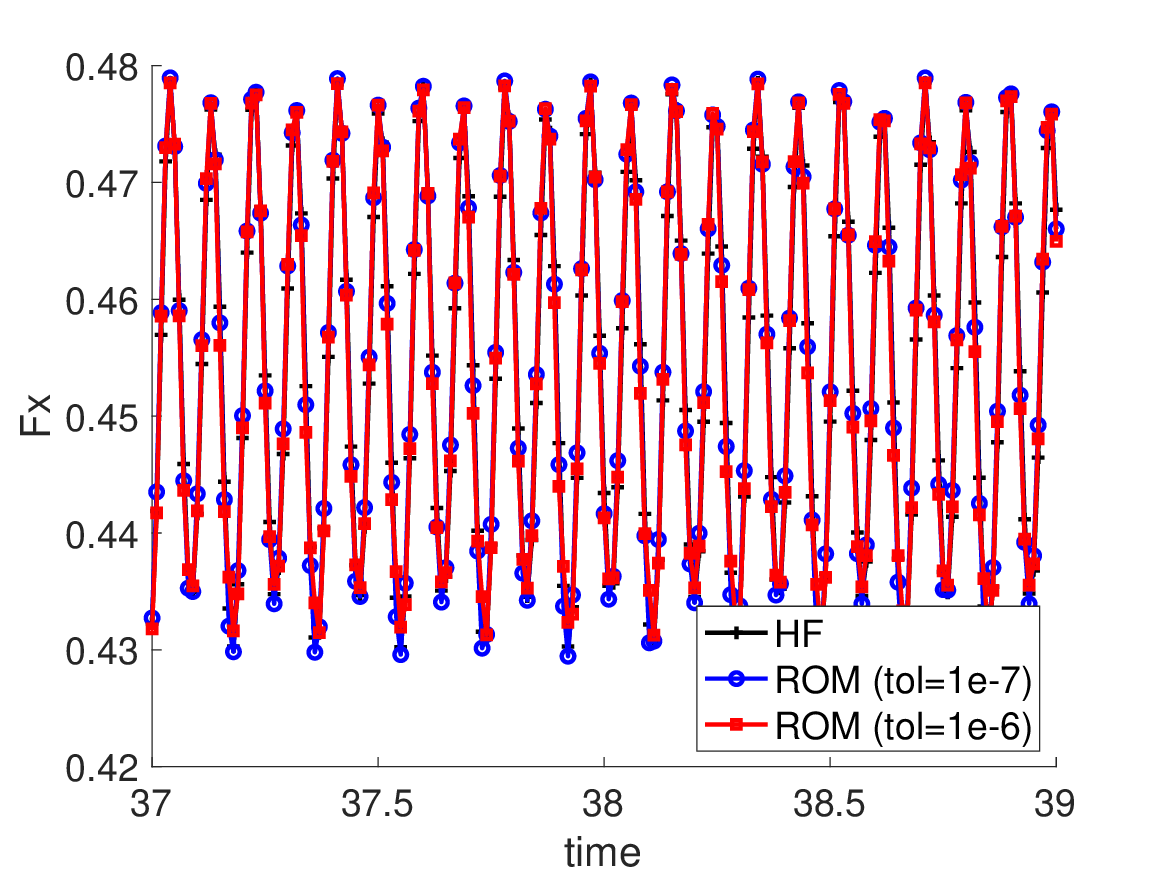}}
~~
\subfloat[$F_y$]{\includegraphics[width=0.45\textwidth]{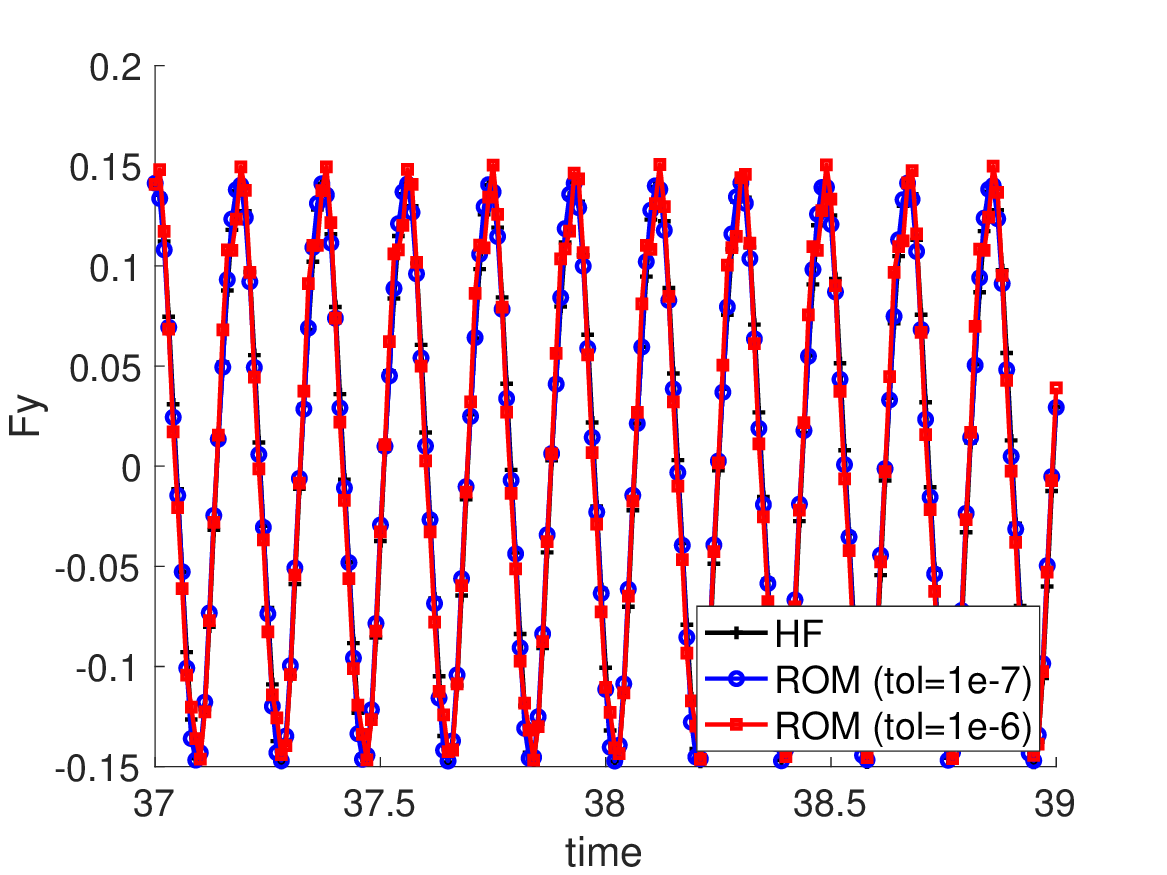}}
\\
\subfloat[$d_x$]{\includegraphics[width=0.45\textwidth]{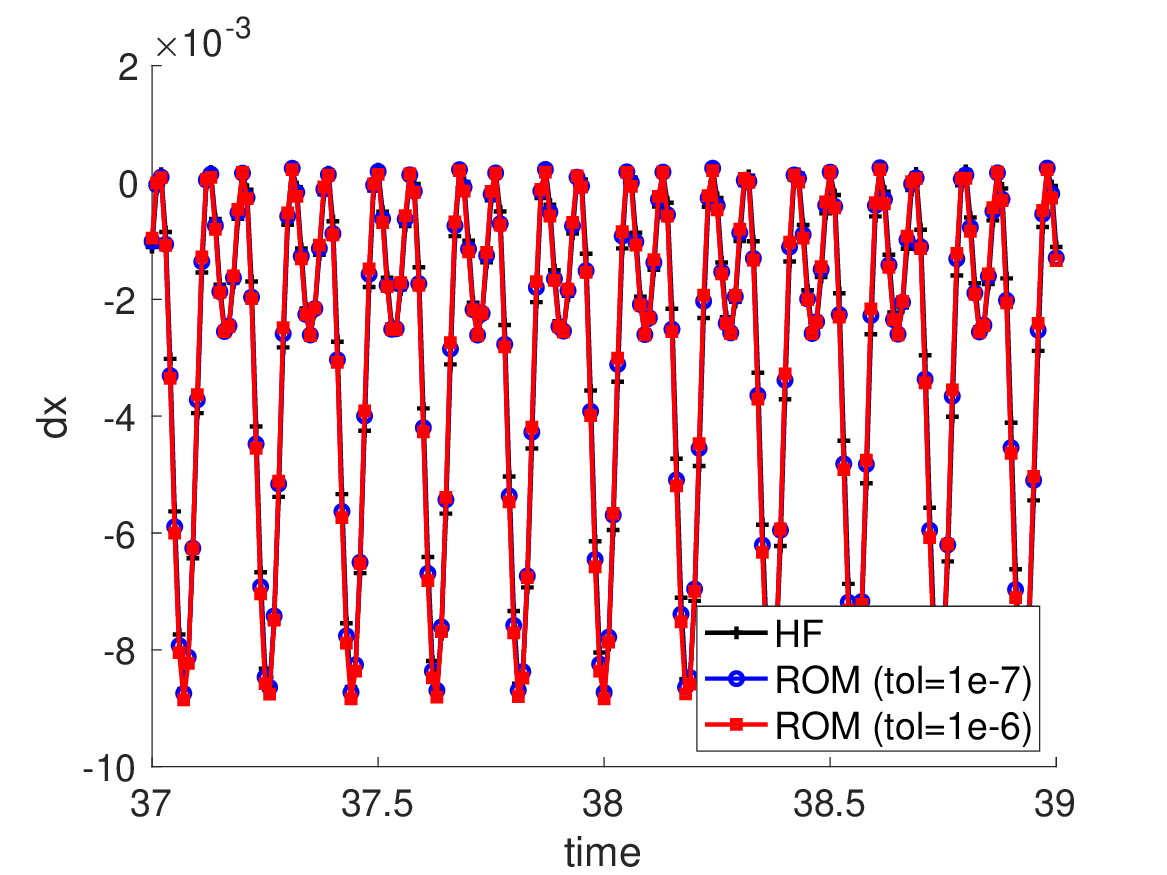}}
~~
\subfloat[$d_y$]{\includegraphics[width=0.45\textwidth]{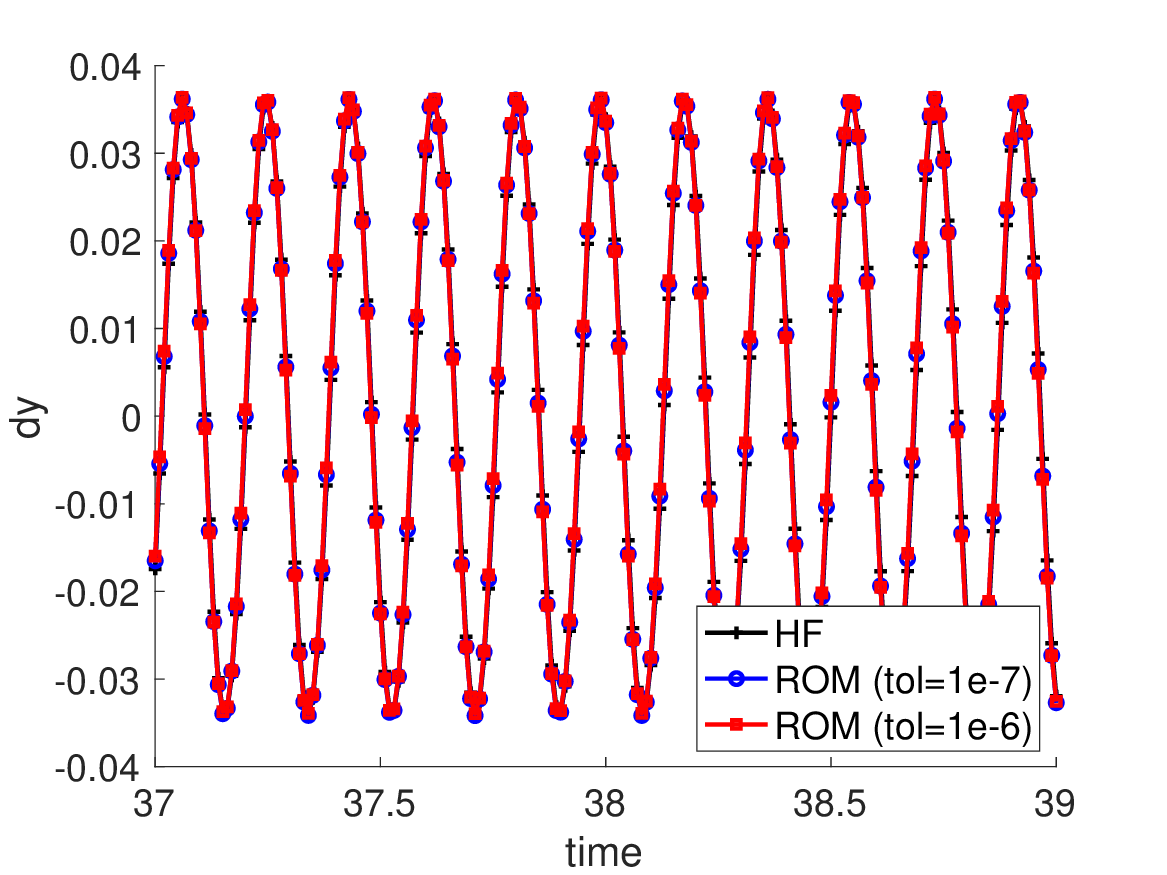}}
\caption{Turek; comparison between HF and two ROM results (test point $2$). }
\label{fig:ROM_turek_irk3_compr_param2}
\end{figure}

\begin{figure}[H]
\centering
\subfloat[$F_x$]{\includegraphics[width=0.45\textwidth]{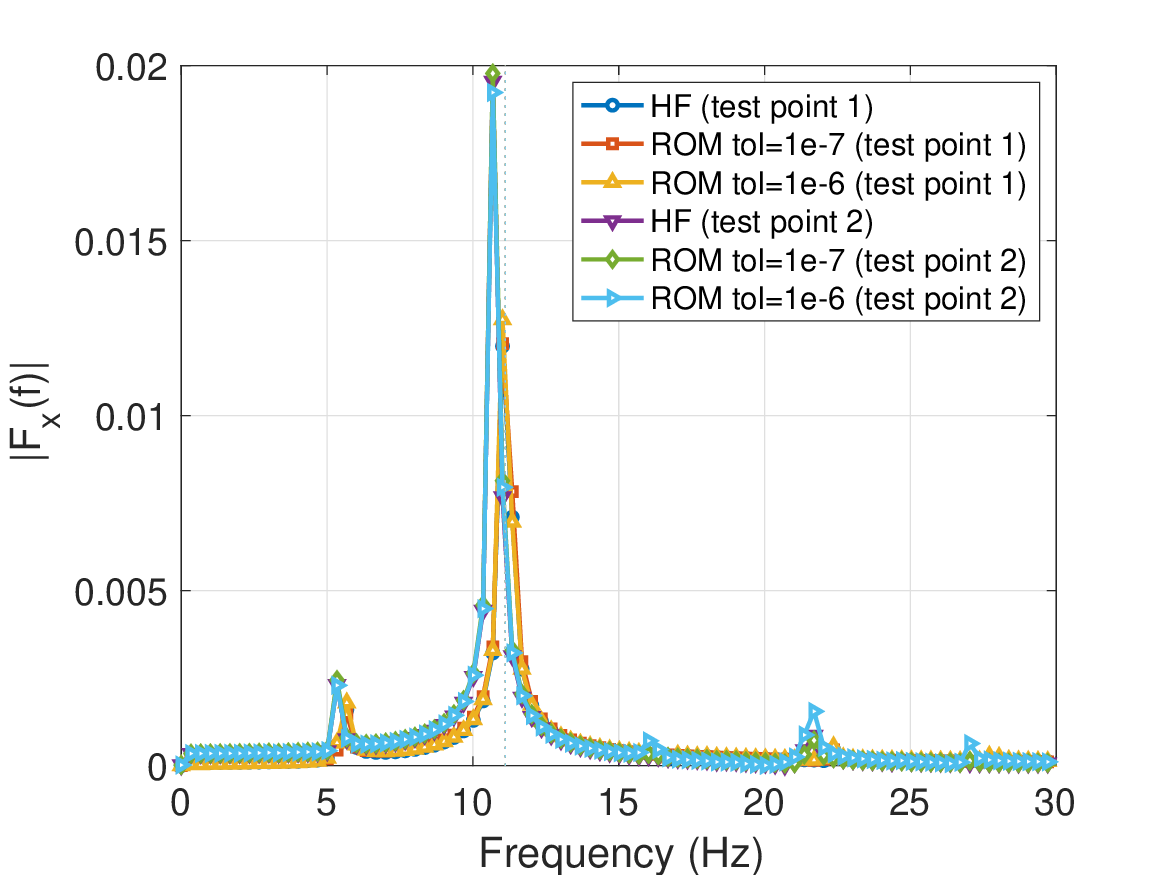}}
~~
\subfloat[$F_y$]{\includegraphics[width=0.45\textwidth]{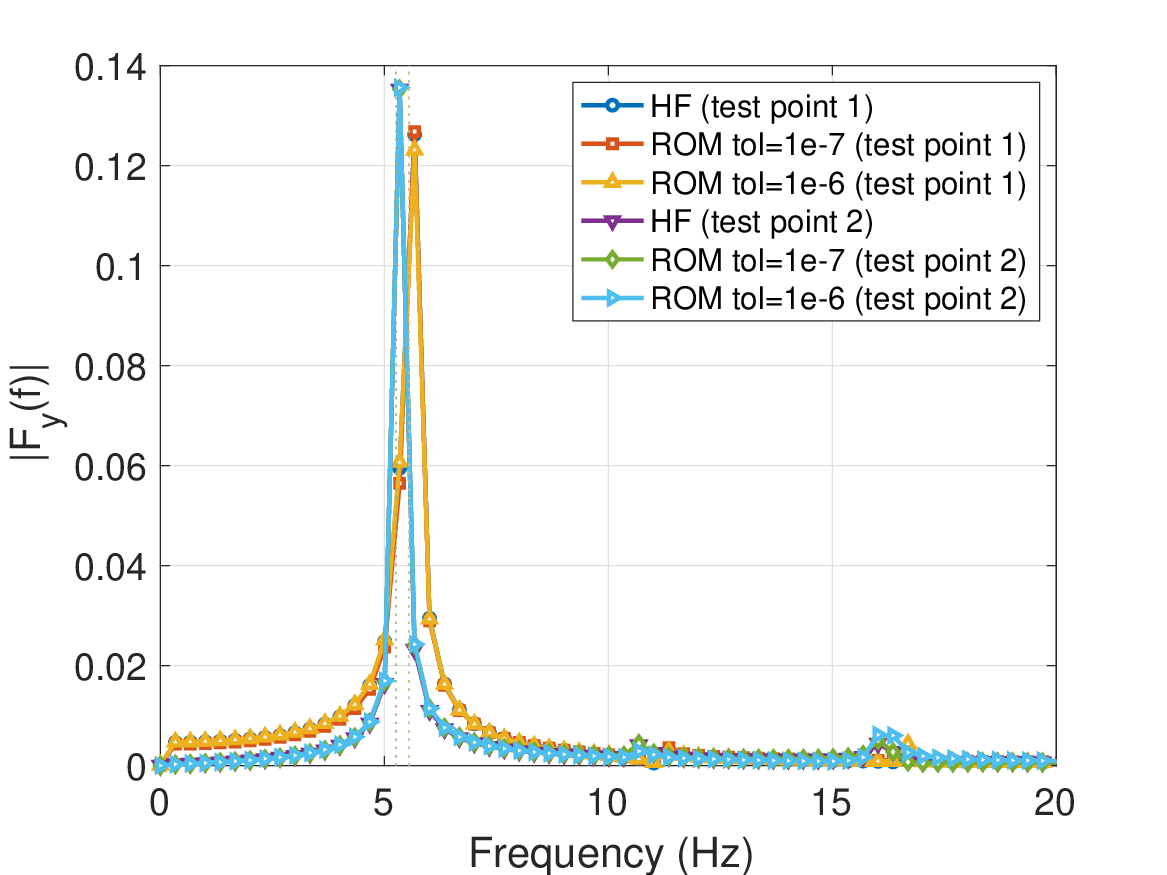}}
\caption{Turek; comparison of spectra of drag and lifted forces obtained by ROMs for different $\text{tol}_{\text{POD}}$ at two representative test points.}
\label{fig:ROM_turek_spectrum_param}
\end{figure}

\begin{figure}[H]
\centering{\includegraphics[width=0.45\textwidth]{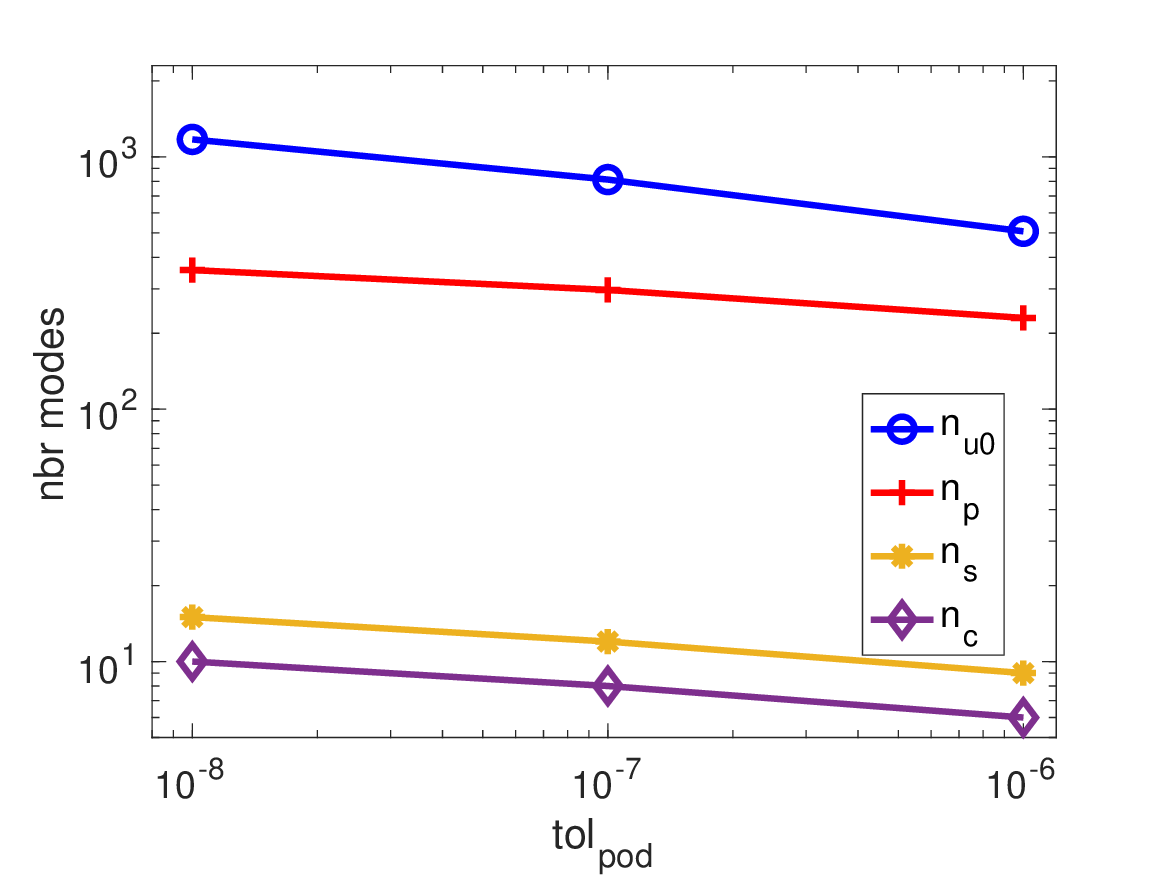}}
\caption{Turek; number of modes in terms of POD tolerance in the parametric case.}
\label{fig:turek_ROM_nb_modes_param}
\end{figure}

\textbf{Parametric case}. 
Finally, we perform a parametric study in which we vary the fluid viscosity $\mu_{\rm f} \in [0.75, 1.3]$ and the mass ratio $m^{*}:=\frac{\widetilde{\rho}_{\rm s}}{\rho_{\rm f}} \in [0.5, 2]$, while keeping $\rho_{\rm f}$ fixed. We consider 
$s=2,\Delta t=0.01\,$s for the IRK time integrator. We randomly generate $20$ training points and $5$ test points. For each training parameter, 
we collect HF snapshots over the time window   $I_{\rm train}=[33,36]$[s], which are used to construct the reduced basis. The resulting ROMs are then 
run and subsequently assessed over the extended time interval $I_{\rm test}=[33,39]$s.
Figures \ref{fig:ROM_turek_irk3_compr_param1} and \ref{fig:ROM_turek_irk3_compr_param2} 
compare the high-fidelity and ROM predictions at two representative test points. 
 In both cases, the ROM solutions remain stable, and  accurately reproduce the amplitude
 of the drag and lift forces, as well as the horizontal and vertical displacements of a probe point on the beam. 
Figure \ref{fig:ROM_turek_spectrum_param} shows the single-sided amplitude spectra of the drag ($F_x$) and lift ($F_y$) forces computed by ROMs at two representative test points. 
The two points exhibit different dominant frequencies, 
and the ROMs consistently capture the dominant frequency and the overall spectral characteristics across the considered POD tolerances.  
 These results demonstrate the robustness and predictive capability of the ROM in the parametric regime.
Figure \ref{fig:turek_ROM_nb_modes_param} reports the number of POD modes required to achieve a given tolerance in the parametric setting. We observe that the parametric dependence significantly increases the 
number of POD modes that are required to meet the target accuracy: this clearly deteriorates the efficiency of the linear-subspace ROM and ultimately motivates the development of nonlinear model reduction strategies to further enhance computational efficiency in parametric FSI problems.

\section{Conclusion}

\label{sec:conclusions}

In this paper, we introduced a 
new component-based model order reduction method 
for incompressible fluid-structure interaction (FSI) problems,  that leverages implicit Runge-Kutta (IRK) schemes and a displacement-based interface control strategy.
We rigorously verified the semi-discrete energy balance properties at the reduced-order level. Numerical experiments demonstrated that the combination of high-order Radau-IIA IRK schemes with the displacement-based interface coupling
  significantly improves the numerical stability of the reduced-order model, and ultimately  enables accurate and stable predictions of the  long-time response.
  
  For future research, we plan to extend the current framework toward fully energy-conserving numerical schemes at both the full-order and reduced-order levels. 
Second, we
 plan to combine our formulation with nonlinear approximation techniques to mitigate the issue of the slow decay of the Kolmogorov $n$-width.
Third, we plan to apply our formulation to  larger-scale systems with many fluid and solid components: towards this end, we plan to investigate localized training strategies to avoid simulations of the full system and also effective hyper-reduction techniques to speed up online calculations.

\appendix

\section{Alternative derivation of the fully discrete FSI formulation}

\label{sec:alt_derive_fully_discrete}

For completeness, we presents an equivalent derivation
of \eqref{eq:FSI_fully_discrete_stage} obtained by first applying the IRK scheme to the uncoupled fluid and solid semi-discrete equations and then assembling the coupled stage system using the mask matrices as in \eqref{eq:FSI_semidiscrete}.

The fluid subsystem is represented by the  differential-algebraic equation (DAE): 
\begin{equation}
\left\{
\begin{array}{l}
\displaystyle{
\dfrac{d\mathbf M_{\rm f}\mathbf u_{\rm f}}{dt} + \mathbf R_{\rm f}'(t,\mathbf u_{\rm f},
\mathbf p_{\rm f},
 \mathbf{P}_{\rm s,\Gamma}
\mathbf d_{\rm s}) = 0,} \\
\mathbf B_{\rm f}  (
 \mathbf{P}_{\rm s,\Gamma}
\mathbf d_{\rm s}
)  \, \mathbf u_{\rm f} = 0. 
\end{array}
\right.
\end{equation}
Here, the solid displacement $\mathbf{d}_{\rm s}$ (and ALE mesh velocity) are treated as prescribed at the IRK stages. 
Applying the IRK scheme \cite{sanderse2013energy} yields, for each stage $i=1,\ldots,s,$ 
\begin{equation}
\left\{
\begin{array}{l}
\mathbf R_{\rm f}^{(n+c_i)} = 0, \\
\mathbf B_{\rm f}^{(n+c_i)}\mathbf u_{\rm f}^{(n+c_i)} = 0,
\end{array}
\right.
\end{equation} 
where the fully discrete fluid stage residual
$$
\mathbf R_{\rm f}^{(n+c_i)}:= \mathbf M_{\rm f}^{(n+c_i)}\mathbf u_{\rm f}^{(n+c_i)} - 
\mathbf M_{\rm f}^{(n)}\mathbf u_{\rm f}^{(n)}
+\Delta t\sum_{j=1}^sa_{ij}
\mathbf R_{\rm f}'\left(t^{(n+c_j)},\mathbf u_{\rm f}^{(n+c_j)},
\mathbf p_{\rm f}^{(n+c_j)},\mathbf d_{\rm s}^{(n+c_j)}\right). 
$$
For the solid, we transform the second-order equation into a first-order system by introducing the velocity variable $\mathbf u_{\mathrm s}=\dot{\mathbf d}_{\mathrm s}$, obtaining the ODEs:
\begin{equation}
\left\{
\begin{array}{l}
\displaystyle{
\dfrac{d\mathbf M_{\rm s}\mathbf u_{\rm s}}{dt} + \mathbf R_{\rm s}'(t,\mathbf d_{\rm s}) = 0,} \\[2mm]
\dot{\mathbf d}_{\rm s}
 = \mathbf u_{\rm s}. 
\end{array}
\right.
\end{equation}
The IRK discretization reads
\begin{equation}
\left\{
\begin{array}{l}
\mathbf R_{\rm s}^{(n+c_i)}  = 0, \\
{\mathbf d}_{\rm s}^{(n+c_i)} - {\mathbf d}_{\rm s}^{(n)} - \Delta t\sum_{j=1}^sa_{ij}\mathbf u_{\rm s}^{(n+c_j)} = 0, 
\end{array}
\quad\mathrm{~for~}i=1,\ldots,s,
\right.
\label{eq:irk_solid}
\end{equation}
where the fully discrete solid stage residual 
$$
\mathbf R_{\rm s}^{(n+c_i)}:=\mathbf M_{\rm s}\mathbf u_{\rm s}^{(n+c_i)} - \mathbf M_{\rm s}\mathbf u_{\rm s}^{(n)}  +\Delta t\sum_{j=1}^sa_{ij}
\mathbf R_{\rm s}'\left(t^{(n+c_j)},\mathbf d_{\rm s}^{(n+c_j)}\right).$$
The stage velocities ${\mathbf u}_{\rm s}^{(n+c_i)}$ in the solid stage residual are then eliminated by substituting the kinematic relation
\eqref{eq:irk_solid}$_2$ into \eqref{eq:irk_solid}$_1$, resulting in a system that involves only the unknown variable ${\mathbf d}_{\rm s}^{(n+c_i)}$.

By combining fluid and solid residuals defined above, and incorporating interface coupling conditions and boundary conditions through Boolean mask matrices, we obtain the fully discrete coupled FSI system at each IRK stage \eqref{eq:FSI_fully_discrete_stage}.

\section{Additional numerical results}
\label{app:numerics}

\subsection{Vertical beam: ROM results for $s=3$, $\Delta t = 0.1$}

\label{app:vbeam_rom_s3}

\begin{figure}[H]
  \centering
\subfloat[$t=1$]{\includegraphics[width=0.33\textwidth]{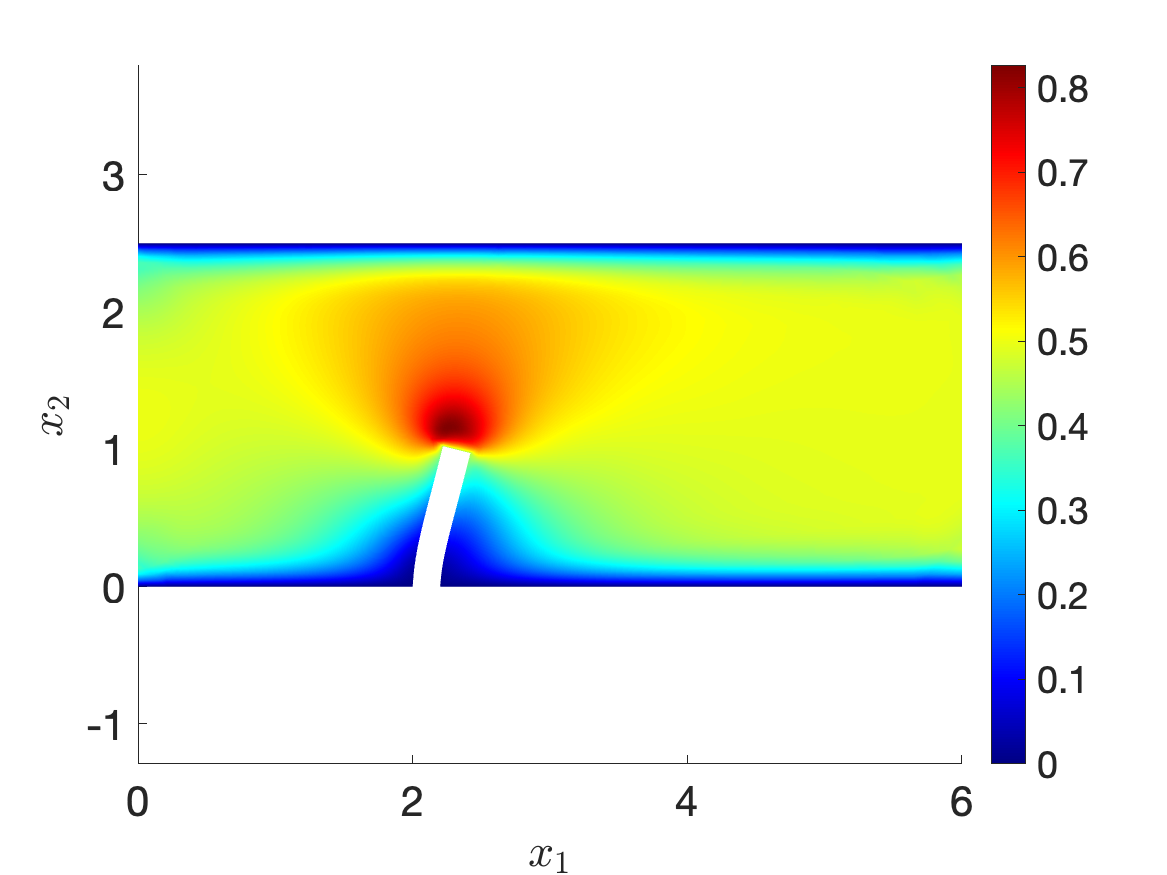}}
~~
\subfloat[$t=2$]{\includegraphics[width=0.33\textwidth]{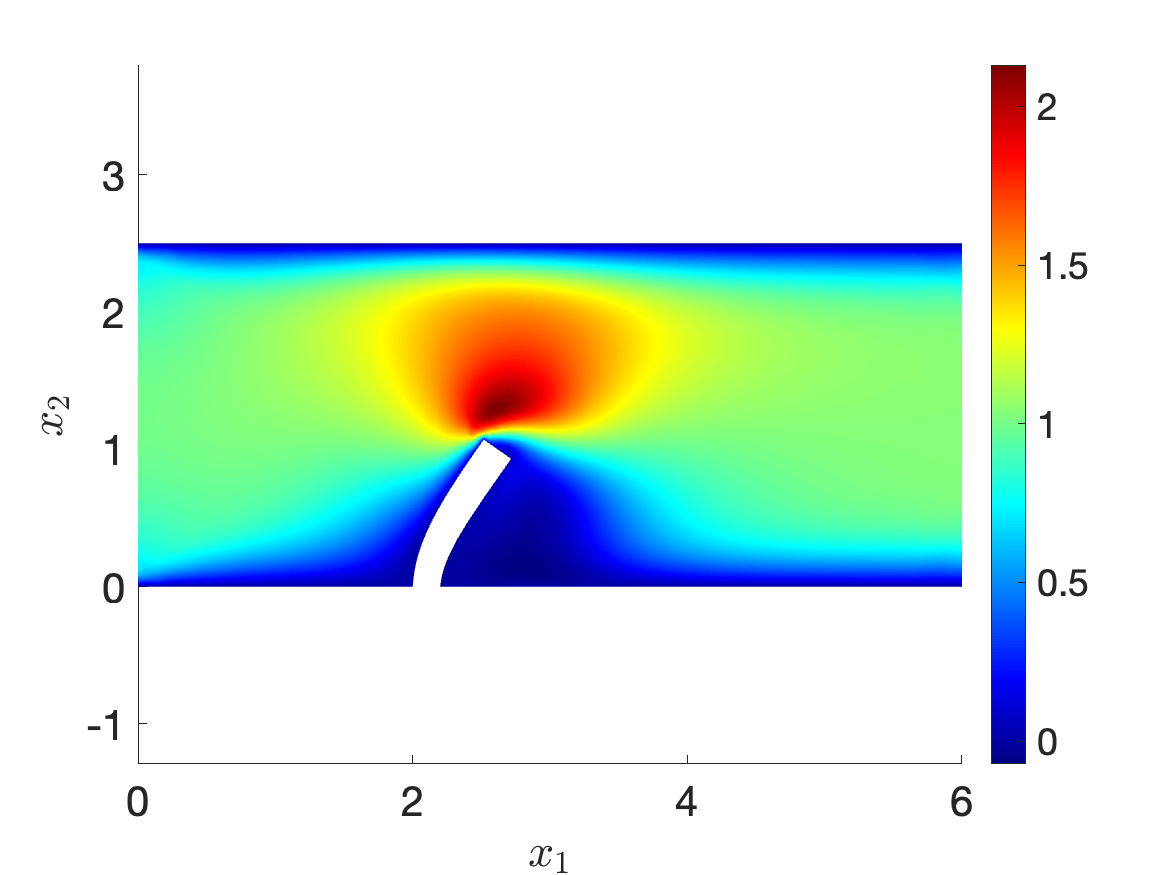}}
~~
\subfloat[$t=6$]{\includegraphics[width=0.33\textwidth]{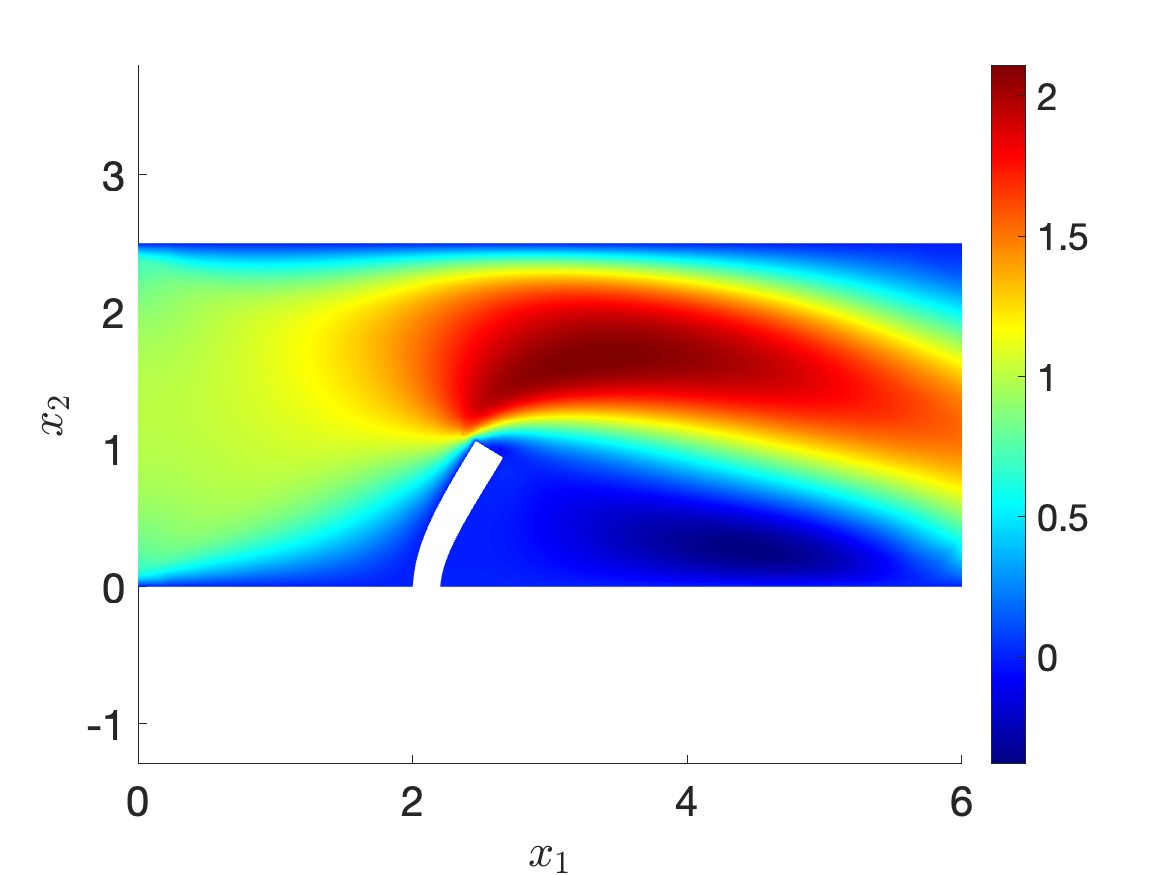}}
\\[2mm]
\subfloat[$t=1$ (error)]{\includegraphics[width=0.33\textwidth]{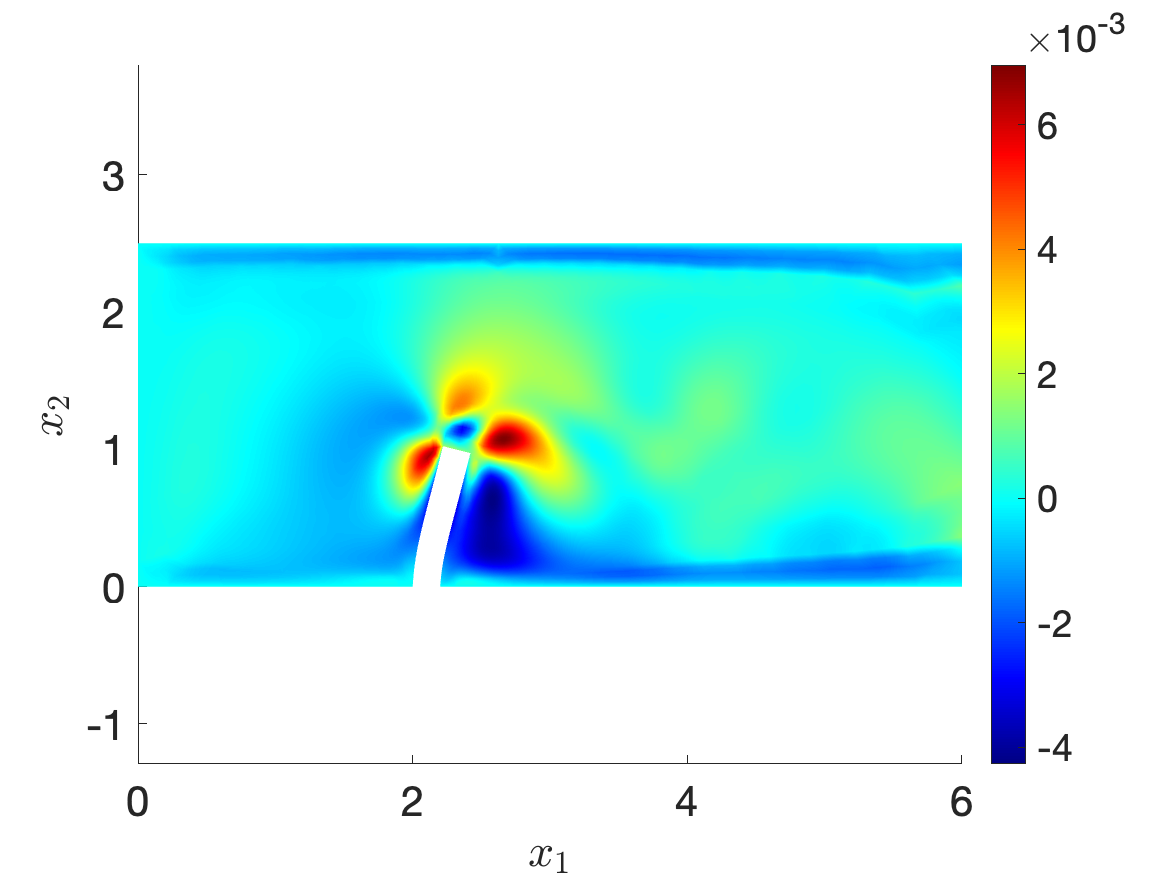}}
~~
\subfloat[$t=2$ (error)]{\includegraphics[width=0.33\textwidth]{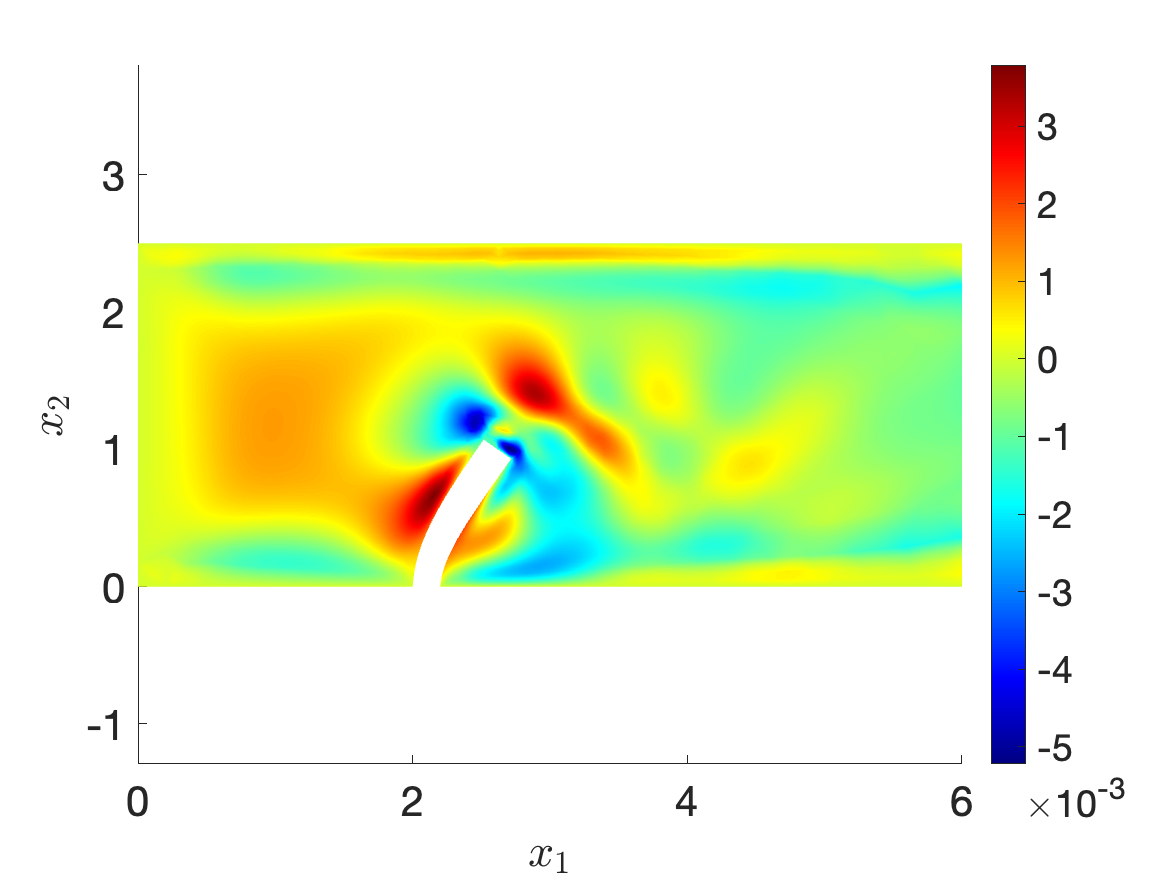}}
~~
\subfloat[$t=6$ (error)]{\includegraphics[width=0.33\textwidth]{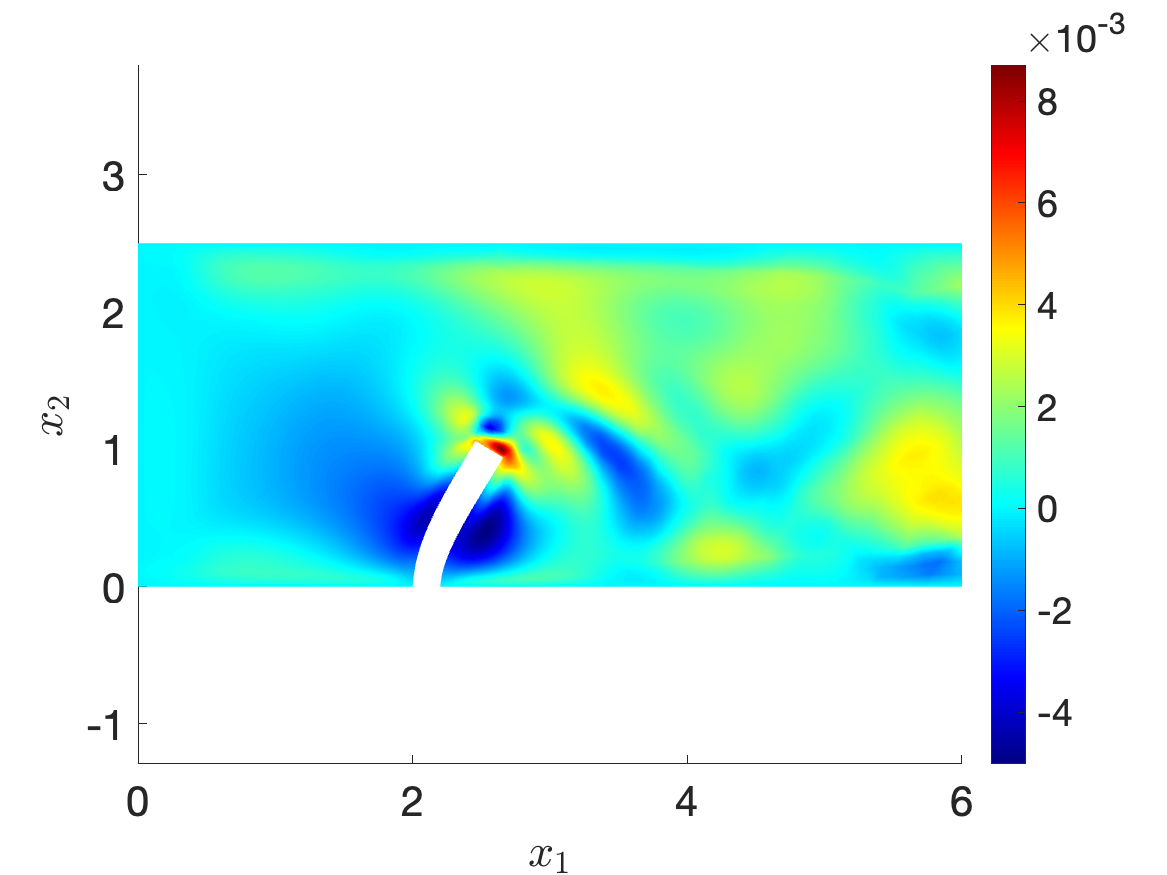}}
\caption{Vertical beam; ROM horizontal velocity and pointwise error for IRK with $s=3$, $\Delta t=0.1$\,s and $\mathrm{tol}_{\rm POD}=10^{-5}$.}
\label{fig:vbeam_ROM_irk5_stab}
\end{figure}

We report here the ROM results for IRK with $s=3$ and $\Delta t = 0.1$\,s, complementing the $s=2$, $\Delta t = 0.05$\,s case discussed in the main text. The POD basis is built from IRK snapshots with $t \in [0,6]$\,s, and the ROM is advanced over the same time interval.
Figure \ref{fig:vbeam_ROM_irk5_stab} shows the ROM horizontal velocity field and corresponding pointwise error at $t=1,2,6$\,s for $\mathrm{tol}_{\rm POD}=10^{-5}$. The ROM remains stable and provides an accurate approximation of the HF solution.

\begin{figure}[H]
\centering
\subfloat[ROM errors]{\includegraphics[width=0.45\textwidth]{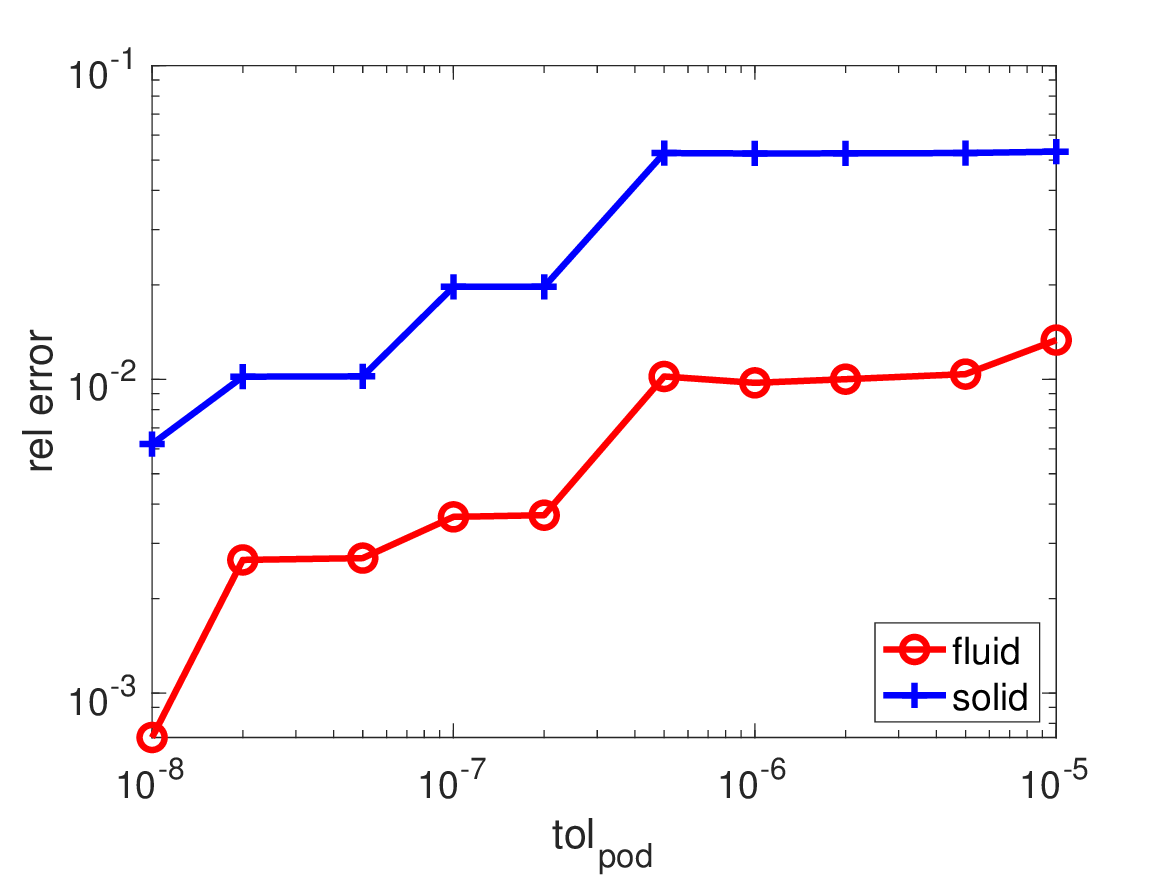}}
~~
\subfloat[Number of modes]{\includegraphics[width=0.45\textwidth]{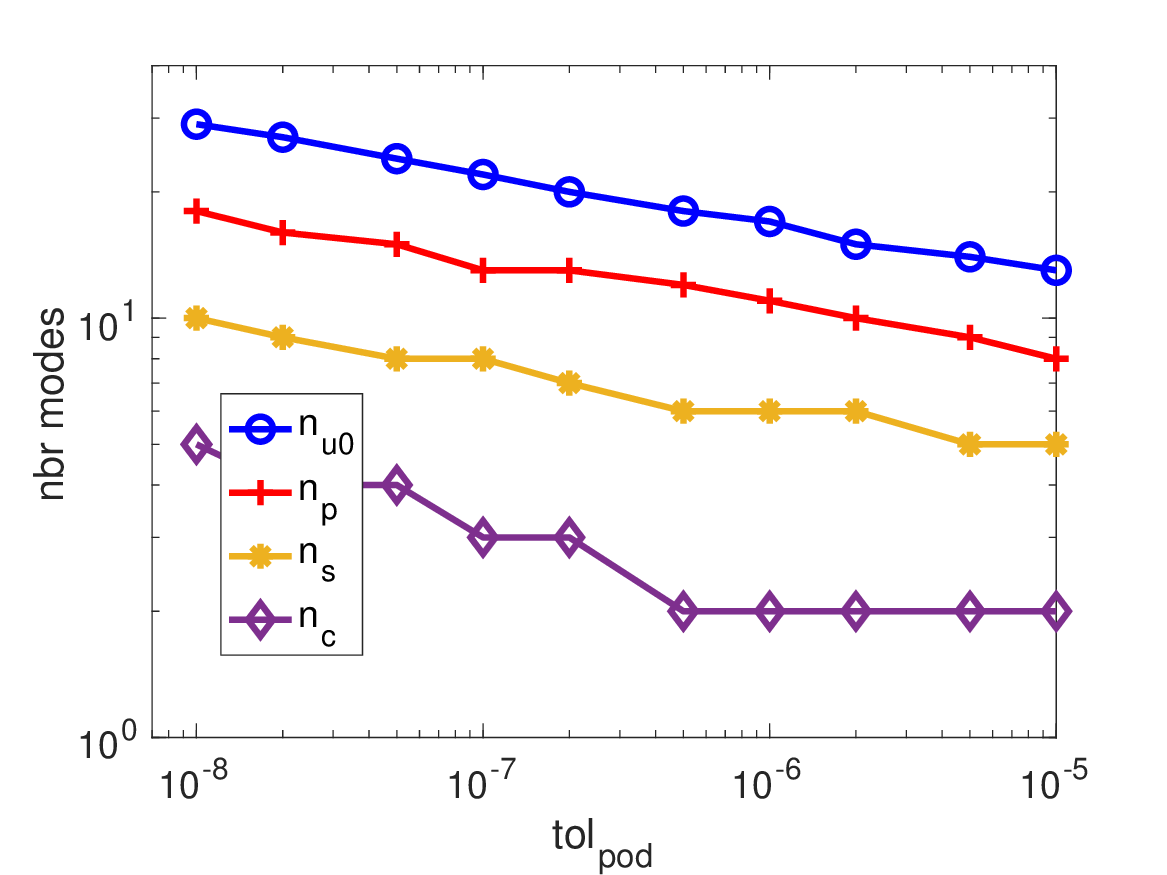}}
\caption{Vertical beam; ROM errors and number of modes versus POD tolerance ($s=3$, $\Delta t=0.1$\,s).}
\label{fig:vbeam_ROM_err_s3_0d1}
\end{figure}

Figure \ref{fig:vbeam_ROM_err_s3_0d1} shows the $H^1\times L^2$ relative error of the velocity-pressure pair and the number of modes as functions of $\mathrm{tol}_{\rm POD}$. As in the $s=2$ case, relaxing the POD tolerance reduces the number of modes but increases the approximation error.

\begin{figure}[H]
\centering
\subfloat[$F_x$]{\includegraphics[width=0.33\textwidth]{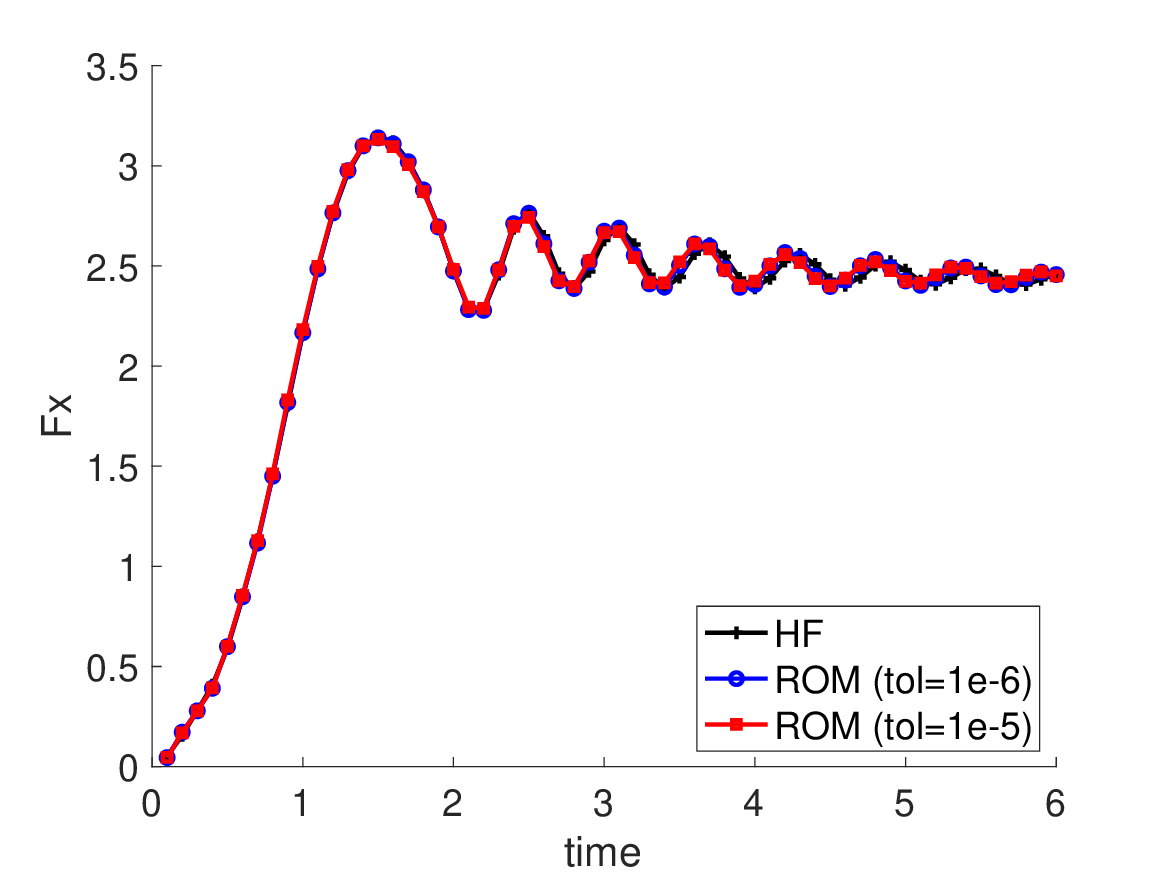}}
~~
\subfloat[$F_y$]{\includegraphics[width=0.33\textwidth]{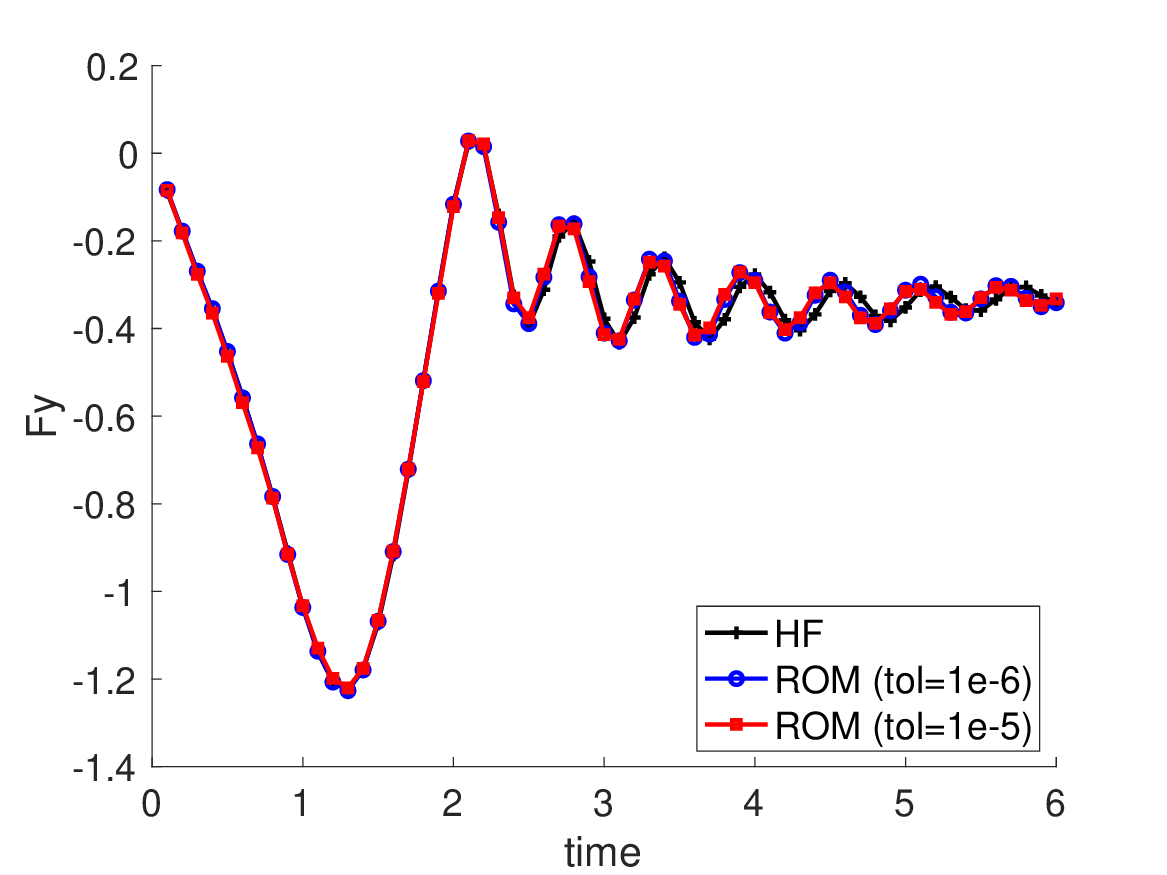}}
~~
\subfloat[Tip $x$-displacement]{\includegraphics[width=0.33\textwidth]{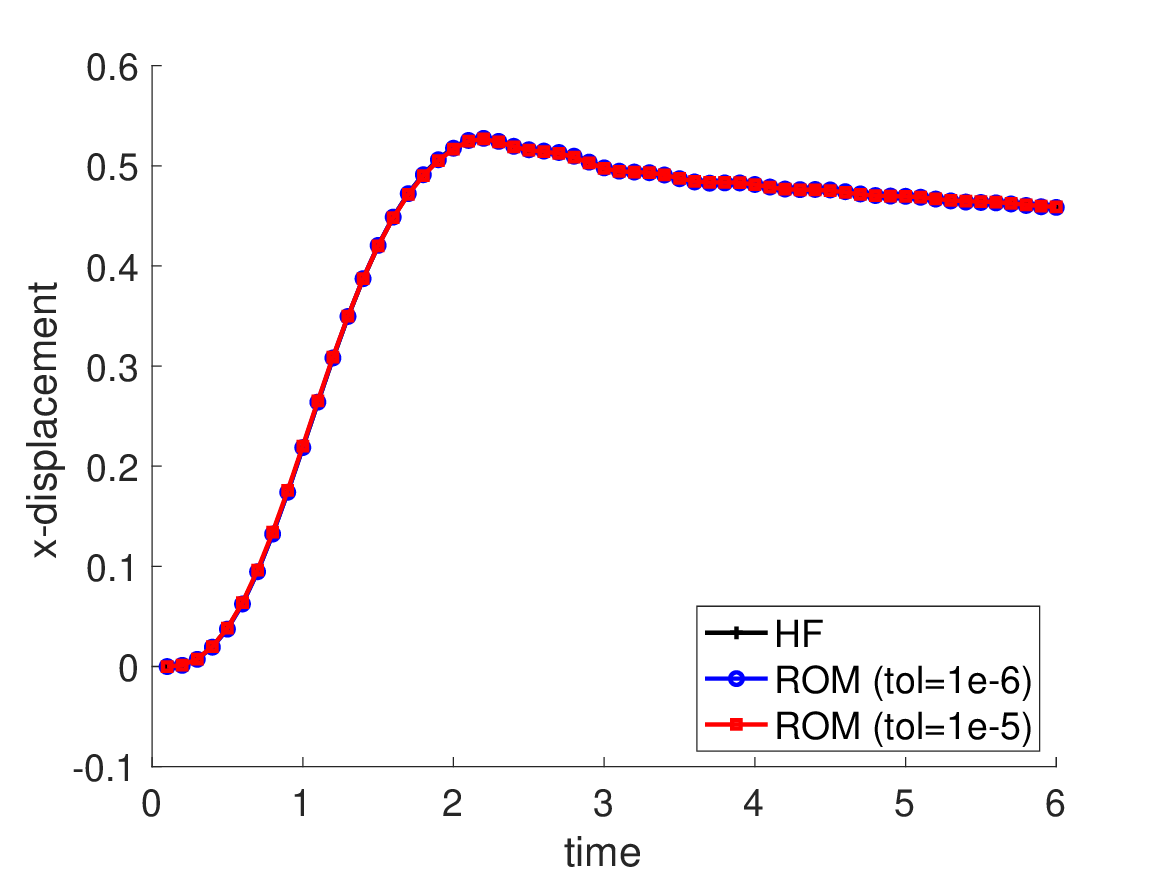}}
\caption{Vertical beam; comparison between HF and ROM predictions for $F_x$, $F_y$, and tip $x$-displacement ($s=3$, $\Delta t=0.1$\,s).}
\label{fig:vbeam_compr_s3_0d1}
\end{figure}

Figure \ref{fig:vbeam_compr_s3_0d1} compares HF and ROM predictions for drag, lift, and tip $x$-displacement for $\mathrm{tol}_{\rm POD}=10^{-6}$ and $10^{-5}$. In both cases, the ROMs reproduce the force and displacement histories with good agreement and no sign of instability, confirming the robustness of the IRK-based bubble-port ROM for this test.

\subsection{Turek: HF results}
\label{app:turek_hf}

\begin{figure}[H]
\centering
\subfloat[$F_x$]{\includegraphics[width=0.45\textwidth]{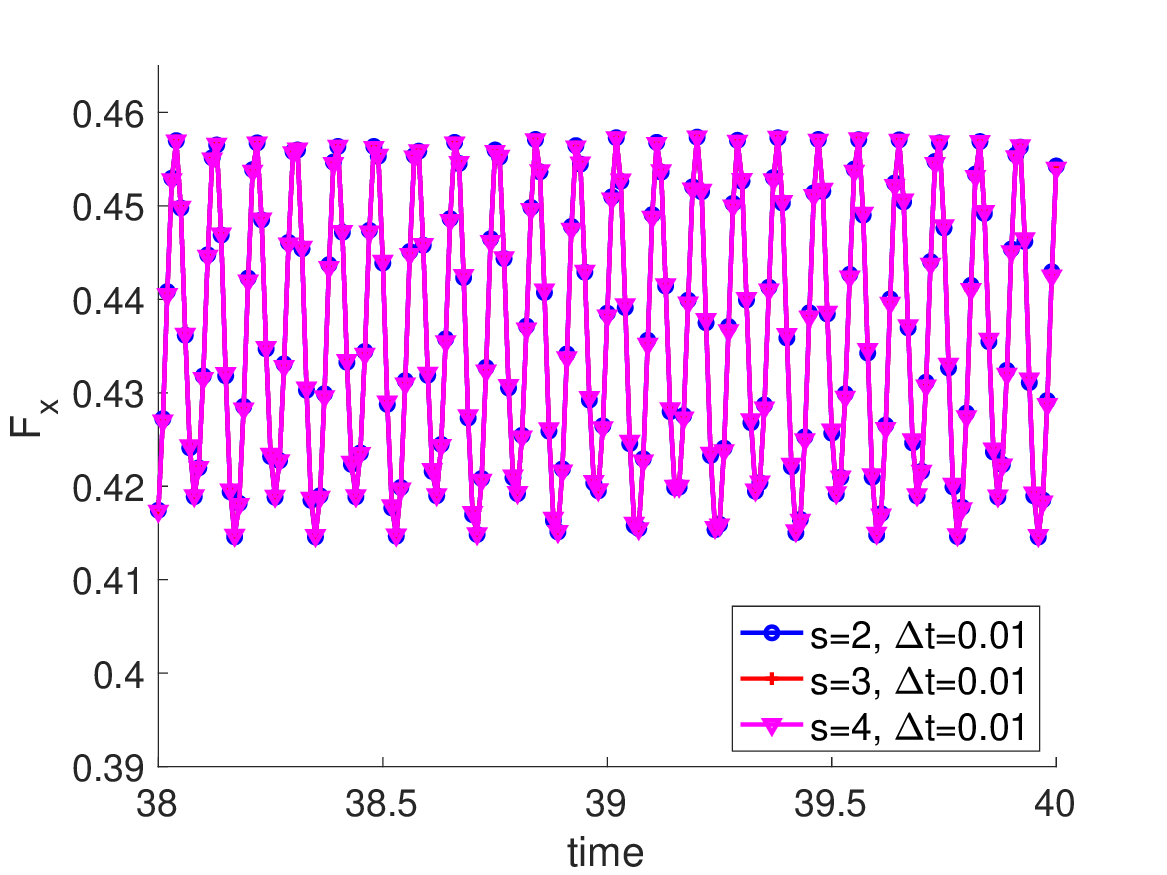}}
~~
\subfloat[$F_y$]{\includegraphics[width=0.45\textwidth]{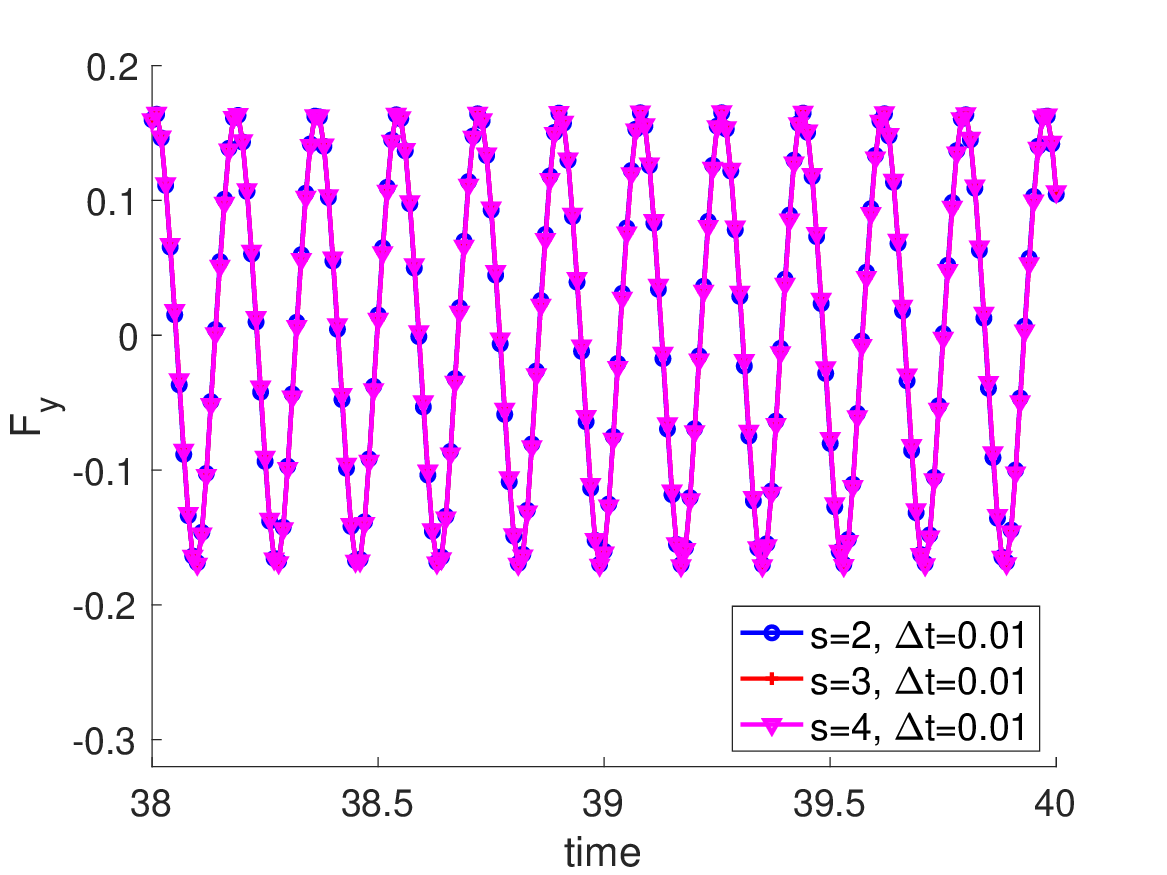}}
\caption{Turek; comparison of drag and lift using $s=2,3,4$ and $\Delta t = 0.01$\,s.}
\label{fig:hf_turek_irk357_0d01_force}
\end{figure}

\begin{figure}[H]
\centering
\subfloat[$F_x$]{\includegraphics[width=0.45\textwidth]{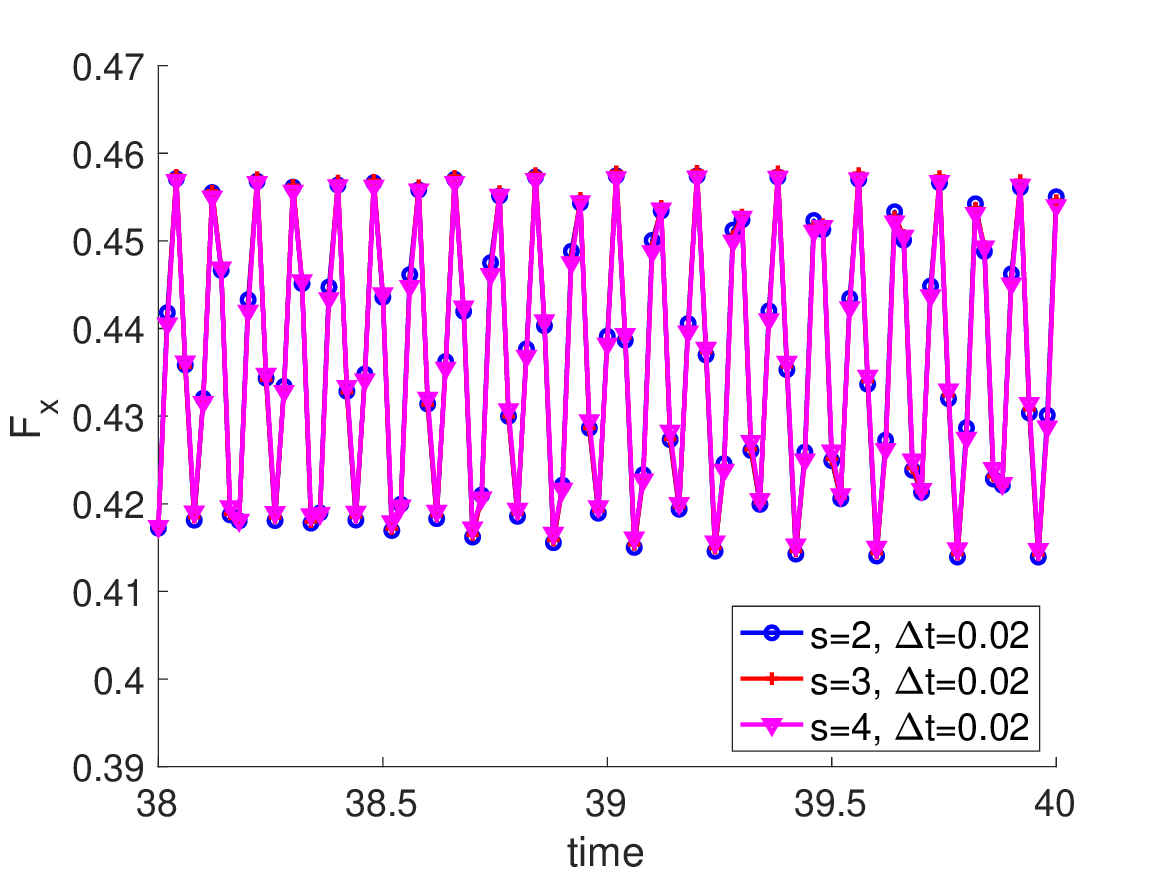}}
~~
\subfloat[$F_y$]{\includegraphics[width=0.45\textwidth]{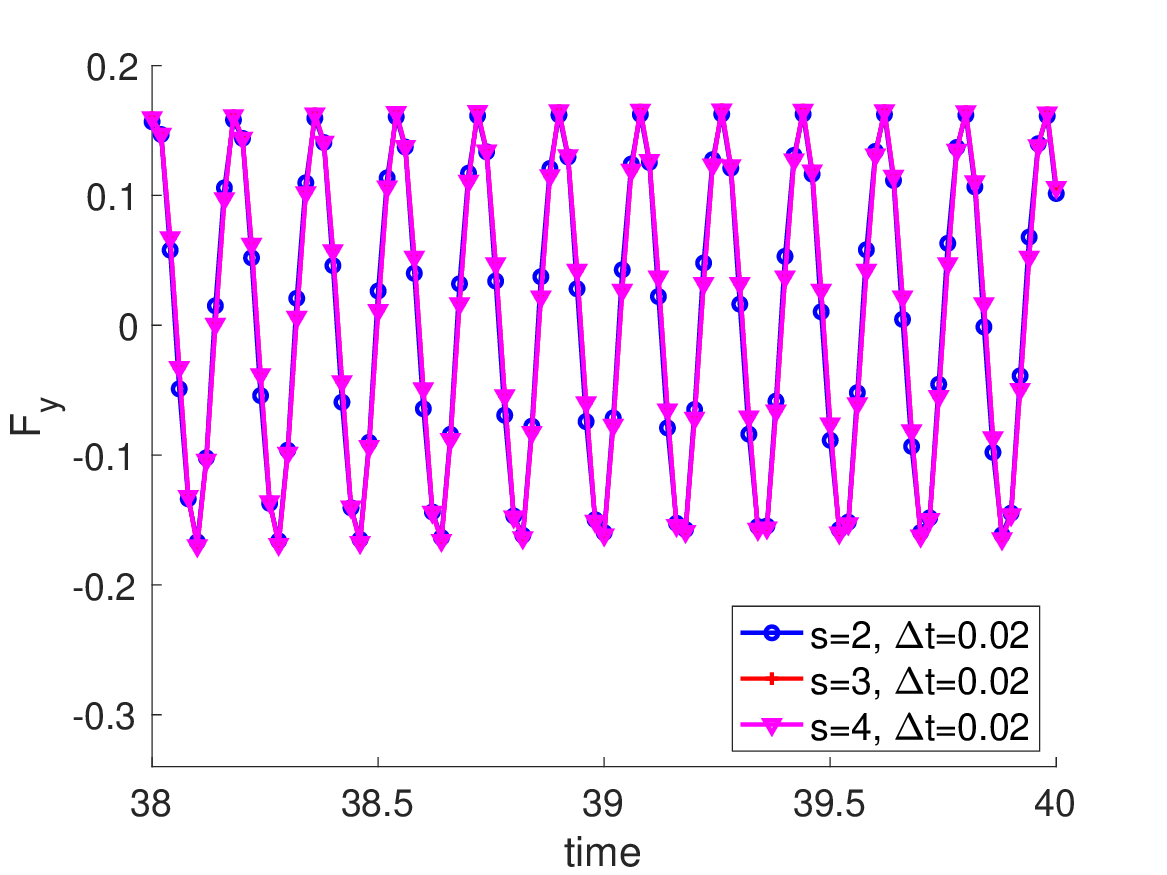}}
\caption{Turek; comparison of drag and lift using $s=2,3,4$ and $\Delta t = 0.02$\,s.}
\label{fig:hf_turek_irk357_0d02_force}
\end{figure}

We examine the influence of the Radau-IIA time integrator on the Turek FSI3 benchmark by varying the number of stages $s=2,3,4$ and the time step $\Delta t = 0.01, 0.02\,$s. The corresponding drag and lift histories are shown in Figures \ref{fig:hf_turek_irk357_0d01_force}–\ref{fig:hf_turek_irk357_0d02_force}. For both time steps, all three IRK schemes yield almost indistinguishable force signals, with negligible differences in amplitude and phase.

\subsection{Turek: ROM results for $s=3$, $\Delta t = 0.02$}
\label{app:turek_rom_s3}

\begin{figure}[H]
  \centering
\subfloat[$t=32$]{\includegraphics[width=0.33\textwidth]{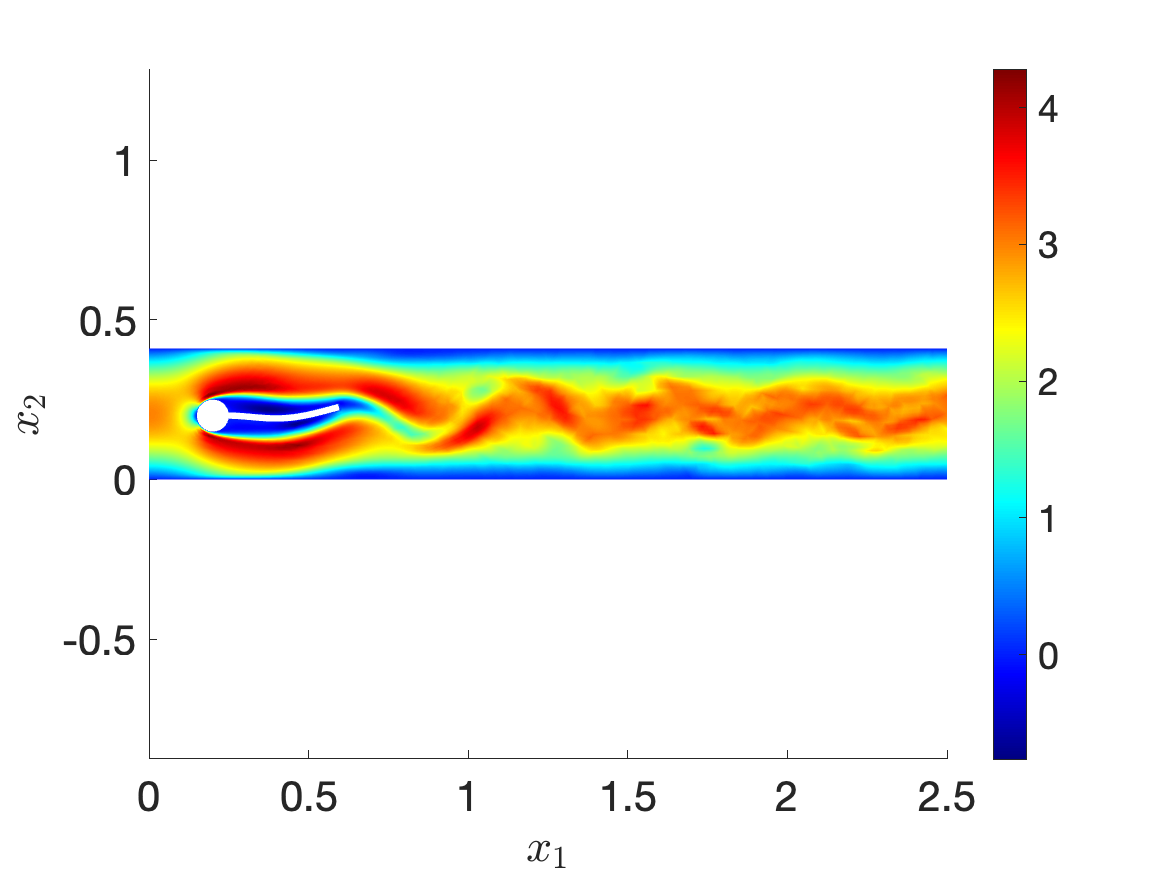}}
~~
\subfloat[$t=35$]{\includegraphics[width=0.33\textwidth]{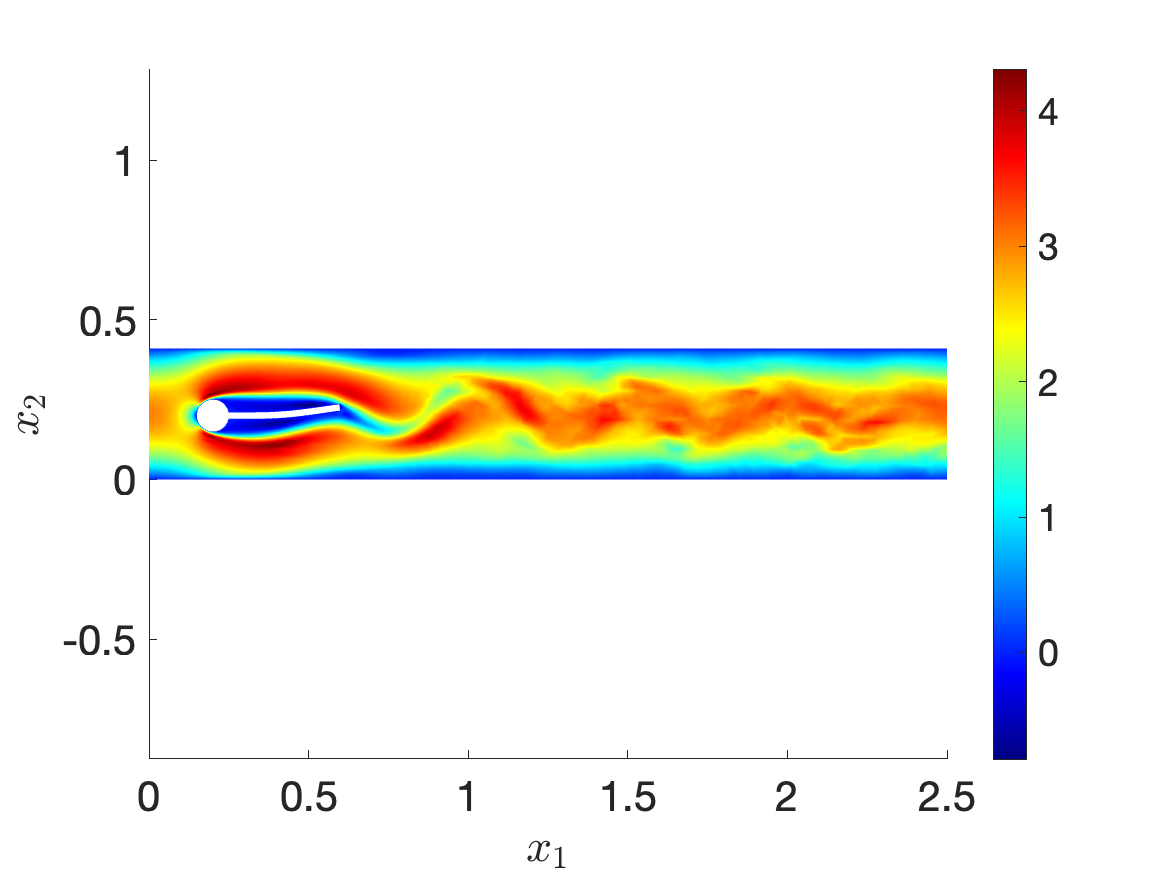}}
~~
\subfloat[$t=38$]{\includegraphics[width=0.33\textwidth]{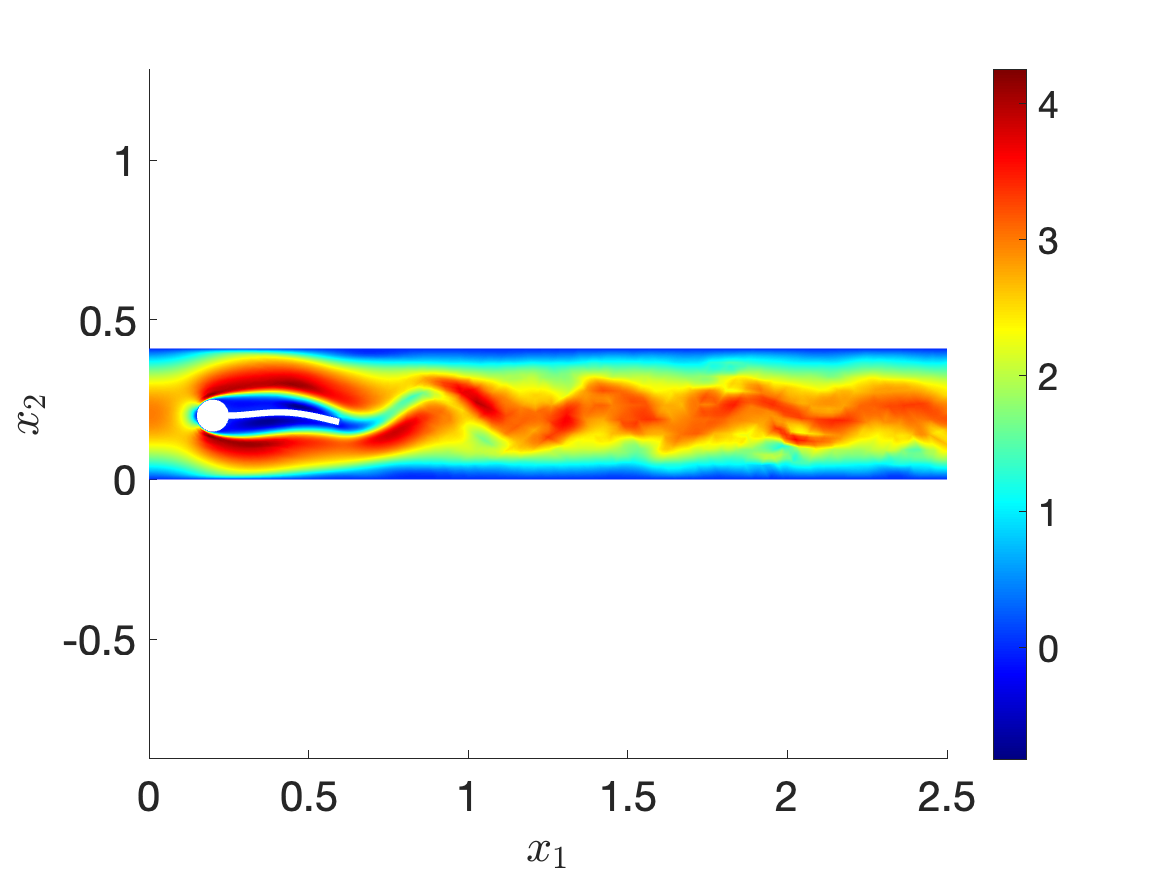}}
\\[2mm]
\subfloat[$t=32$ (error)]{\includegraphics[width=0.33\textwidth]{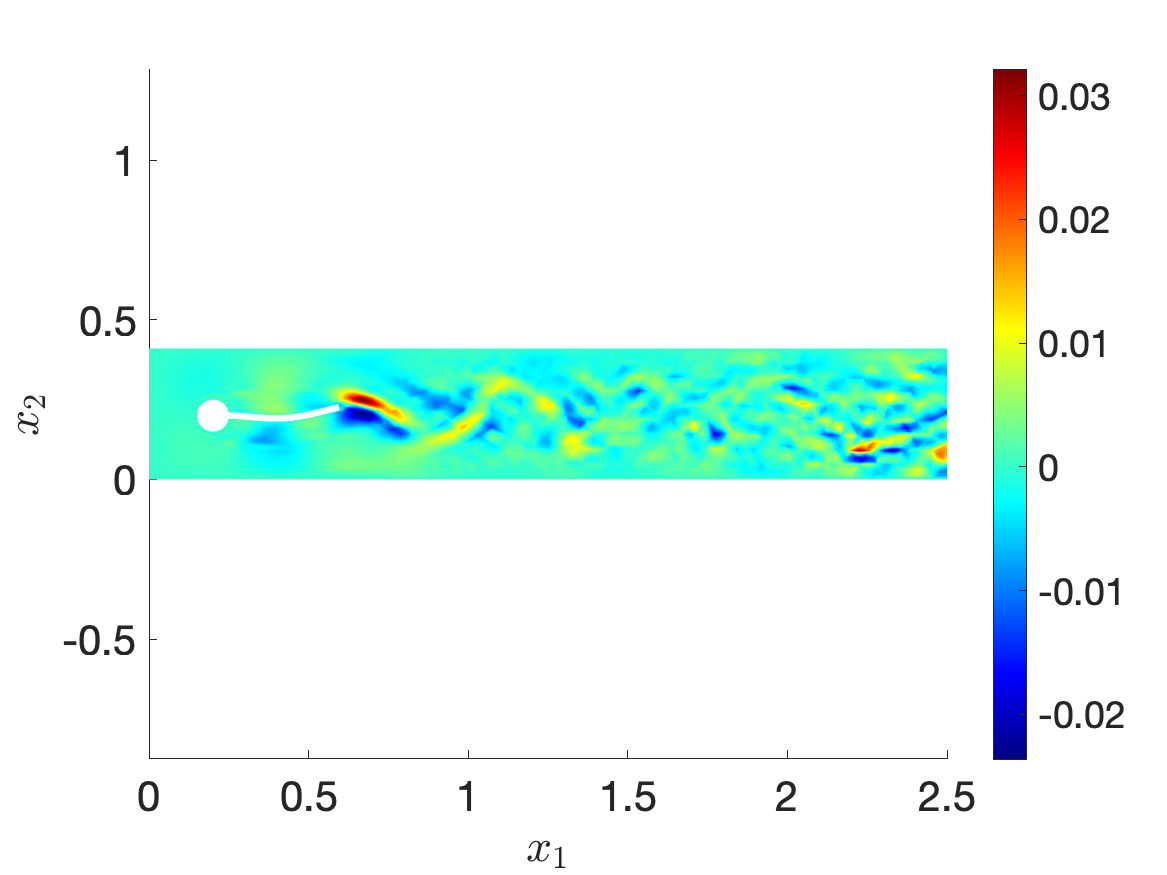}}
~~
\subfloat[$t=35$ (error)]{\includegraphics[width=0.33\textwidth]{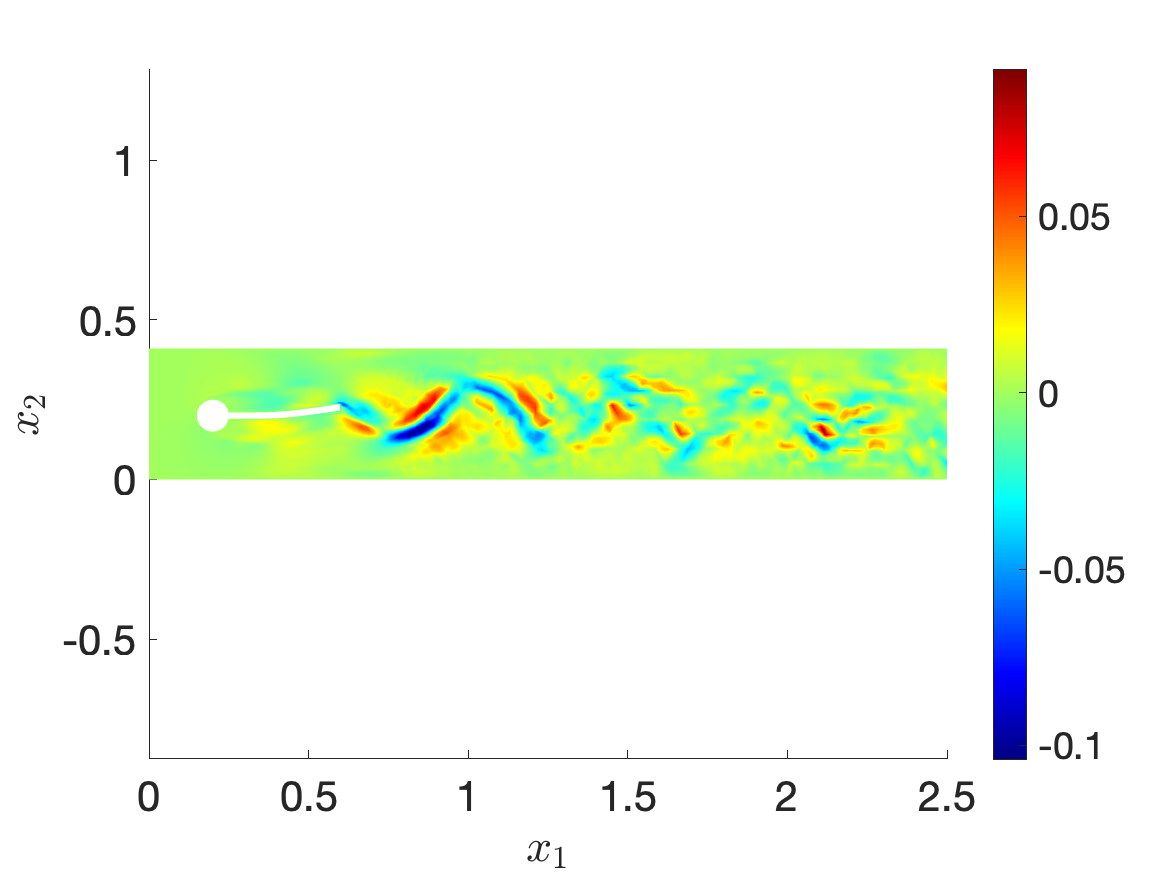}}
~~
\subfloat[$t=38$ (error)]{\includegraphics[width=0.33\textwidth]{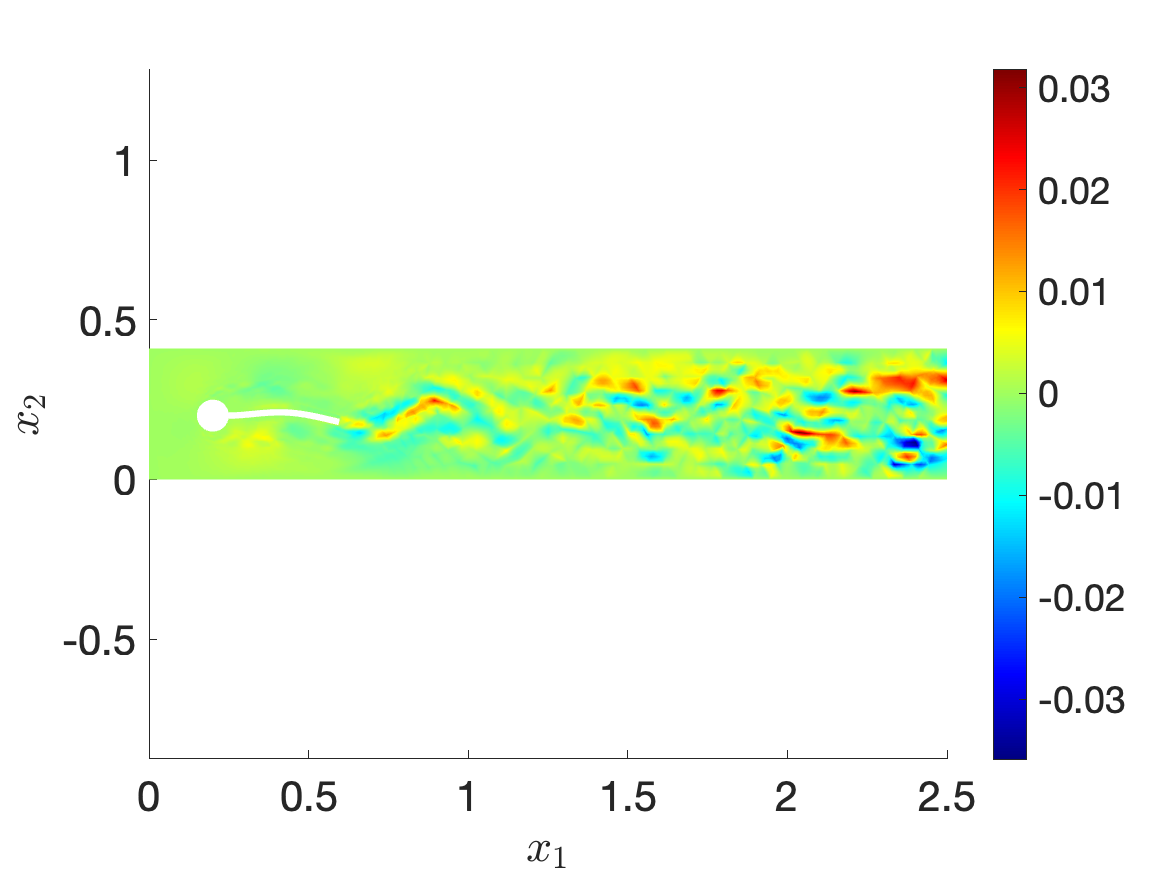}}
\caption{Turek; ROM horizontal velocity and pointwise error for IRK with $s=3$, $\Delta t=0.02$\,s and $\mathrm{tol}_{\rm POD}=10^{-6}$.}
\label{fig:turek_ROM_irk5_stab}
\end{figure}

We report here the detailed ROM results for the IRK scheme with $s=3$ and $\Delta t = 0.02$\,s, which complement the $s=2$, $\Delta t = 0.01$\,s results shown in the main text. The POD basis is again built from snapshots in the interval $t \in [30,33]$\,s, and the ROM is advanced up to $t = 39$\,s.
Figure \ref{fig:turek_ROM_irk5_stab} shows the ROM horizontal velocity field and the corresponding pointwise error at $t = 32, 35, 38$\,s for $\mathrm{tol}_{\rm POD} = 10^{-6}$. The ROM remains stable and accurately captures the flow and beam dynamics.

\begin{figure}[H]
\centering
\subfloat[ROM errors]{\includegraphics[width=0.45\textwidth]{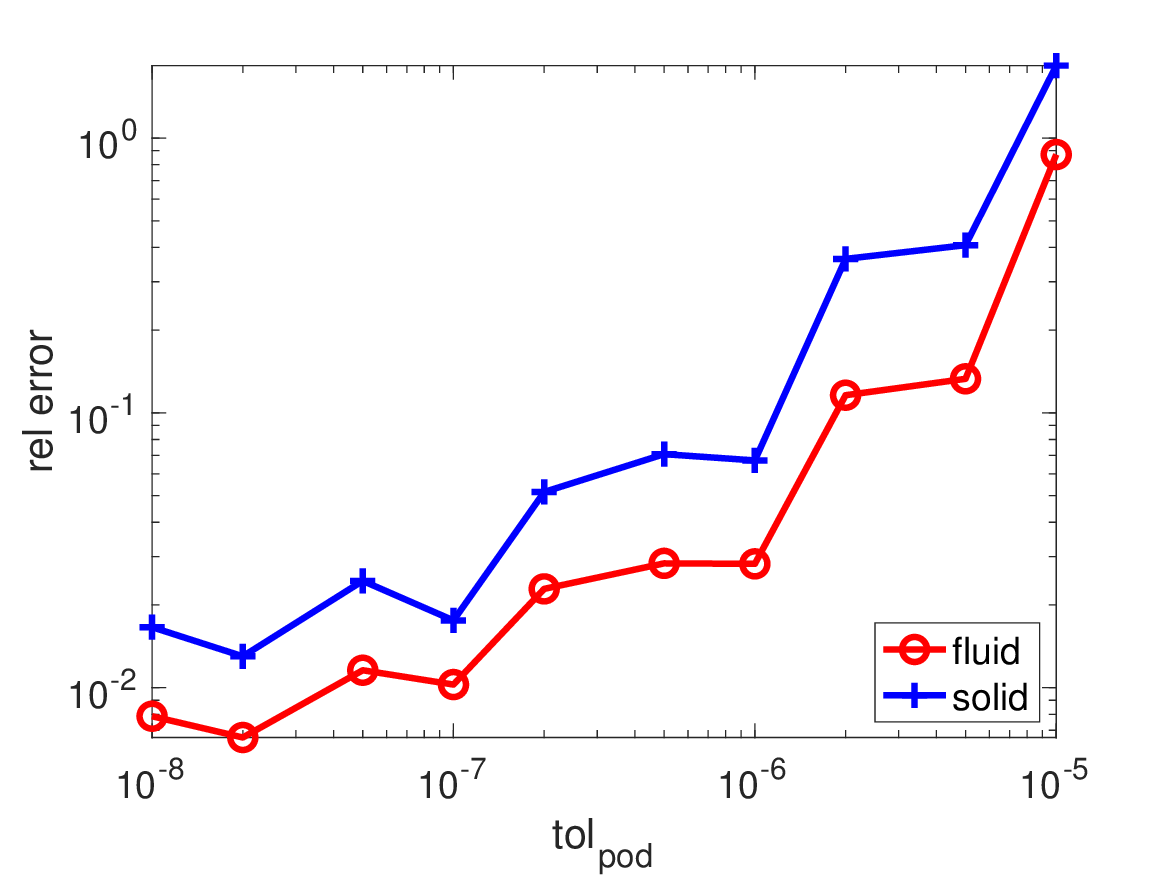}}
~~
\subfloat[Number of modes]{\includegraphics[width=0.45\textwidth]{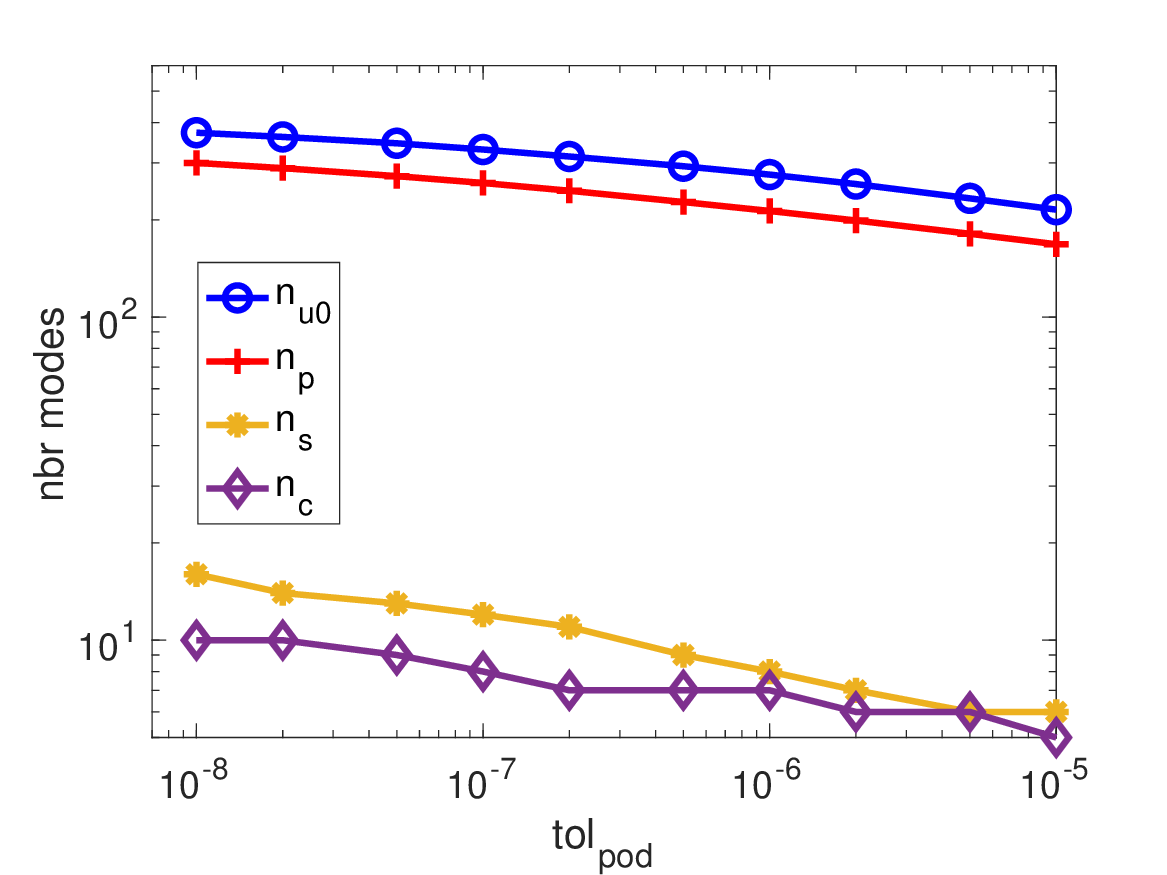}}
\caption{Turek; ROM errors and number of modes versus POD tolerance ($s=3$, $\Delta t=0.02$\,s).}
\label{fig:turek_ROM_err_s3_0d02}
\end{figure}

Figure \ref{fig:turek_ROM_err_s3_0d02} shows the $H^1 \times L^2$ relative error of the velocity-pressure pair versus POD tolerance, together with the corresponding number of reduced modes. The behavior is analogous to the $s=2$ case: smaller tolerances improve accuracy at the cost of larger reduced spaces.

\begin{figure}[H]
\centering
\subfloat[$F_x$]{\includegraphics[width=0.33\textwidth]{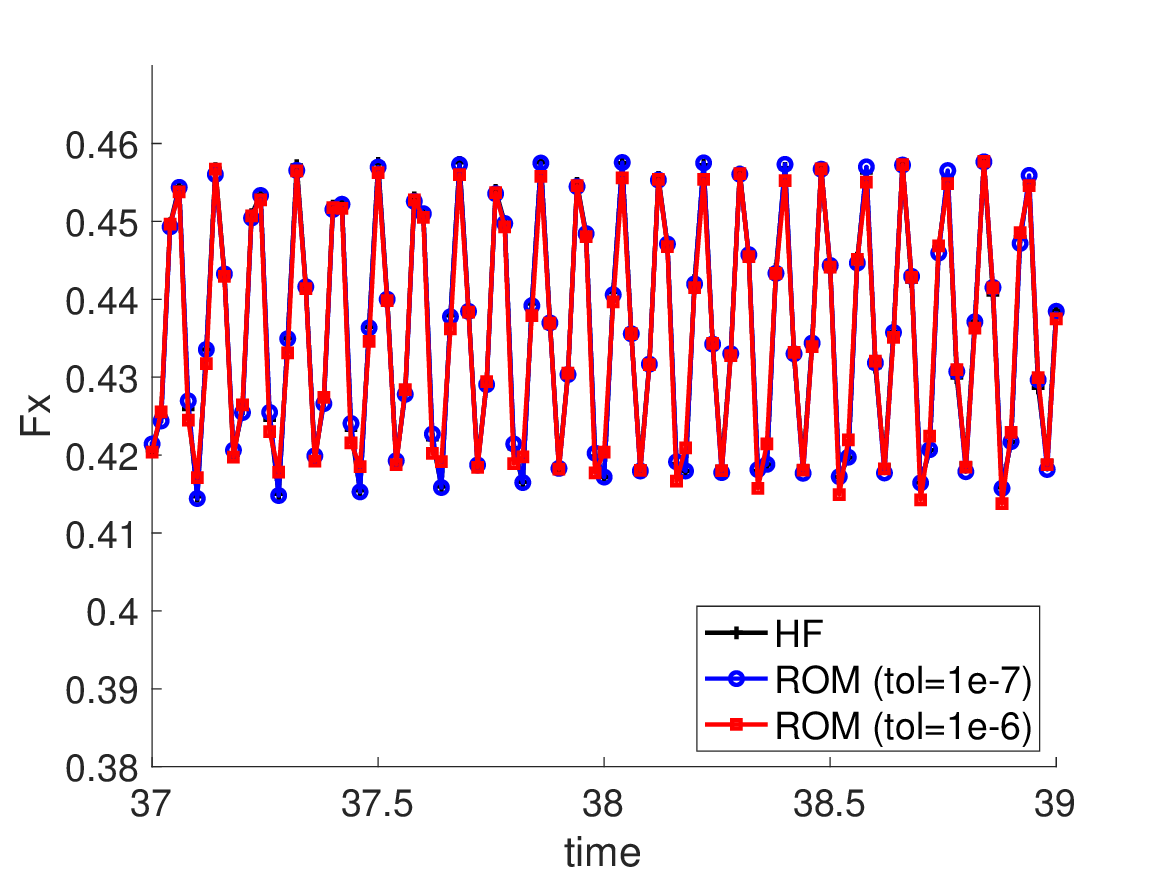}}
~~
\subfloat[$F_y$]{\includegraphics[width=0.33\textwidth]{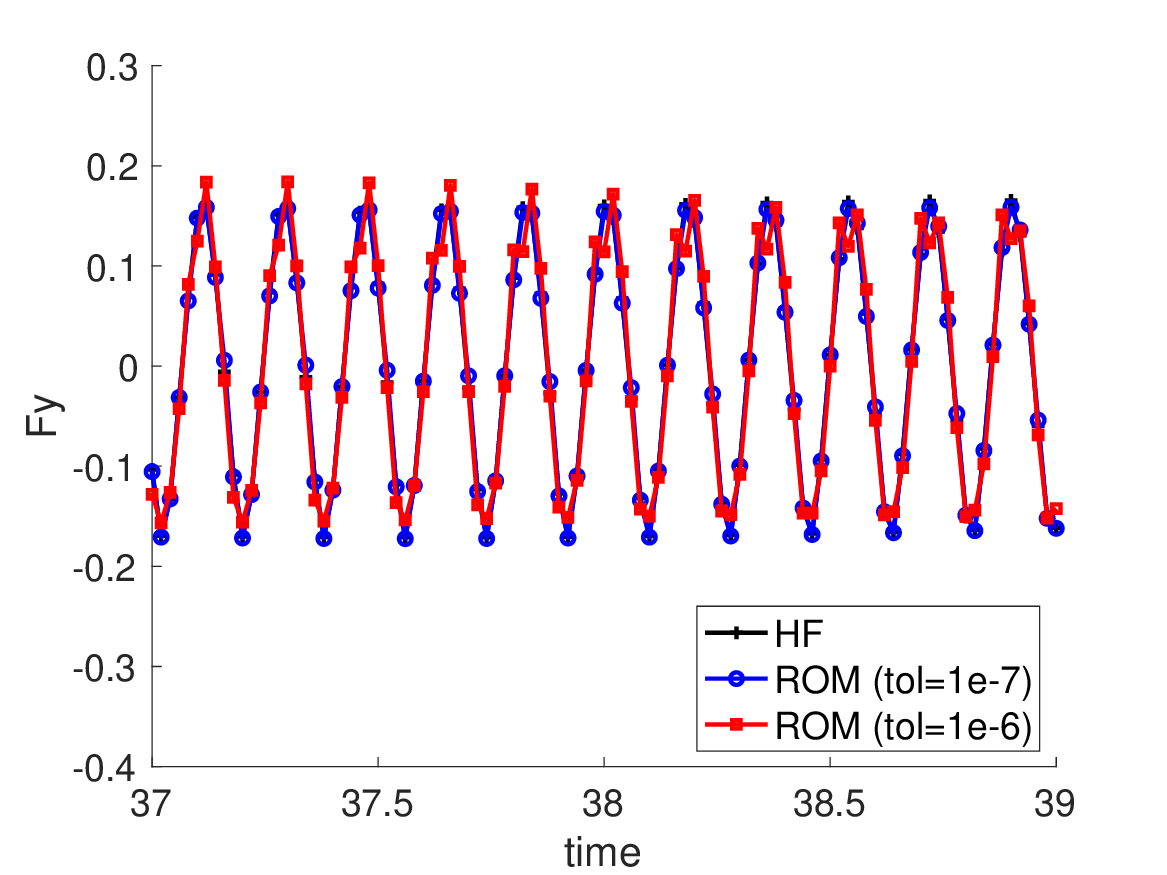}}
~~
\subfloat[Total energy]{\includegraphics[width=0.33\textwidth]{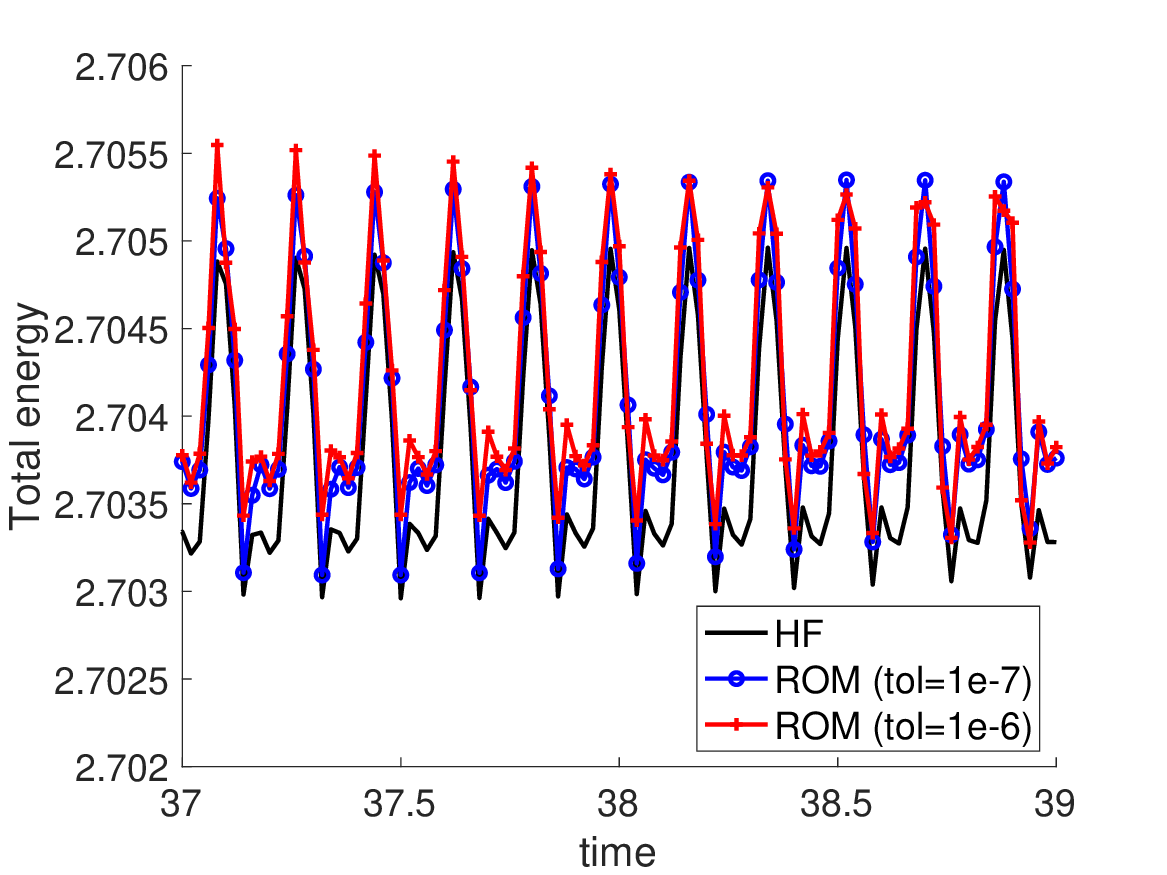}}
\caption{Turek; comparison between HF and ROM results for $F_x$, $F_y$, and total energy ($s=3$, $\Delta t=0.02$\,s).}
\label{fig:turek_compr_s3_0d02}
\end{figure}

In Figure \ref{fig:turek_compr_s3_0d02}, we compare HF and ROM predictions for drag, lift, and total energy for two ROMs with $\mathrm{tol}_{\rm POD} = 10^{-7}$ and $10^{-6}$. As in the $s=2$ case, both ROMs reproduce the force signals and energy evolution with very good agreement and no signs of instability.

\bibliography{all_references}

\end{document}